 \documentclass[times =11]{amsart}

 \usepackage[curve,all]{xy}
 \xyoption{curve}
 \xyoption{line}

 \usepackage{amscd, amsmath, amsthm, amssymb,xy, xypic, multicol}
\usepackage{amsthm}
\newtheorem{theorem}{Theorem}[section]
\newtheorem{lemma}[theorem]{Lemma}
\newtheorem{corollary}[theorem]{Corollary}
\newtheorem{proposition}[theorem]{Proposition}

\theoremstyle{definition}
\newtheorem{definition}[theorem]{Definition}
\newtheorem{example}[theorem]{Example}

\theoremstyle{remark}
\newtheorem{remark}[theorem]{Remark}

\numberwithin{equation}{section}

\newcommand{\calB}{\mathcal{B}}
\newcommand{\calC}{\mathcal{C}}
\newcommand{\calD}{\mathcal{D}}
\newcommand{\calE}{\mathcal{E}}

\newcommand{\calH}{\mathcal{H}}
\newcommand{\calI}{\mathcal{I}}

\newcommand{\calL}{\mathcal{L}}
\newcommand{\calM}{\mathcal{M}}

\newcommand{\calO}{\mathcal{O}}
\newcommand{\calP}{\mathcal{P}}

\newcommand{\calR}{\mathcal{R}}

\newcommand{\be}{\mathcal{\begin{equation}}}
\newcommand{\ee}{\mathcal{\end{equation}}}

\newcommand{\bbF}{\mathbb{F}}
\newcommand{\bbC}{\mathbb{C}}

\newcommand{\bbN}{\mathbb{N}}
\newcommand{\bbP}{\mathbb{P}}
\newcommand{\bbQ}{\mathbb{Q}}
\newcommand{\bbR}{\mathbb{R}}

\newcommand{\bbZ}{\mathbb{Z}}

\newcommand{\bfn}{\mathbf{n}}

\newcommand{\bfF}{\mathbf{F}}

\newcommand{\Or}{\textup{O}}
\newcommand{\SL}{\textup{SL}}

\newcommand{\GL}{\textup{GL}}

\newcommand{\Pic}{\textup{Pic}}
\newcommand{\Ker}{\textup{Ker}}
\newcommand{\PGL}{\textup{PGL}}
\newcommand{\PSL}{\textup{PSL}}

\newcommand{\Sp}{\textup{Sp}}

\newcommand{\Tr}{\textup{Tr}}
\newcommand{\id}{\textup{id}}

\newcommand{\trian}{\bigtriangleup}

\newcommand{\Hom}{\textup{Hom}}

\newcommand{\rank}{\textup{rank}}

\newcommand{\la}{\langle}
\newcommand{\ra}{\rangle}
\newcommand{\half}{\tfrac{1}{2}}

\newcommand{\bb}{\textup{bb}}
\newcommand{\htt}{\textup{ht}}
\newcommand{\Aut}{\textup{Aut}}
\newcommand{\Cr}{\textup{Cr}}

\newcommand{\elm}{\textup{elm}}

\newcommand{\dJJ}{\textup{dJ}}

\newcommand{\djj}{\text{dj}}
\newcommand{\PO}{\textup{PO}}

\newcommand{\pt}{${\SMALL$\bullet$}$}

\begin{document}
 \title{Finite subgroups of the plane Cremona group}

\author{Igor V. Dolgachev}

\address{Department of Mathematics, University of Michigan, 525 E. University Av., Ann Arbor, Mi, 49109}
\email{idolga@umich.edu}
\thanks{The  author was supported in part by NSF grant 0245203.}

\author{Vasily A. Iskovskikh}

\address{Steklov Institute of Mathematics, Gubkina 8, 117966, Moscow GSP-1, Russia}
\email{iskovsk@mi.ras.ru}
\thanks{The  author was supported in part by RFBR 05-01-00353-a
RFBR 08-01-00395-a , grant CRDF RUMI 2692-MO-05 and grant of
NSh 1987-2008.1}


\dedicatory{To Yuri  I. Manin}

\begin{abstract}
This paper completes the classic  and modern results on classification of conjugacy classes of finite subgroups of the  group of birational automorphisms of the complex projective plane.
\end{abstract}

\maketitle
CONTENTS
\begin{itemize}
\item[1.] Introduction
\item[2.] First examples
\item[3.] Rational $G$-surfaces
\item[4.] Automorphisms of minimal rational surfaces
\item[5.] Automorphisms of conic bundles
\item[6.] Automorphisms of Del Pezzo surfaces
\item[7.] Elementary links and the factorization theorem
\item[8.] Birational classes of minimal $G$-surfaces
\item[9.] What is left?
\item[10.] Tables
\end{itemize}

\section{Introduction}
The Cremona group $\Cr_k(n)$ over a field $k$ is the group of birational automorphisms of the projective space $\bbP_k^n$, or, equivalently, the group of $k$-automorphisms of the field $k(x_1,x_2,\ldots,x_n)$ of rational functions in $n$ independent  variables. The group $\Cr_k(1)$ is the group of automorphisms of the projective line, and hence it is isomorphic to the projective linear group $\PGL_k(2)$.  Already in the case $n=2$  the group $\Cr_k(2)$ is not well understood in spite of extensive classical literature (e.g. \cite{Coolidge}, \cite{Hudson}) on the subject as well as some modern research and expositions of classical results (e.g. \cite{Alberich}). Very little is known about the Cremona groups in higher-dimensional spaces.

In this paper we consider  the plane  Cremona group over the field of complex numbers, denoted by $\Cr(2)$. We return to the classical problem of classification of finite subgroups of $\Cr(2)$. The classification of finite subgroups of $\PGL(2)$ is well-known and goes back to F. Klein.  It consists of cyclic,  dihedral, tetrahedral, octahedral and icosahedral groups. Groups of the same type and order constitute a unique   conjugacy class in $\PGL(2)$. Our goal is to find a similar classification in the two-dimensional case.

The history of this problem begins with the work of E. Bertini
\cite{Bertini} who classified conjugacy classes of subgroups of
order 2 in $\Cr(2)$. Already in this case the answer is drastically
different. The set of conjugacy classes is parametrized by a
disconnected  algebraic variety whose connected components are
respectively isomorphic to either the moduli spaces of hyperelliptic curves of genus $g$
 (de Jonqui\`eres
involutions), or  the moduli space of canonical curves of
genus 3 (Geiser involutions), or the moduli space of  canonical curves of genus 4 with
vanishing theta characteristic (Bertini involutions). Bertini's
proof was considered to be incomplete even according to the standards of rigor of nineteenth-century algebraic geometry. A complete
and short proof was published only a few years ago by L. Bayle and
A. Beauville \cite{Bayle-Beauville}.

In 1894 G. Castelnuovo \cite{Castelnuovo}, as an application of his
theory of  adjoint linear systems, proved that any element of finite
order in $\Cr(2)$ leaves invariant either a net of lines, or a
pencil of lines, or a linear system of cubic curves with $n\le 8$
base points. A similar result was claimed earlier by S. Kantor in
his memoir which was awarded  a prize by the Accademia delle Scienze
di Napoli in 1883. However Kantor's arguments, as was pointed out by
Castelnuovo, required justifications. Kantor went much further and
announced a similar theorem for arbitrary finite subgroups of
$\Cr(2)$.  He proceeded to classify possible groups in each case
(projective linear groups, groups of de Jonqui\`eres type, and groups of
type  $M_n$).  A much clearer exposition of his results can be found
in a paper of A. Wiman \cite{Wiman}.  Unfortunately, Kantor's
classification, even with some correction made by Wiman, is
incomplete for the following two  reasons. First, only maximal groups were
considered and even some of them were missed. The most notorious
example is a cyclic group of order 8 of automorphisms of a cubic
surface, also missed by B. Segre \cite{Segre} (see \cite{Hosoh1}).
Second, although Kantor was aware of the problem of conjugacy of
subgroups, he did not attempt to fully investigate this problem.

The goal of our work is to complete Kantor's classification. We use
a modern approach to the problem initiated in the works of Yuri  Manin
and the second author (see a survey of their results in
\cite{Isk3}). Their work gives a clear understanding of the conjugacy
problem via the concept of a rational $G$-surface. It is a pair $(S,G)$ consisting of a nonsingular rational projective surface and a subgroup $G$ of its automorphism group. A birational  map
$S- \to \bbP_k^2$ realizes $G$ as a finite subgroup of
$\Cr(2)$.  Two birational isomorphic
$G$-surfaces define conjugate subgroups of $\Cr(2)$, and
conversely a conjugacy class of a finite subgroup $G$ of $\Cr(2)$
can be realized as a birational isomorphism class of $G$-surfaces.
In this way classification of  conjugacy classes of subgroups of
$\Cr(2)$ becomes equivalent to the birational classification of
$G$-surfaces. A $G$-equivariant analog of a minimal surface allows
one to concentrate on the study of minimal $G$-surfaces, i.e.
surfaces which cannot be $G$-equivariantly birationally and
regularly mapped to another $G$-surface. Minimal $G$-surfaces turn
out to be $G$-isomorphic either to the projective plane, or a conic bundle, or
a Del Pezzo surface of degree $d =9-n\le 6$ and $d = 8$. This leads to groups of
projective transformations, or groups of de Jonqui\`eres type, or
groups of type $M_n$, respectively. To complete the classification
one requires
\begin{itemize}
\item  to classify all finite groups $G$ which may occur in a minimal $G$-pair $(S,G)$;
\item to determine when two minimal $G$-surfaces are birationally isomorphic.
\end{itemize}

To solve the first part of the problem one has to compute the full
automorphism group of a conic bundle surface or a Del Pezzo surface
(in the latter case this was essentially accomplished by Kantor and
Wiman), then make a list of all finite subgroups which act minimally
on the surface (this  did not come up in the works of Kantor and
Wiman). The second part is less straightforward. For this we  use
the ideas from Mori  theory to decompose a  birational map of
rational $G$-surfaces  into elementary links. This theory was
successfully applied in the arithmetic case, where the analog of the group $G$ is  the Galois group of the base field (see \cite{Isk3}). We
borrow these results with obvious modifications adjusted to the
geometric case. Here we use the analogy between $k$-rational points
in the arithmetic case (fixed points of the Galois action) and fixed
points of the $G$-action. As an  important implication of the
classification of elementary $G$-links is the rigidity property of
groups of type $M_n$ with $n\ge 6$: any minimal Del Pezzo surface
$(S,G)$ of degree $d\le 3$ is not isomorphic to a minimal
$G$-surface of different type. This allows us to avoid much of  the painful analysis
of possible conjugacy for a lot of  groups.

The large amount of  group-theoretical computations needed for the classification of finite subgroups of groups of automorphisms of conic bundles and Del Pezzo surfaces makes us expect some possible gaps in our  classification. This seems to be a destiny of enormous classification problems. We hope that  our hard work will be useful for the future faultless classification of conjugacy classes of finite subgroups of $\Cr(2)$.

It is appropriate to mention some recent work on classification of conjugacy classes of subgroups of $Cr(2)$. We have already mentioned the work of L. Bayle and A. Beauville on groups of order 2. The papers \cite{BB},\cite{deFernex}, \cite{Zhang} study  groups of  prime orders, Beauville's paper \cite{Beauville} classifies $p$-elementary groups, and a thesis of J. Blanc \cite{Blanc} contains a classification of all finite abelian groups.  The second author studies two non-conjugate classes of subgroups isomorphic to  $S_3\times \bbZ/2\bbZ$. In the work of S. Bannai and H. Tokunaga  \cite{BT} examples are given of non-conjugate subgroups isomorphic to $S_4$ and $A_5$.

This paper is partly based on the lectures by the first author in
workshops on Cremona transformations in Torino in September 2005 and
Lisbon in May 2006.  He takes the opportunity to thank the
organizers for the invitation and for providing a good audience. We
 like to thank A. Beauville, Chenyang Xu and, especially, J. Blanc for pointing out some errors in the previous versions of our paper.

This paper is dedicated to Yuri Ivanovich Manin, to whom both authors are grateful for initiating them into algebraic geometry more than 40 years ago. Through his seminars, inspiring lectures, and as the second author's thesis adviser he was an immeasurable influence on our mathematical lives.

\section{First examples}
\subsection{Homaloidal linear systems} 

We will be working over the field of complex numbers. Recall that a dominant rational  map $\chi:\bbP^2-\to \bbP^2$ is given by a 2-dimensional linear system $\calH$ equal to the proper transform of the linear system of lines $\calH' = |\ell|$  in the target plane. A choice of a basis in $\calH$  gives an explicit formula for the map in terms of homogeneous coordinates
$$(x_0',x_1',x_2') = (P_0(x_0,x_1,x_2),P_1(x_0,x_1,x_2),P_2(x_0,x_1,x_2)),$$
where $P_0,P_1,P_2$ are linear independent homogeneous polynomials
of  degree $d$, called the (algebraic) \emph{degree} of the
map. This is the smallest number $d$ such that $\calH$ is contained
in the complete linear system $|\calO_{\bbP^2}(d)|$ of curves of degree $d$ in
the plane.  By definition of the proper transform, the linear
system $\calH$ has no  fixed components, or, equivalently, the
polynomials $P_i$'s  are mutually coprime.
The birational map $\chi$ is not a projective transformation if and
only if the degree is larger than 1, or, equivalently, when $\chi$
has \emph{base points}, the common zeros of the members of the
linear system. A linear system defining a birational map is called a
\emph{homaloidal} linear  system. Being proper transform of a
general line under a birational map, its general member is an
irreducible  rational curve. Also two general curves from the linear
system intersect outside the base points at one point. These two
conditions characterize homaloidal linear systems (more about this
later).

\subsection{Quadratic transformations} A quadratic Cremona transformation is a birational map $\chi:\bbP^2-\to \bbP^2$ of degree 2. The simplest example is the \emph{standard quadratic transformation} defined by the formula
\begin{equation}\label{sct}
\tau_1:(x_0,x_1,x_2) \mapsto  (x_1x_2, x_0x_2,x_0x_1).
\end{equation}
In affine coordinates this is given by $\tau_1: (x,y) \mapsto
(\frac{1}{x},\frac{1}{y}).$ It follows from the definition that
$\tau_1^{-1} = \tau_1$, i.e., $\tau_1$ is a birational involution of
$\bbP^2$.  The base points of $T$ are the points $p_1 = (1,0,0), p_2
= (0,1,0), p_3 = (0,0,1)$. The transformation maps an open subset of
the  coordinate line $x_i = 0$ to the  point $p_i$.  The homaloidal
linear system defining $\tau_1$ is the linear system of conics
passing through the points $p_1, p_2, p_3$.

The Moebius transformation $x\mapsto x^{-1}$ of $\bbP^1$ is conjugate to the transformation $x\mapsto -x$ (by means of the map $x\mapsto \frac{x-1}{x+1}$). This shows that the standard Cremona transformation $\tau_1$ is conjugate in $\Cr(2)$ to a projective transformation given  by
$$(x_0,x_1,x_2) \mapsto (x_0,-x_1,-x_2).$$

When  we change the homaloidal linear system defining $\tau_1$ to the
homaloidal linear system of conics passing through the point $p_2,p_3$
and tangent at $p_3$  to the line $x_0 = 0$ we obtain the
transformation
\begin{equation}
\label{st2} \tau_2:(x_0,x_1,x_2) \mapsto (x_1^2,x_0x_1,x_0x_2).
\end{equation}
In affine coordinates it is given by $(x,y) \mapsto (\frac{1}{x}, \frac{y}{x^2})$.
The transformation $\tau_2$ is also a birational  involution conjugate to a projective involution. To see this we define a rational map $\chi:\bbP^2-\to \bbP^3$ by the formula
$(x_0,x_1,x_2)\mapsto  (x_1^2,x_0x_1,x_0x_2,x_1x_2)$. The Cremona transformation $\tau_2$ acts on $\bbP^3$ via this transformation by $(u_0,u_1,u_2,u_3)\mapsto (u_1,u_0,u_3,u_2)$.  Composing with the projection of the image from the fixed point $(1,1,1,1)$ we get a birational map
$(x_0,x_1,x_2)\mapsto (y_0,y_1,y_2) = (x_1(x_0-x_1),x_0x_2-x_1^2,x_1(x_2-x_1))$. It defines the  conjugation of $\tau_2$ with the projective transformation $(y_0,y_1,y_2)\mapsto (-y_0,y_2-y_0,y_1-y_0)$.

Finally, we could further ``degenerate'' $\tau_1$  by considering
the linear system of conics passing through the point $p_3$ and
intersecting each other at this point with multiplicity 3. This
linear system defines a birational involution
\begin{equation}
\label{st3} \tau_3:(x_0,x_1,x_2) \mapsto
(x_0^2,x_0x_1,x_1^2-x_0x_2).
\end{equation}
Again it can be shown that $\tau_3$ is conjugate to a projective involution.

Recall that a birational transformation is not determined by the
choice of a homaloidal linear system, one has  to choose
additionally a basis of the linear system. In the above examples,
the basis is chosen to make the transformation an involution.

\subsection{De Jonqui\`eres  involutions}\label{dejon}
Here we exhibit a series of birational involutions which are not conjugate to each other and not conjugate to a projective involution. In affine coordinates they are given by the formula
\begin{equation}\label{dj}
\djj_{P}: (x,y) \mapsto (x,\frac{P(x)}{y}),
\end{equation}
where $P(x)$ a polynomial of degree $2g+1$  or $2g+2$ without
multiple roots. The conjugation by the transformation $(x,y)
\mapsto (\frac{ax+b}{cx+d},y)$ shows that the conjugacy class of
$\djj_P$ depends only on the orbit of the set of roots of $P$ with
respect to the group $\PGL(2)$, or, in other words, on the birational
class of the hyperelliptic curve
\begin{equation}\label{hc}
y^2+P(x) = 0.
\end{equation}

The transformation $\djj_P$ has the following beautiful geometric
interpretation. Consider the projective model $H_{g+2}$ of the
hyperelliptic curve \eqref{hc} given by the homogeneous  equation of
degree $g+2$
\begin{equation}\label{hyp}
T_2^2F_g(T_0,T_1)+2T_2F_{g+1}(T_0,T_1)+F_{g+2}(T_0,T_1) = 0,
\end{equation}
where
$$D = F_{g+1}^2-F_gF_{g+2} = T_0^{2g+2}P(T_1/T_0)$$
is the homogenization of the polynomial $P(x)$. The curve $H_{g+2}$
has an ordinary singular point of multiplicity $g$ at $q = (0,0,1)$
and the projection from this point to $\bbP^1$ exhibits the curve as
a double cover of $\bbP^1$ branched over the $2g+2$ zeroes of  $D$.

Consider the affine set $T_2 = 1$ with affine coordinates $(x,y) =
(T_0/T_2, T_1/T_2)$. A general line $y = kx$ intersects the curve
$H_{g+2}$ at the point $q = (0,0)$ with multiplicity $g$ and at two
other points $(\alpha,k\alpha)$ and $(\alpha',k\alpha'),$ where
$\alpha,\alpha'$ are the roots of the quadratic equation
$$t^2F_{g+2}(1,k)+2tF_{g+1}(1,k)+F_{g}(1,k) = 0.$$
Take a general point $p = (x,kx)$ on the line and define the point
$p' = (x',kx')$ such that the pairs
$(\alpha,k\alpha),(\alpha',k\alpha')$ and $(x,kx), (x',kx')$ are
harmonic conjugate. This means that $x,x'$ are the roots of the
equation $at^2+ 2bt+c= 0,$ where
$aF_g(1,k)+cF_{g+2}(1,k)-2bF_{g+1}(1,k) = 0$. Since $x+x' = -2b/a,
xx' = c/a$ we get $F_g(1,k)+xx'F_{g+2}(1,k)+(x+x')F_{g+1}(1,k) = 0$.
We express $x'$ as $(ax+b)/(cx+d)$ and solve for $(a,b,c,d)$ to
obtain
$$x' = \frac{-F_{g+1}(1,k)x-F_g(1,k)}{xF_{g+2}(1,k)+ F_{g+1}(1,k)}.$$
Since $k = y/x$, after changing the affine coordinates $(x,y) = (x_0/x_2,x_1/x_2)$ to
$(X,Y) = (x_1/x_0,x_2/x_0) = (y/x,1/x)$, we get
\begin{equation}\label{dj2}
IH_{g+2}:  (X,Y) \mapsto (X',Y') := \bigl(X, \frac{-YP_{g+1}(X)-P_{g+2}(X)}{P_g(X)Y+P_{g+1}(X)}\bigr),
\end{equation}
where $P_i(X) = F_i(1,X)$. Let $T:(x,y) \mapsto (x,yP_g+P_{g+1})$.
Taking $P(x) = P_{g+1}^2-P_gP_{g+2}$, we check that
$T^{-1}\circ \djj_P\circ T = IH_{g+2}$.
 This shows that our geometric de Jonqui\`eres involution $IH_{g+2}$ given
 by \eqref{dj2} is conjugate to the de Jonqui\`eres involution $\djj_P$ defined by \eqref{dj}.

Let us rewrite  \eqref{dj2} in homogeneous coordinates:
\begin{eqnarray}\label{f1}
x_0' & = & x_0(x_2F_{g}(x_0,x_1)+F_{g+1}(x_0,x_1)), \\\notag
x_1' & = & x_1(x_2F_{g}(x_0,x_1)+F_{g+1}(x_0,x_1)),\\\notag
x_2'&=& -x_2F_{g+1}(x_0,x_1)-F_{g+2}(x_0,x_1).\notag
\end{eqnarray}
Now it is clear that the homaloidal linear system defining
$IH_{g+2}$ consists of curves of degree $g+2$ which pass through the
singular point $q$ of the hyperelliptic curve \eqref{hyp} with
multiplicity $g$. Other base points satisfy
$$x_2F_g(x_0,x_1)+F_{g+1}(x_0,x_1) = -x_2F_{g+1}(x_0,x_1)-F_{g+2}(x_0,x_1) = 0.$$
Eliminating $x_2$, we get the equation $F_{g+1}^2-F_gF_{g+2} = 0$ which defines the set of the $2g+2$ ramification points $p_1,\ldots,p_{2g+2}$ of the projection $H_{g+2}\setminus \{q\} \to \bbP^1$.

Let
$$\Gamma:T_2F_g(T_0,T_1)+F_{g+1}(T_0,T_1) = 0$$
be the first polar $\Gamma$ of $H_{g+2}$
with respect to the point $q$. The  transformation $IH_{g+2}$ blows down $\Gamma$ and the lines $\la q,p_i\ra$  to points.
It follows immediately from \eqref{f1} that the set of fixed points of the involution  $IH_{g+2}$ outside the base locus is  equal to the hyperelliptic curve \eqref{hyp}.  Also we see that the pencil of lines through $q$ is invariant with respect to $IH_{g+2}$.

Let $\sigma:S\to \bbP^2$ be the blowup of the point $q$ and the
points $p_1,\ldots,p_{2g+2}$. The full preimage of the line $\ell_i = \la
q,p_i\ra$ consists  of two irreducible components, each  isomorphic to $\bbP^1$. They
intersect transversally at one point. We will call such a reducible curve, a \emph{bouquet} of two $\bbP^1$'s.  One component is the
exceptional curve $R_i = \sigma^{-1}(p_i)$ and  another one is the
proper  transform $R_i'$ of the line $\ell_i$. The proper
transform of $H_{g+2}$ intersects $\sigma^{-1}(\ell_i)$ at its
singular point. Thus the proper transform $\bar{H}_{g+2}$ of the hyperelliptic curve
$H_{g+2}$ intersects the exceptional curve $E= \sigma^{-1}(q)$ at
the same points where the proper transform of lines $\ell_i$
intersect $E$. The proper transform $\bar{\Gamma}$ of
$\Gamma$ intersects $R_i$ at one nonsingular point, and intersects
$E$ at $g$ points, the same points where the proper inverse
transform $\bar{H}_{g+2}$ of $H_{g+2}$ intersects $E$. The
involution $IH_{g+2}$ lifts to a biregular automorphism $\tau$ of
$S$. It switches the components $R_i$ and $R_i'$ of
$\sigma^{-1}(\ell_i)$, switches $E$ with $\bar{\Gamma}$ and fixes
the curve $\bar{H}_{g+2}$ pointwise. The pencil of lines through
$q$ defines a morphism $\phi:S\to \bbP^1$ whose fibres over the
points corresponding to the lines $\ell_i$ are isomorphic to a
bouquet of two $\bbP^1$'s. All other fibres are isomorphic to
$\bbP^1$. This is an example of a \emph{conic bundle} or a Mori
fibration (or in the archaic terminology  of \cite{Isk2}, a minimal
rational surface with a pencil of rational curves).

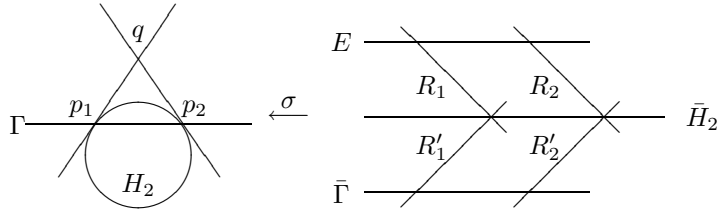
\begin{figure}[h]
\xy (-35,0)*{};
(5,5)*{};(-10.5,-18)**{}**\dir{-};
(-5,5)*{};(10.9,-18)**{}**\dir{-};
(0,-15)*\cir<20pt>{};
(-15,-11)*{};(15,-11)**{}**\dir{-};
(30,0)*{};(60,0)**{}**\dir{-};
(30,-20)*{};(60,-20)**{}**\dir{-};
(35,-22)*{};(49,-8)**{}**\dir{-};
(35,2)*{};(49,-12)**{}**\dir{-};
(50,-22)*{};(64,-8)**{}**\dir{-};
(50,2)*{};(64,-12)**{}**\dir{-};
(-16,-11)*{\Gamma};
(-7.5,-9)*{p_1};
(7.5,-9)*{p_2};
(0,1)*{q};
(27,0)*{E};
(30,-10)*{};(70,-10)**{}**\dir{-};
(75,-10)*{\bar{H}_2};
(27,-20)*{\bar{\Gamma}};
(39,-6)*{R_1};
(39,-14)*{R_1'};
(54,-6)*{R_2};
(54,-14)*{R_2'};
(20,-10)*{\longleftarrow};
(20,-8)*{\sigma};
(0,-19)*{H_2};
\endxy
\end{figure}

To show that the birational involutions $IH_{g+2}, g > 0,$ are not conjugate to each other or to a projective involution we use the following.

\begin{lemma}\label{fix} Let $G$ be a finite subgroup of $\Cr(2)$ and let $C_1,\ldots,C_k$ be non-rational irreducible curves on $\bbP^2$ such that each of them contains an open subset $C_i^0$ whose points are fixed under all $g\in G$. Then the set of birational isomorphism classes of the curves $C_i$ is an invariant of the conjugacy class of $G$ in $\Cr(2)$.
\end{lemma}

\begin{proof} Suppose $G = T\circ G'\circ T^{-1} $ for some  subgroup $G'$ of $\Cr(2)$ and some $T\in \Cr(2)$. Then, replacing $C_i^0$ by a smaller open subset we may assume that $T^{-1}(C_i^0)$ is defined and consists of fixed points of $G'$. Since $C_i$ is not rational, $T^{-1}(C_i^0)$ is not a point, and hence its  Zariski closure is a rational irreducible curve $C_i'$ birationally isomorphic to $C_i$ which  contains an open subset of fixed points of $G'$.
\end{proof}

Since a connected component of the  fixed locus of a finite
group of projective transformations is a line or a point, we see
that $IH_{g+2}$ is not conjugate to a subgroup of projective
transformations for any $g > 0$. Since  $IH_{g+2}$ is conjugate to
some involution \eqref{dj}, where $P(x)$ is determined by the
birational isomorphism class of $H_{g+2}$, we see from the previous
lemma that $IH_{g+2}$ is conjugate to $IH'_{g'+2}$ if and only if $g
= g'$ and the curves $H_{g+2}$ and $H'_{g+2}$ are birationally
isomorphic. Finally, let us look at the involution $IH_2$. It is a
quadratic transformation which is conjugate to the quadratic transformation
$\tau_2:(x,y) \mapsto (x,x/y)$.

A Jonqui\`eres involution \eqref{dj} is a special case of  a Cremona transformation of the form
$$(x,y) \mapsto (\frac{ax+b}{cx+d}, \frac{r_1(x)y+r_2(x)}{r_3(x)y+r_4(x)}),$$
where $a,b,c,d\in \bbC, ad-bc\ne 0$ and $r_i(x)\in \bbC(x)$ with
$r_1(x)r_4(x)-r_2(x)r_3(x) \ne 0$. These transformations form a
subgroup of $\Cr(2)$ called a  {\it de Jonqui\`eres subgroup} and
denoted by $\dJJ(2)$. Of course, its definition requires a choice of
a transcendence basis of the field $\bbC(\bbP^2)$.  If we identify
$\Cr(2)$ with the group $\Aut_\bbC(\bbC(x,y))$, and consider the
field $\bbC(x,y)$ as a field $K(y)$, where $K = \bbC(x)$, then
$$\dJJ(2) \cong  \PGL_{\bbC(x)}(2)\rtimes \PGL(2)$$
where  $\PGL(2)$ acts on $ \PGL_{\bbC(x)}(2)$ via Moebius transformations of the variable $x$.

It is clear that all elements from $\dJJ(2)$ leave the pencil of lines parallel to the $y$-axis invariant. One can show that a subgroup of $\Cr(2)$ which leaves a pencil of rational curves invariant is conjugate to $\dJJ(2)$.

\subsection{Geiser and Bertini involutions}\label{geiser}
The classical definition of a \emph{Geiser involution} is as follows
\cite{Geiser}. Fix 7 points $p_1,\ldots,p_7$  in $\bbP^2$ in general
position (we will make this more precise later). The linear system
$L$ of cubic curves through the seven points is two-dimensional. Take a
general point $p$ and consider the pencil of curves from $L$ passing  through $p$. Since a general  pencil of cubic curves has 9 base points, we can
define $\gamma(p)$ as the ninth base point of  the
pencil. One can show that the algebraic degree of a Geiser involution is equal to 8.
Another way to see a Geiser involution is  as follows. The linear
system $L$ defines a rational map of degree 2
$$f:\bbP^2-\to |L|^*\cong \bbP^2.$$
The points $p$ and $\gamma(p)$ lie in the same fibre. Thus $\gamma$
is a birational deck transformation of this cover. Blowing up the
seven points, we obtain a Del Pezzo surface $S$ of degree 2 (more
about this later), and a regular map of degree 2 from $S$ to
$\bbP^2$. The Geiser involution $\gamma$ becomes an automorphism of
the surface $S$.

It is easy to see that the fixed points of a Geiser involution lie on the ramification curve of $f$. This curve is a curve of degree 6 with double points at the  points $p_1,\ldots,p_7$. It is birationally isomorphic to a canonical curve of genus 3. Applying Lemma \ref{fix}, we obtain that a Geiser involution is not conjugate to any de Jonqui\`eres involution $IH_{g+2}$. Also, as we will see later, the conjugacy classes of Geiser involutions are in a bijective correspondence with the moduli space of canonical curves of genus 3.

\medskip
To define a \emph{Bertini involution} we  fix 8 points in $\bbP^2$ in general position and consider
the pencil of cubic curves through these points. It has the ninth
base point $p_9$. For any general point $p$ there will be a unique
cubic curve $C(p)$ from the pencil which passes through $p$. Take
$p_9$ for the  zero in the group law of the cubic $C(p)$ and define
$\beta(p)$ as the negative $-p$ with respect to the group law. This
defines a birational involution on $\bbP^2$, a Bertini
involution \cite{Bertini}. One can show that the algebraic degree of a Bertini involution is equal
to 17. We will see later that the fixed points
of a Bertini involution lie on a canonical curve of genus 4 with
vanishing theta characteristic (isomorphic to a nonsingular
intersection of a cubic surface and a quadric cone in $\bbP^3$).
So, a Bertini involution is not conjugate to a Geiser involution or
a de Jonqui\`eres involution. It can be realized as an automorphism
of the blowup of the eight points (a Del Pezzo surface of degree
1), and the quotient by this involution is isomorphic to a quadratic cone in $\bbP^3$.

\section{Rational $G$-surfaces}
 \subsection{Resolution of indeterminacy points}
 Let $\chi:S -\to S'$ be a birational map of nonsingular projective surfaces. It is well-known (see \cite{Hartshorne}) that there exist birational morphisms of nonsingular surfaces  $\sigma:X\to S$ and $\phi:X\to S'$ such that the following diagram is
 commutative
\begin{equation}\label{hironaka}
\xymatrix{&X{}\ar[ld]_\sigma\ar[rd]^\phi&\\
S\ar@{-->}[rr]&&S'.}
\end{equation}
It is called a \emph{resolution of indeterminacy points} of $\chi$. Recall also  that  any birational morphism can be factored into a finite sequence  of blowups of  points. Let
\begin{equation}\label{decom}
\sigma:X = X_N\overset{\sigma_N}{\longrightarrow} X_{N-1}\overset{\sigma_{N-1}}{\longrightarrow} \cdots \overset{\sigma_2}{\longrightarrow} X_1\overset{\sigma_1}{\longrightarrow} X_0 = S
\end{equation}
be such a factorization.  Here $\sigma_i:X_i\to X_{i-1}$ is the blowup of a point $x_i\in X_{i-1}$. Let
\begin{equation}\label{excl1}
E_i = \sigma_i^{-1}(x_i),\quad \calE_i = (\sigma_{i+1}\circ \ldots\sigma_{N})^{-1}(E_i).
\end{equation}
Let $H'$ be a very ample divisor class  on $S'$ and $\calH'$ be the
corresponding complete linear system $|H'|$. Let $\calH_N =
\phi^*(\calH')$. Define $m(x_N)$ as the smallest positive number
such that $\calH_N+m(x_N)E_N = \sigma_N^*(\calH_{N-1})$ for some
linear system $\calH_{N-1}$  on $X_{N-1}$. Then proceed inductively
to define linear systems $\calH_k$ on each $X_k$ such that
$\calH_{k+1} +m(x_{k+1})E_{k+1} = \sigma_{k+1}^*(\calH_k)$, and
finally a linear system $\calH = \calH_0$ on $S$ such that
$\calH_1+m(x_1)E_1 = \sigma_1^*(\calH)$. It follows from the
definition that
\begin{equation}\label{linsys}
\phi^*(\calH') = \sigma^*(\calH)-\sum_{i=1}^Nm(x_i)\calE_i.
\end{equation}
The proper  transform of $\calH'$ on $S$ under $\chi$  is contained in the linear system $\calH$. It consists of curves which pass through the points $x_i$ with multiplicity $\ge m_i$. We denote it by
$$\chi^{-1}(\calH') = |H-m(x_1)x_1-\cdots-m(x_N)x_N|,$$
where $\calH \subset |H|$. Here for a curve on $S$ to pass through a
point $x_i\in X_{i-1}$ with multiplicity $\ge m(x_i)$ means that the
proper  transform of the curve on $X_{i-1}$  has $x_i$ as a
point of multiplicity $\ge m(x_i)$. The  divisors $\calE_i$ are
called the \emph{exceptional curves} of the resolution
$\sigma:X\to S$ of the birational map $\chi$. Note that $\calE_i$ is
an irreducible curve if and only if $\sigma_{i+1}\circ \ldots\circ
\sigma_N:X\to X_i$ is an isomorphism over $E_i =\sigma^{-1}(x_i)$.

The set of points $x_i\in X_i, i = 1,\ldots, N$, is called the set
of \emph{indeterminacy points}, or {\it base points}, or {\it
fundamental points} of $\chi$. Note that, strictly speaking,  only one
of them, $x_1$, lies in $S$. However, if $\sigma_{1}\circ \ldots
\circ\sigma_i:X_i\to S$ is an isomorphism in a neighborhood of
$x_{i+1}$ we can identify this point with a point in $S$. Let
$\{x_i, i\in I\}$ be the set of such points. Points $x_j, j\not\in
I,$ are {\it infinitely near points}. A precise meaning of this
classical notion is as follows.

Let $S$ be a nonsingular projective surface and $\calB(S)$ be the category  of birational morphisms $\pi:S'\to S$ of nonsingular projective surfaces. Recall that a morphism from $(S'\overset{\pi'}{\to} S)$ to $(S''\overset{\pi''}{\to} S)$ in this category is a regular map $\phi:S'\to S''$ such that
$\pi''\circ \phi = \pi'$.

\begin{definition}\label{bubble}The \emph{bubble space} $S^{\bb}$ of a nonsingular surface $S$ is the factor set
$$ S^{\bb} = \Bigl(\bigcup_{(S'\overset{\pi'}{\to} S)\in \calB(S)}S'\Bigr)/R,$$
where $R$ is the following equivalence relation: $x'\in S'$ is equivalent to $x''\in S''$ if
the rational map $\pi''{}^{-1}\circ \pi':S'-\to S''$  maps isomorphically an open neighborhood of $x'$ to an open neighborhood of $x''$.

\end{definition}

It is clear that for any $\pi:S'\to S$ from $\calB(S)$ we have an injective map $i_{S'}:S'\to S^{\bb}$.
We will identify points of $S'$ with their images. If $\phi:S''\to S'$ is a morphism in $\calB(S)$ which is isomorphic in $\calB(S')$ to  the blowup of a point $x'\in S'$, any point $x''\in \phi^{-1}(x')$ is called {\it infinitely near point}\index{infinitely near point} to $x'$ of the first order. This is denoted by
$x''\succ x'$. By induction, one defines an infinitely near point of order $k$, denoted by $x''\succ_k x'$.  This defines a partial order on $S^{\bb}$.

We say that a point $x\in S^{\bb}$ is of height $k$, if $x\succ_k x_0$ for some $x_0\in S$. This defines the {\it height function}\index{height function} on the bubble space
$$\htt_S:S^{\bb}\to \bbN.$$
Clearly, $S = \htt^{-1}(0)$.

It follows from the known behavior of the canonical class under a blowup that
\begin{equation}\label{can}
K_X = \sigma^*(K_S) + \sum_{i=1}^N\calE_i.
\end{equation}
The intersection theory on a nonsingular surface gives
\begin{eqnarray}\label{for2}
\calH'^2 &= &(\phi^*(\calH'))^2 =
(\sigma^*(\calH)-\sum_{i=1}^Nm(x_i)\calE_i)^2 =
\calH^2-\sum_{i=1}^Nm(x_i)^2,\\ \notag 
K_{S'}\cdot \calH'
&=&K_S\cdot \calH+\sum_{i=1}^Nm(x_i).\notag
\end{eqnarray}

\begin{example} Let $\chi:\bbP^2\ -\to \bbP^2$ be a Cremona transformation, $\calH' = |\ell|$ be the linear system of lines in $\bbP^2$, and $\calH \subset |n\ell|$. The formulas \eqref{for2} give
\begin{eqnarray}\label{eqq}
n^2-\sum_{i=1}^Nm(x_i)^2& = &1,\\ \notag
3n -\sum_{i=1}^N m(x_i) &=
&3.\notag
\end{eqnarray}
The linear system $\calH$ is written in this situation as $\calH =
|n\ell -\sum_{i=1}^Nm_ix_i|.$ For example, a quadratic
transformation with 3 base points $p_1,p_2,p_3$ is given by the
linear system $|2\ell-p_1-p_2-p_3|.$ In the case of the standard
quadratic transformation $\tau_1$ the curves
$\calE_1,\calE_2,\calE_3$ are irreducible,  the map $\sigma_1:
X_1\to X_0 = \bbP^2$ is an isomorphism in a neighborhood of
$p_2,p_3$ and the map $\sigma_2:X_2\to X_1$ is an isomorphism in a
neighborhood of $\sigma^{-1}(p_3)$. This shows that  we can identify
$p_1,p_2,p_3$ with points on $\bbP^2$. In the case of the
transformation \eqref{st2}, we have $\sigma_1(p_2) = p_1$ and $p_3$
can be identified with a point on $\bbP^2$. So in this case
$p_2\succ p_1$. For the transformation~(\ref{st3}) we have $p_3\succ
p_2 \succ  p_1$.

For a Geiser involution (resp. a Bertini involution) we have
$\calH = |8\ell-3p_1-\cdots-3p_7|$ (resp. $\calH = |17\ell-6p_1-\cdots-6p_8|$).
\end{example}

\subsection{$G$-surfaces} Let $G$ be a finite group. A  $G$-surface is a pair $(S,\rho)$, where $S$ is a nonsingular projective surface and $\rho$ is an isomorphism from $G$ to a group of automorphisms of $S$. A morphism of the pairs $(S,\rho)\to (S',\rho')$ is defined to be a morphism of surfaces $f:S\to S'$ such that $\rho'(G') = f \circ \rho(G)\circ f^{-1}$. In particular, two subgroups of $\Aut(S)$ define isomorphic $G$-surfaces if and only if they are conjugate inside of $\Aut(S)$.  Often, if no confusion arises, we will denote a $G$-surface by $(S,G)$.

Let $\chi:S\to S'$ be a birational map of $G$-surfaces. Then one can
$G$-equivariantly resolve $\chi$, in the sense that one can find the
diagram \eqref{hironaka} where all maps are morphisms of  $G$-surfaces.
The group  $G$ acts on the surface $X$ permuting the exceptional
configurations $\calE_i$ in such a way that $\calE_i \subset
\calE_j$ implies $g(\calE_i) \subset g(\calE_j)$. This defines an
action of $G$ on the set of indeterminacy points of $\chi$ ($g(x_i)
= x_j$ if $g(\calE_i) = g(\calE_j)$). The action preserves the
order, i.e. $x_i\succ x_j$ implies $g(x_i)\succ g(x_j)$, so the
function $\htt:\{x_1,\ldots,x_N\}\to \bbN$ is constant on each orbit
$Gx_i$.

Let $\calH' = |H'|$ be an ample linear system on $S'$  and
$$\phi^*(\calH' )= \sigma^*(\calH)-\sum_{i=1}^Nm(x_i)\calE_i$$
be its  inverse transform on $X$ as above. Everything here is $G$-invariant, so $\calH$ is a $G$-invariant linear system on $S$ and the multiplicities $m(x_i)$ are constant on the $G$-orbits. So we can rewrite the system in the form
$$\phi^*(\calH' ) = \sigma^*(\calH)-\sum_{\kappa\in \calI}m(\kappa)\calE_{\kappa},$$
where $\calI$ is the set of $G$-orbits of indeterminacy points. For any $\kappa\in \calI$ we set $m(\kappa) = m(x_i)$, where  $x_i\in \kappa$ and
$\calE_{\kappa}, = \sum_{x_i\in \kappa} \calE_i$. Similarly one can rewrite the proper transform of $\calH'$ on $S$
\begin{equation}\label{kappa}
|H-\sum_{\kappa\in \calI}m(\kappa)\kappa|.\end{equation}
Now we can  rewrite the intersection formula \eqref{for2} in the form
\begin{eqnarray}\label{gfor2}
H'{}^2&= &H^2-\sum_{\kappa\in \calI}m(\kappa)^2d(\kappa)\\ \notag
K_{S'}\circ H' &=&K_S\cdot H+\sum_{\kappa\in \calI}m(\kappa)d(\kappa),\notag
\end{eqnarray}
where $d(\kappa) = \#\{i:x_i\in \kappa\}$.

\begin{remark} In the arithmetical analog of the previous theory all the notation become  more natural. Our maps are maps over a perfect ground field $k$. A blowup is the blowup of a closed point in the scheme-theoretical sense, not necessary $k$-rational. An exceptional curve is defined over $k$ but when we replace $k$ with an algebraic closure $\bar{k}$, it splits into the union of conjugate exceptional curves over $\bar{k}$.   So, in the above notation, $\kappa$ means a closed point on $S$ or an infinitely near point. The analog of  $d(\kappa)$ is of course the degree of a point, i.e. the extension degree $[k(x):k]$, where $k(x)$ is the residue field of $x$.
\end{remark}

\subsection{The $G$-equivariant bubble space} Here we recall Manin's formalism of the theory of linear systems with base conditions in its $G$-equivariant form (see \cite{Manin2}).

First we define the \emph{$G$-equivariant bubble space} of a $G$-surface $S$ as a $G$-equivariant version $(S,G)^{\bb}$ of Definition \ref{bubble}. One replaces the category $\calB(S)$ of birational morphisms $S'\to S$ with the category $\calB(S,G)$ of birational morphisms of $G$-surfaces. In this way the group $G$ acts on the bubble space $(S,G)^{\bb}$.  Let
\begin{equation}
Z^*(S,G) =\underrightarrow{\lim} \Pic(S'),
\end{equation}
where the inductive limit is taken with respect to the functor $\Pic$ from the category $\calB(S,G)$ with values in the category of abelian groups defined by $S'\to \Pic(S')$. The group $Z^*(S,G)$ is equipped with a natural structure of $G$-module. Also it is equipped with the following natural structures.
\begin{itemize}
\item[(a)] A symmetric $G$-invariant pairing
$$Z^*(S,G) \times Z^*(S,G)\to \bbZ$$
induced by the intersection pairing on each $\Pic(S')$.
\item[(b)]   A distinguished cone of effective divisors classes in $Z^*(S,G)$
$$Z_+^*(S,G) = \underrightarrow{\lim}\Pic_+(S'),$$
where $\Pic_+(S')$ is the cone of effective divisor classes on each  $S'$ from $\calB(S,G)$.
\item[(c)] A distinguished $G$-equivariant homomorphism
$$K:Z^*(S,G) \to \bbZ, \quad K(z) = K_{S'}\cdot z,\  \text{for any}\ S'\to S\  \text{from}\ \calB(S,G).$$
\end{itemize}

Let $f:S'\to S$ be a morphism from $\calB(S,G)$ and $\calE_1,\ldots,\calE_n$ be its  exceptional curves. We have a natural splitting
$$\Pic(S') = f^*(\Pic(S))\oplus \bbZ[\calE_1]\oplus \cdots \oplus \bbZ[\calE_n].$$
Now let $Z_0(S,G) = \bbZ^{(S,G)^{\bb}}$ be the free abelian group generated by the set $(S,G)^{\bb}$. Identifying exceptional curves with points in the bubble space, and passing to the limit we obtain  a natural splitting
\begin{equation}\label{split}
Z^*(S,G) =  Z_0(S,G)\oplus \Pic(S).\end{equation}
Passing to invariants we get the splitting
\begin{equation}\label{gsplit}
Z^*(S,G)^G =  Z_0(S,G)^G\oplus \Pic(S)^G.
\end{equation}
Write an element of $Z^*(S,G)^G $ in the form
$$z = D-\sum_{\kappa\in O}  m(\kappa)\kappa,$$
where $O$ is the set of $G$-orbits  in $Z_0(S,G)^G$ and $D$ is a $G$-invariant divisor class on $S$. Then \begin{itemize}
\item[(a)]
$z\cdot z' = D\cdot D' -\sum_{\kappa\in O} m(\kappa)m'(\kappa)d(\kappa);$
\item[(b)] $z\in Z_+^*(S,G)$ if and only if $D\in \Pic_+(S)^G, m(\kappa)\ge 0$ and
$m(\kappa')\le m(\kappa) $ if $\kappa'\succ \kappa$;
\item[(c)] $K(z) = D\cdot K_S+\sum_{\kappa\in O} m(\kappa)d(\kappa).$
\end{itemize}

Let $\phi:S'\to S$ be an object of $\calB(S,G)$. Then we have a natural map
$\phi_{\bb}:(S',G)^{\bb}\to (S,G)^{\bb}$ which induces an isomorphism $\phi_{\bb}^{*}:Z(S,G) \to Z(S',G)$. We also have a natural  isomorphism
$\phi_*^{\bb}: Z(S',G)\to Z(S,G)$. None of these maps   preserves the splitting \eqref{split}. Resolving indeterminacy points of any birational map $\chi:(S,G)-\to (S',G')$ we can define
\begin{itemize}
\item \emph{proper direct transform map} $\chi_*:Z^*(S,G)\overset{\sim}{\to} Z^*(S',G)$;
\item \emph{proper inverse transform map} $ \chi^*:Z^*(S',G)\overset{\sim}{\to} Z^*(S,G)$.
\end{itemize}
The group $Z^*(S,G)$ equipped with all above structures  is one of the main $G$-birational invariants of $S$. It can be viewed as the Picard group of the bubble space $(S,G)^{\bb}$.

The previous machinery gives a convenient way to consider the linear
systems defining rational maps of surfaces. Thus we can rewrite
\eqref{linsys} in the form $|z|$, where $z= H-\sum m_ix_i$ is
considered as an element of $Z_+^*(S,G)$.  The condition that $|z|$ be 
homaloidal becomes equivalent to the conditions
\begin{eqnarray}
z^2& = &H^2-\sum m_i^2 = H'{}^2,\\ \notag 
K(z)&= &H\cdot K_S+\sum m_i
= H'\cdot K_{S'}.\notag
\end{eqnarray}
When $S = S' =\bbP^2$ we get the equalities \eqref{eqq}.

\subsection{Minimal rational $G$-surfaces}\label{minimal} Let $(S,\rho)$ be a rational $G$-surface. Choose a birational map $\phi:S- \to \bbP^2$. For any $g\in G$, the map
$\phi\circ  g\circ \phi^{-1}$ belongs to  $\Cr(2)$. This defines an injective  homomorphism
\begin{equation}
\iota_\phi:G\to \Cr(2).
\end{equation}
Suppose $(S',\rho')$ is another  rational $G$-surface and $\phi':S'-\to \bbP^2$ is a birational map.

The following lemma is obvious.

\begin{lemma} The subgroups $\iota_\phi(G)$ and $\iota_{\phi'}(G)$ of $\Cr(2)$ are conjugate if and only if there exists a birational map of $G$-surfaces $\chi:S'-\to S$.
\end{lemma}

The lemma shows that a birational isomorphism class of $G$-surfaces defines a conjugacy class of subgroups of $\Cr(2)$ isomorphic to $G$.
The next lemma shows that any conjugacy class is obtained in this way.

\begin{lemma} Suppose  $G$ is a finite subgroup of $\Cr(2)$, then there exists a rational $G$-surface
$(S,\rho)$  and a birational map $\phi:S\to \bbP^2$ such that
$$G= \phi\circ\rho(G)\circ \phi^{-1}.$$
\end{lemma}

\begin{proof} We give two proofs. The fist one is after A. Verra. Let $D = \cap_{g\in G}\text{dom}(g)$, where $\text{dom}(g)$ is an open subset on which $g$ is defined. Then $U = \cap_{g\in G} g(D)$ is an open invariant subset of $\bbP^2$ on which $g\in G$ acts biregularly. Order $G$ in some way and consider a copy of $\bbP_g^2$ indexed by $g\in G$. For any $u\in U$ let $g(u)\in \bbP_g^2$. We define a morphism
$$\phi:U\to \prod_{g\in G} \bbP_g^2, \ u\mapsto (g(u))_{g\in G}.$$
Define an action of $G$ on $\phi(U)$  by $g'((x_g)_{g\in G}) = (x_{gg'})_{g\in G}$. Then $\phi$ is obviously $G$-equivariant. Now define $V$ as the Zariski closure of $\phi(U)$ in the product. It is obviously a $G$-invariant surface  which contains an open $G$-invariant subset $G$-isomorphic to $U$. It remains to replace $V$ by its $G$-equivariant resolution of singularities (which always exists).

The second proof is standard. Let $U$ be as above and $U' = U/G$ be
the orbit space. It is a normal algebraic surface. Choose any normal
projective completion $X'$ of $U'$. Let $S'$ be the normalization of
$X'$ in the field of rational functions of $U$. This is a normal
projective surface on which $G$ acts by biregular transformations.
It remains to define $S$ to be a $G$-invariant resolution of
singularities (see also \cite{dFe-Ein}).
\end{proof}

Summing up, we obtain the following result.

\begin{theorem} There is a natural bijective correspondence between birational isomorphism classes of rational $G$-surfaces and conjugate classes of subgroups of $\Cr(2)$ isomorphic to $G$.
\end{theorem}

So our goal is to classify $G$-surfaces $(S,\rho)$ up to birational isomorphism of $G$-surfaces.

\begin{definition} A \emph{minimal} $G$-surface is a $G$-surface $(S,\rho)$ such that any birational morphism of $G$-surfaces $(S,\rho)\to (S',\rho')$ is an isomorphism. A group $G$ of automorphisms of a rational surface $S$ is called a {\it minimal group of automorphisms} if the pair $(S,\rho)$ is minimal.
\end{definition}

Obviously, it is enough to classify minimal rational $G$-surfaces up to birational isomorphism of $G$-surfaces.

Before we state the next fundamental result, let us recall some terminology.

A \emph{conic bundle structure} on a rational $G$-surface $(S,G)$ is
a $G$-equivariant morphism $\phi:S\to \bbP^1$ such that the fibres
are isomorphic to a reduced conic in $\bbP^2$. A \emph{Del Pezzo
surface} is a surface with ample anti-canonical divisor $-K_S$.

\begin{theorem}\label{enriques} Let $S$ be a minimal rational $G$-surface. Then either $S$ admits a structure of a conic bundle with $\Pic(S)^G \cong \bbZ^2$, or $S$ is isomorphic to a Del Pezzo surface with $\Pic(S)^G\cong \bbZ$.
\end{theorem}

An analogous result  from  the classical literature is proven by
using  the method of the termination of adjoints, first introduced
for linear system of plane curves in the work of G. Castelnuovo. It
consists in replacing a linear system $|D|$ with the linear system
$|D+K_S|$  and repeat doing this to stop only if the next step leads
to the empty linear system. The application of this method to
finding a $G$-invariant linear system of curves on the plane was
initiated in the works of S. Kantor \cite{Kantor}, who essentially
stated the theorem above but without the concept of minimality. In
arithmetical situation this method was first applied by F. Enriques
Yu. Manin \cite{Manin} and by the second author \cite{Isk2} (an earlier
proof of Manin used the assumption that $G$ is an abelian group).
Nowadays the theorem follows easily from a $G$-equivariant version
of Mori  theory (see \cite{Kollar-Mori}, Example 2.18) and the
proof can be found in literature (\cite{BB}, \cite{deFernex}). For
this reason we omit the proof.

 Recall the classification of  Del Pezzo surfaces (see \cite{Demazure}, \cite{Manin2}).  The number $d = K_S^2$ is called the \emph{degree}. By Noether's formula, $1 \le d\le 9$. For $d\ge 3$, the anti-canonical linear system $|-K_S|$ maps $S$  onto a nonsingular surface of degree $d$ in $\bbP^d$. If $d = 9$, $S\cong \bbP^2$. If $d = 8$, then $S\cong \bbP^1\times \bbP^1$, or $S\cong \bfF_1$, where as always we denote by $\bfF_{n}$ the minimal ruled surface $\bbP(\calO_{\bbP^1}\oplus \calO_{\bbP^1}(n))$. For
$d\le 7$, a Del Pezzo surface $S$ is isomorphic to the blowup of $n= 9-d$ points in $\bbP^2$ in general position that means that 
\begin{itemize}\label{gen}
\item no three are on a line;
\item no six are on a conic;
\item  if $n=8$, then then the points are not  on a  plane cubic which has one of them as its singular point.
\end{itemize}
For $d=2$, the linear system $|-K_S|$ defines a finite morphism  of degree 2 from $S$ to $\bbP^2$ with a nonsingular quartic as the branch curve. Finally,  for $d = 1$, the linear system $|-2K_S|$ defines a finite morphism of degree 2 onto a quadric cone $Q\subset \bbP^{3}$ with the branch curve cut out by a cubic.

For a minimal Del Pezzo $G$-surface the group $\Pic(S)^G$ is generated by $K_S$ if $S$ is not isomorphic to $\bbP^2$ or $\bbP^1\times \bbP^1$. In the latter cases it is generated by $\frac{1}{3}K_S$ or $\frac{1}{2}K_S$, respectively.

A conic bundle surface is either isomorphic to $\bfF_n$ or to a surface obtained from $\bfF_n$  by blowing up a finite set of points, no two lying in a fibre of a  ruling. The number of blowups is equal to the number of singular fibres of the conic bundle fibration. We will exclude the surfaces $\bfF_0$ and $\bfF_1$, considering them as Del Pezzo surfaces.

There are minimal conic bundles with ample $-K_{S}$ (see Proposition \ref{cdp}).

\section{Automorphisms of minimal ruled surfaces}\label{sec3}
\subsection{Some of group theory} We employ the standard notations for groups used by group-theorists (see  \cite{Atlas}):
\smallskip

$C_n$,  a cyclic group of order $n$,

$n = C_n$ if no confusion arises;

$n^r = C_n^r$, the direct sum of $r$ copies of $C_n$ (not to be confused with cyclic group of order $n^r$);

$S_n$, the \emph{permutation group} of degree $n$;

$A_n$, the \emph{alternating group}  of degree $n$;

$D_{2n}$, the \emph{dihedral group} of order $2n$;

$Q_{4n} = \la a,b\mid a^{2n} =1, b^2 = a^n, b^{-1}ab = a^{-1}\ra$, \emph{dicyclic group} of order $4n$, a \emph{generalized quaternion group} if $n = 2^k$;

$H_n(p)$, the \emph{Heisenberg group} of  unipotent $n\times n$-matrices with entries in $\bbF_p$.

$\GL(n) = \GL(n,\bbC)$, general linear group over $\bbC$,

$\SL(n) = \SL(n,\bbC)$, special linear group over $\bbC$,

$\PGL(n) = \GL(n,\bbC)/\bbC^*$, general projective linear group over $\bbC$,

$\Or(n)$, the orthogonal  linear group over $\bbC$,

$\PO(n)$, the projective orthogonal  linear group over $\bbC$,

$L_n(q) = \PSL(n,\bbF_q),$ where $q = p^r$ is a power of a prime
number $p$;

$T \cong A_4, \ O \cong S_4 \cong \PGL(2,\bbF_3), \ I \cong A_5 \cong L_2(5)\cong L_2(2^2)$, tetrahedral, octahedral , icosahedral  subgroups of $\PGL(2)$;

$ \overline{T} \cong \SL(2,\bbF_3), \ \bar{O} \cong \GL(2,\bbF_3), \ \bar{I} \cong \SL(2,\bbF_5), \overline{D}_{2n} \cong Q_{4n}$,
binary tetrahedral, binary octahedral ,binary icosahedral, binary dihedral   subgroups of $\SL(2)$;

$A${\SMALL$\bullet$}$B$ is an upward extension of $B$ with help of a normal subgroup $A$;

$A: B$ is a split extension $A\pt B$, i.e.  a semi-direct product $A\rtimes B$ ((it is defined by a homomorphism  $\phi: B\to \Aut(A));$

$A^\bullet B$
  is a non-split extension $A\pt B$,

$nA = n^\bullet A$, where the normal group $n$ is equal to the center.

$p^{a+b} = C_p^a\pt C_p^b, $ where $p$ is prime;

$A\trian B$ (or $A\trian_D B$), the diagonal product of $A$ and $B$ over their common homomorphic image $D$ (i.e. the subgroup of $A\times B$ of pairs $(a,b)$ such that $\alpha(a) = \beta(b)$ for some  surjections $\alpha:A\to D, \beta: B\to D$). When $D$ is omitted it means that $D$ is the largest possible;

$\frac{1}{m} [A\times B] = A\trian_D B$, where $\# D = m$;

$A\wr S_n $,  the wreath product, i.e. $A^n:S_n$, where $S_n$ is the symmetric group acting on $A^n$ by permuting the factors;


$\mu_n$,  the group of $n$th roots of unity with generator $\epsilon_n = e^{2\pi i/n}$.
\medskip

We will often use  the following simple result from group theory which is known as Goursat's Lemma.

\begin{lemma}\label{goursat} Let $G$ be a finite subgroup of the product $A\times B$ of two groups $A$ and $B$. Let $p_1:A\times B \to A, p_2:A\times B\to B$ be the projection homomorphisms. Let
$G_i = p_i(G), H_i = \Ker(p_j|G), i\ne j = 1,2. $ Then $H_i$ is a normal subgroup in $G_i$. The map $\phi:G_1/H_1\to G_2/H_2$ defined by
 $\phi (aH_1) = p_2(a)H_2$ is an isomorphism, and
$$G = G_1\trian_DG_2,
$$
where $D = G_1/H_1, \alpha:G_1\to D$ is the projection map to the
quotient, and $\beta$ is the composition of the projection $G_2\to
G_2/H_2$ and $\phi^{-1}$.
\end{lemma}

Note some special cases:
$$G\trian_{1}G'\cong G\times G',  \quad G\trian_{G^\prime}G' = \{(g,\alpha(g))\in G\times G', g\in G_1\},$$
where $\alpha:G\to G'$  is a surjection and $G'\to G'$ is the identity.

We will be dealing with various group extensions. The following  lemma  is known  in group theory as the Schur-Zassenhaus Theorem. Its proof can be found in \cite{Gorenstein}, 6.2.

\begin{lemma}\label{zas} Let $A\pt B$ be an extension of groups. Suppose that the orders of $A$ and $B$ are coprime. Then the extension splits. If,  moreover, $A$ or $B$ is solvable, then all subgroups of $A:B$ defining splittings are conjugate.
\end{lemma}

We will often use the following simple facts, their proofs are left to the reader (or can be found in www.planetmath.org).

\begin{lemma} \label{dihedral} A subgroup of $D_{2n}  = \la a,b\mid a^{n} = b^2 = b^{-1}aba = 1\ra$ is either cyclic or dihedral. A normal subgroup $H$ is either cyclic subgroup  $\la a\ra$, or  $n = 2k$ and $H$ is one of the following two subgroups $ \la a^2,b\ra, \la a^2,ab\ra$ of index 2. These two subgroups are interchanged under the outer automorphism $a\mapsto a, b\mapsto ab$. If $H$ is cyclic of order $n/k$, the quotient group is isomorphic to $D_{2n/k}$.

The  group of $\Aut(D_{2n})$ is isomorphic to $(\bbZ/n)^*:n$ and it is generated by the transformations $a\mapsto a^{s}, b\mapsto a^tb$. The subgroup of inner automorphisms is generated by transformations $a\mapsto a^{-1}, b\mapsto b$ and $a\mapsto a, b\mapsto a^2b$.
\end{lemma}

It will be convenient to list all isomorphism classes of non-abelian groups of order  $16$.

 \begin{table}[h]
\begin{center}
{\small\begin{tabular}{||r |r | r|r|r ||}
\hline
Notation&Center&LCS&Extensions&Presentation\\ \hline
$2\times D_8$&$2^2$&16,\ 2,\ 1&$2^{1+3},2^{2+2},2^{3+1}$,&$a^4=b^2=c^2=1,$\\
&&&$(2\times 4):2$&$[a,b]a^2 = [a,c]= [b,c] = 1$\\ \hline
$2\times Q_8$&$2^2$&16,\  2,\ 1&$2^{2+2},$&$a^4 = a^2b^{-2}=a^2[a,b] = 1$\\
&&&$(2\times 4)^\bullet 2$&\\
\hline
$D_{16}$&$2$&16,\ 4,\  2,\ 1&$8:2, 2D_8,$&$a^8 = b^{2}=a^2[a,b] = 1$\\
&&&$(2^2)^\bullet 4, D_8^\bullet 2$&\\\hline
$SD_{16}$&$2$&16,\ 4,\  2,\ 1&$8:2,D_8^\bullet 2,$&$a^8 = b^{2}=[a,b]a^{-2} = 1$\\
&&&$2D_8, (2^2)^\bullet 4$&\\\hline
$Q_{16}$&$2$&16,\ 4,\  2,\ 1&$8^\bullet 2,2D_8, $&$a^8 = a^4b^{-2}=[a,b]a^2 = 1$\\
&&&$4^\bullet(2^2)$&\\\hline
$AS_{16}$&$4$&16,\  2,\ 1&$2^{1+3},D_8:2$&$a^4 = b^2 = c^2 =[a,b] = 1$\\
&&&$4(2^2),(2\times 4):2$&$[c,b]a^{-2} = [c,a]= 1$\\ \hline
$K_{16}$&$2^2$&16,\ 2,\ 1&$2^{2+2},(2\times 4)^\bullet 2$&$a^4 = b^4 = [a,b] a^2 =1$\\
&&&$2^\bullet (2\times 4), 4: 4$&\\ \hline
$L_{16}$&$2^2$&16,\  2,\ 1&$2^2:4,2^\bullet (2\times 4)$&$a^4=b^2 =c^2 = 1,$\\
&&&$(2\times 4):2$&$ [c,a]b= [a,b] = [ c,b] =1$\\ \hline
$M_{16}$&$4$&16,\  2,\ 1&$8:2,\  4(2^2)$&$a^8 = b^2 = 1, [a,b]a^{4}=1$\\
&&&$2^\bullet (2\times 4)$&\\ \hline

\end{tabular}}
\end{center}
\caption{Non-abelian groups of order 16}\label{order8}
\end{table}

Recall that there are two non-isomorphic non-abelian groups of order 8, $D_8$ and $Q_8$.

Finally we describe central extension of polyhedral and binary polyhedral groups. Recall that the isomorphism classes of central extensions $A\pt G$, where $A$ is an abelian group, are parametrized by the 2-cohomology group $H^2(G,A)$. We will assume that $A \cong p$, where $p$ is prime. We will use the following facts about the cohomology groups  of polyhedral and binary polyhedral groups which can be found in  text-books on group cohomology (see, for example, \cite{Adem}). 

\begin{lemma}\label{coh} Let $G$ be a polyhedral group or a binary polyhedral group. If $G \cong n$ is cyclic, then $H^2(G,p) \cong p$ if $p\mid n$ and zero otherwise. If $G$ is not cyclic, then $H^2(G,p) = 0$ if $p\ne 2,3$. Moreover
\begin{itemize}
\item[(i)] If $G$ is a polyhedral group, then
$$H^2(G,2) \cong \begin{cases}
      2&\ \text{if}\  G\cong D_{2n},  n \ \text{odd},\\
      2^3&\ \text{if}\  G\cong D_{2n},  n \ \text{even},\\
        2&\ \text{if}\  G\cong T,  \\
          2^2&\ \text{if}\  G\cong O,  \\
    2&\ \text{if}\ G\cong I.
\end{cases}$$
$$H^2(G,3) \cong \begin{cases}
     3& \ \text{if}\ G\cong T,\\
     1& \text{otherwise}.
\end{cases}
$$
\item[(ii)] If $G$ is a binary polyhedral group, then
$$H^2(G,2) \cong \begin{cases}
     2&\ \text{if}\  G\cong \overline{D}_{2n}, n \ \text{odd},\\
2^2& \ \text{if}\ G\cong \overline{D}_{2n}, n \ \text{even},\\
 2& \ \text{if}\ G\cong \overline{O},\\
  1&\  \text{otherwise}.
\end{cases}
$$
$$H^2(G,3) \cong \begin{cases}
     3&\ \text{if}\  G\cong \overline{T},\\
     1&\  \text{otherwise}.
\end{cases}
$$
\end{itemize}
\end{lemma}

\subsection{Finite groups of projective automorphisms} We start with the case $S = \bbP^2$, where $\Aut(S) \cong \PGL(3)$. To save space we will often denote a projective transformation
$$(x_0,x_1,x_2)\mapsto (L_0(x_0,x_1,x_2),L_1(x_0,x_1,x_2),L_2(x_0,x_1,x_2))$$ by
$[L_0(x_0,x_1,x_2),L_1(x_0,x_1,x_2),L_2(x_0,x_1,x_2)]$.

Recall some standard terminology from the theory of linear groups. Let $G$ be a subgroup of the general linear group $\GL(V)$ of a complex vector space $V$. The group $G$ is called \emph{intransitive}  if the representation of $G$ in $V$ is reducible. Otherwise it is called  \emph{transitive}. A transitive group $G$ is called \emph{imprimitive} if it contains an intransitive normal subgroup $G'$.  In this case $V$ decomposes into a direct sum of $G'$-invariant proper subspaces, and elements from $G$ permute them. A group is \emph{primitive} if it is neither intransitive, nor imprimitive. We reserve this terminology for subgroups of $\PGL(V)$ keeping in mind that each such group can be represented by a subgroup of $\GL(V)$.

Let $G'$ be a finite  intransitive subgroup of $\GL(3)$ and $G$ be its image in $\PGL(3)$. Then $G'$ is conjugate to a subgroup $\bbC^*\times \GL(2)$ of block matrices.

To classify such subgroups we have to classify subgroups of $\GL(2)$. We will use the well-known  classification of  finite subgroups of $\PGL(2)$. They are isomorphic to one of the following {\it polyhedral groups}
\begin{itemize}
\item a cyclic group $C_n$;
\item a dihedral group $D_{2n}$ of order $2n\ge 2$;
\item the tetrahedral group $T \cong A_4$ of order 12;
\item the octahedral  group $O\cong S_4$ of order 24;
\item the icosahedral group $I \cong A_5$ of order 60.
\end{itemize}
Two isomorphic subgroups are conjugate subgroups of $\PGL(2)$.

The preimage of such group in $\SL(2,\bbC)$ under the natural map 
$$\SL(2) \to \PSL(2) = \SL(2)/(\pm 1)\cong \PGL(2)$$
 is a  double extension $\overline{G} = 2\pt G$. The group $\overline{G} = 2\pt G$ is called  a {\it binary polyhedral group}. A cyclic group of odd order is isomorphic to a subgroup $\SL(2)$ intersecting trivially the center.

Consider the natural surjective homomorphism of groups
$$\beta:\bbC^*\times \SL(2) \to \GL(2), \ \ (c,A)\mapsto cA.$$
Its kernel is the subgroup $\{(1,I_2), (-1,-I_2)\}$.

Let $G$ be a finite subgroup of $\GL(2)$ with center $Z(G)$. Since $cA = (-c)(-A)$ and $\det(cA) = c^2\det A$, we see that $\widetilde{G} = \beta^{-1}(G)$ is a subgroup of $\mu_{2m}\times \bar{G'}$, where $\overline{G'}$ is a binary polyhedral group with nontrivial center whose image $G'$ in $\PGL(2)$ is isomorphic to $G/Z(G)$. The homomorphism $\beta$ defines an isomorphism from  the kernel $H_2$ of the first projection $\widetilde{G}\to \mu_{2m}$ onto the subgroup
$G_0 = \Ker(\det:G\to \bbC^*)$. Also it  defines an isomorphism from  the kernel $H_1$ of the second projection $\widetilde{G}\to \bar{G'}$  onto $Z(G)$.   Applying Lemma \ref{goursat}, we obtain
$$\widetilde{G} \cong \mu_{2m}\trian_D \overline{G}', \ D = \bar{G'}/G_0.$$

\begin{lemma}\label{bri} Let $G$ be a finite non-abelian subgroup of $\GL(2)$. Then $G = \beta(\widetilde{G})$, where $\widetilde{G}\subset \bbC^*\times \SL(2,\bbC)$ is conjugate to one of the following groups
\begin{itemize}
\item[(i)] $\widetilde{G} =\mu_{2m}\times \bar{I}, \ G\cong m\times \bar{I}$;
\item [(ii)] $\widetilde{G} =  \mu_{2m}\times \overline{O}, \ G\cong m\times \overline{O}$;
\item [(iii)]  $\widetilde{G} =  \mu_{2m}\times \overline{T}, \ G\cong   m\times \overline{T}$;
\item[(iv)] $\widetilde{G} = \mu_{2m}\times Q_{4n},\ G\cong   m\times Q_{4n}$;
\item [(v)] $ \widetilde{G} = \half [\mu_{4m}\times  \overline{O}], \ G\cong 2m\pt O \cong (m\times \overline{T})\pt 2$ (split if $m = 1,2$);
\item [(vi)] $\widetilde{G} = \frac{1}{3}[\mu_{6m}\times\overline{T}], \  G \cong 2m\pt T \cong (m\times 2^2)\pt  3$ (split if $m = 1,3$);
\item [(vii)] $\widetilde{G} \cong \half[\mu_{4m}\times Q_{8n}], G \cong 2m\pt D_{4n}\cong (m\times Q_{4n})\pt 2$ (split if $m = 1,2$);
\item [(viii)] $\widetilde{G} = \half [\mu_{4m} \times Q_{4n}], \   G \cong 2m \pt D_{2n}\cong (m\times 2n)\pt 2 $ (split if $m = 1,2$);
\item [(ix)] $\widetilde{G} = \frac{1}{4}[\mu_{4m}\times Q_{4n}], n $ is odd,\   $G\cong m\pt D_{2n} \cong (m\times n)\pt 2$ (split if $m = 1,2$).

\end{itemize}
\end{lemma}
Note that, although $Q_{8n}$ has two different non-cyclic subgroups of index 2, they are conjugate under  an element of $\SL(2)$, so they lead to conjugate subgroups in $\GL(2)$.

Lemma \ref{coh} gives us some information when some of these extensions split.

An abelian subgroup $G\subset \GL(2)$ is conjugate to a
subgroup of diagonal matrices  of the form $(\epsilon_{m}^a,
\epsilon_{n}^b)$, where $\epsilon_m,\epsilon_n$ are primitive roots
of unity and $a,b\in \bbZ$. Let $d = (m,n), m = du, n = dv, d = kq$
for some fixed positive integer $k$. Let $H_1 = \la \epsilon_m^k\ra
\subset \la \epsilon_m\ra, H_2 = \la \epsilon_n^k\ra \subset \la
\epsilon_n\ra$ be cyclic subgroups of index $k$. Applying Lemma
\ref{goursat} we obtain
$$G\cong \la \epsilon_m\ra\trian_k\la \epsilon_n\ra$$
where the homomorphisms $\la \epsilon_m\ra \to k, \la \epsilon_n\ra
\to k$  differ by  an automorphism of the cyclic group $ \la
\epsilon_k\ra \cong k$ defined by a choice of a new generator
$\epsilon_m^s, (s,k) = 1$. In this case
\begin{equation}\label{abgr}
G = (\la \epsilon_m^k\ra\times \la \epsilon_n^k\ra_\bullet \la
\epsilon_k\ra
\end{equation}
is of order
$mn/k = uvkq^2$. In other words, $G$ consists of diagonal matrices
of the form $(\epsilon_m^a,\epsilon_n^b)$, where $a\equiv sb \mod
k$.

\begin{corollary} Let $G$ be an intransitive finite subgroup of $\GL(3)$. Then its image in $\PGL(3)$ consists of transformations $[ax_0+bx_1,cx_0+dx_1,x_2]$,  where the matrices $\left(\begin{smallmatrix}a&b\\
c&d\end{smallmatrix}\right)$ form a non-abelian finite subgroup $H$ of $\GL(2)$ from Lemma \ref{bri} or an abelian group of the form \eqref{abgr}.
\end{corollary}

Now suppose $G$ is transitive but imprimitive subgroup of
$\PGL(3)$.  Let $G'$ be its largest intransitive normal
subgroup. Then $G/G'$ permutes transitively the invariant subspaces
of $G'$, hence we may assume that all of them are one-dimensional.
Replacing $G$ by a conjugate group we may assume that $G'$ is a
subgroup of diagonal matrices. We can represent its elements by
diagonal matrices  $g = (\epsilon_m^a,\epsilon_n^b,1)$, where
$a\equiv sb \mod k$ as in \eqref{abgr}. The group  $G$ contains  a
cyclic permutation $\tau$  of  coordinates. Since $G'$ is a normal
subgroup of $G$, we get $\tau^{-1}g\tau =
(\epsilon_n^{-b},\epsilon_n^{-b}\epsilon_m^a,1)\in G'$. This implies
that $n|bm, m|an$, hence $u|b, v|a$. Since
$(\epsilon_m,\epsilon_n^s,1)$ or $(\epsilon_m^{s'},\epsilon_n,1),
ss'\equiv 1 \mod k$, belongs to $G$ we must have $u=v = 1$, i.e. $m
= n = d$. Therefore $G'$ consists of diagonal matrices $g =
(\epsilon_d^a,\epsilon_d^{sa},1)$. Since $\tau^{-1}g\tau =
(\epsilon_d^{-sa},\epsilon_d^{a-sa},1)\in G'$, we get $a-sa\equiv
-s^2a \mod k$ for all $a\in \bbZ/m\bbZ$. Hence the integers $s$
satisfy the congruence $s^2-s+1\equiv 0 \mod k$. If, moreover,
$G/G'\cong S_3$, then we have an additional condition $s^2\equiv
1\mod k$, and hence either $k = 1$ and $G'= \mu_n\times \mu_n$ or
$k=3$, $s=2$ and $G'=n \times n/k$.

This gives the following.

\begin{theorem}\label{imprim} Let $G$ be a transitive imprimitive finite subgroup of $\PGL(3)$. Then $G$ is conjugate to one of  the following groups
\begin{itemize}
\item $G \cong n^2:3$ generated by transformations
$$[\epsilon_nx_0,x_1,x_2],\ [x_0,\epsilon_nx_1,x_2],\
[x_2,x_0,x_1];$$
\item $G \cong n^2:S_3$ generated by transformations
$$[\epsilon_nx_0,x_1,x_2],\  [x_0,\epsilon_nx_1,x_2],\  [x_0,x_2,x_1],\  [x_2,x_0,x_1]; $$
\item $G = G_{n,k,s}\cong (n\times \frac{n}{k}):3,$ where $k> 1, k|n$ and $s^2-s+1 = 0\mod k$. It is generated by  transformations
$$[\epsilon_{n/k}x_0,x_1,x_2], \ [\epsilon_{n}^sx_0,\epsilon_nx_1,x_2],\  [x_2,x_0,x_1].$$
\item $G \simeq (n \times \frac{n}{3}):S_{3}$ generated by transformations
$$[\epsilon_{n/3}x_{0},x_{1},x_{2}],\
[\epsilon_{n}^{2}x_{0},\epsilon_{n}x_{1},x_{2}],\
[x_{0},x_{2},x_{1}],\ [x_{1},x_{0},x_{2}].$$
\end{itemize}
\end{theorem}

The next theorem is a well-known result of Blichfeldt \cite{Blichfeldt}.

\begin{theorem} Any primitive finite subgroup $G$ of $\PGL(3)$ is conjugate to one of the following groups.
\begin{enumerate}
\item The icosahedral group $A_5$ isomorphic to $L_2(5)$. It leaves invariant a nonsingular conic.
\item The \emph{Hessian group}  of order 216 isomorphic to $3^2:\overline{T}$. It is realized as the group of automorphisms of the Hesse pencil of cubics
$$x^3+y^3+z^3 + txyz = 0.$$ 
\item The \emph{Klein group}  of order 168 isomorphic to $L_2(7)$ (realized as the full group of automorphisms of the Klein quartic $x^3y+y^3z+z^3x  = 0$).
\item The \emph{Valentiner group} of order 360 isomorphic to $A_6$. It can be realized as the full group of automorphisms of the nonsingular plane sextic
$$10x^3y^3+9zx^5+y^5-45x^2y^2z^2-135xyz^4+27z^6 = 0.$$
\item Subgroups of the Hessian group:
\begin{itemize}
\item $3^2:4$;
\item $3^2:Q_8$.
\end{itemize}

\end{enumerate}
\end{theorem}

\subsection{Finite groups of automorphisms of $\bfF_0$}
Since $\bfF_0$ is isomorphic to a nonsingular quadric in $\bbP^3$, the group $\Aut(\bfF_0)$ is isomorphic to the projective orthogonal group $\PO(4)$. The classification of finite subgroups of $O(4)$ is due to \'E. Goursat \cite{Goursat} (in the real case see a modern account in \cite{Conway2}). Goursat'  Lemma   \ref{goursat} plays an important role in this classification.

Obviously,
$$\Aut(\bfF_0)  \cong\PGL(2)\wr S_2,$$

First we classify subgroups of $\PGL(2)\times \PGL(2)$ by  applying Goursat's Lemma.

Observe the following  special subgroups of $\PGL(2)\times \PGL(2)$.
\begin{enumerate}
\item  $G =  G_1\times G_2$ is the {\it product subgroup}.
\item
$G\trian_1 G =  \{(g_1,g_2)\in G\times G: \alpha(g_1) = g_2\} \cong G$
 is a {\it $\alpha$-twisted  diagonal subgroup}. If $\alpha = \text{id}_G$, we get  the {\it diagonal subgroup}.
\end{enumerate}
Note that  $\alpha$-twisted diagonal groups are conjugate in $\Aut(\bfF_0)$ if $\alpha(g) = xgx^{-1}$for some $x$ in the normalizer of $G$ inside $\Aut(\bbP^1)$. In particular, we may always assume that $\alpha$ is an exterior automorphism of $G$.

We will use the notation $[p_1,\ldots,p_r]$ for the Coxeter group defined by the Coxeter diagram

\begin{figure}[h]
\xy @={(0,0),(10,0),(20,0),(40,0),(50,0)}@@{*{\bullet}};
(0,0)*{};(25,0)**{}**\dir{-};
(35,0)*{};(50,0)**{}**\dir{-};
(-35,0)*{};(5,3)*{p_1};
(15,3)*{p_2};
(45,3)*{p_r};
(30,0)*{\cdots};
\endxy
\end{figure}

Following \cite{Conway2} we write $[p_1,\ldots,p_r]^+$ to denote the index 2 subgroup of even length words in standard generators of the Coxeter group.  If exactly one of the numbers $p_1,\ldots,p_r$ is even, say $p_k$, there are two other subgroups of index 2  denoted by $[p_1,\ldots,p_r^+]$ (resp.  $[{}^+p_1,\ldots,p_r])$. They consist of words which contain each generator  $R_1,\ldots,R_{k-1}$ (resp. $R_{k+1},\ldots,R_{r}$) even number of times. The intersection of these two subgroups is denoted by $[{}^+p_1,\ldots,p_r^+]$. For example,
$$D_{2n} = [n], \ T = [3,3]^+,\ O = [3,4]^+,\ I = [3,5]^+.$$

Recall that each group $[p_1,\ldots,p_r]$ has  a natural linear representation  in $\bbR^{r}$ as a reflection group. If $r= 3$, the corresponding representation defines a subgroup of $\PO(4)$. If $r = 2$, it defines a subgroup of $\PO(3)$ which acts diagonally on $\bbP^2\times \bbP^2$ and on $\bbP^1\times \bbP^1$ embedded in $\bbP^2\times \bbP^2$ by the product of the Veronese maps. We denote $\overline{[p_1,\ldots,p_r]}$ the quotient of $[p_1,\ldots,p_r]$ by its center. Similar notation is used for the even subgroups of $[p_1,\ldots,p_r]$.

\begin{theorem}\label{clas} Let $G$ be a finite subgroup of $\PGL(2)\times \PGL(2)$ not conjugate to the product $A\times B$ of subgroups of $\PGL(2)$. Then $G$ is conjugate to one of  the following groups or its image under the switching of the factors.
\begin{itemize}
\item $\frac{1}{60}[I\times I] \cong I \cong [3,5]^+$;
\item $\frac{1}{60}[I\times I] \cong I \cong [3,3,3]^+$;
\item $\frac{1}{24}[O \times O] \cong O \cong [3,4]^+$;
\item $\frac{1}{24}[O \times O] \cong O \cong [2,3,3]^+$;
\item $\frac{1}{12}[T \times T] \cong T \cong [3,3]^+$;
\item $\half[O\times O]\cong (T\times T):2 \cong \overline{[3,4,3]^+}$;
\item $\frac{1}{6}[O\times O] \cong 2^4:S_3 \cong \overline{[3,3,4]^+}$;
\item $\frac{1}{3}[T\times T] \cong 2^4:3 \cong \overline{[{}^+3,3,4^+]};$
\item $\half [D_{2m}\times D_{4n}] \cong (m\times D_{2n})^\bullet 2 $ ($m,n\ge 2$);
\item $\frac{1}{4}[D_{4m}\times D_{4n}] \cong (m\times n): 4$ ($m,n$  odd);
\item $\frac{1}{2k}[D_{2mk}\times D_{2nk}]_s \cong (m\times n): D_{2k}, \ (s,k) = 1$);
\item $\frac{1}{2k}[D_{2mk}\times D_{2nk}]_s \cong (m\times n): D_{2k}, \ (s,2k) = 1, m,n$  odd);
\item $\frac{1}{k}[C_{mk}\times C_{nk}]_s\cong (m\times n) \pt k \ ((s,k) = 1)$;
\item $\frac{1}{k}[C_{mk}\times C_{nk}]_s\cong (m\times n) \pt k$ \ ($(s,2k) = 1, m,n$  odd;);
\item $\half [D_{2m}\times O] \cong (m\times T):2$;
\item $\half [D_{4m}\times O] \cong (D_{2m}\times T): 2  (m\ge 2$);
\item $\frac{1}{6} [D_{6n}\times O] \cong (m\times 2^2):S_3  (m\ge 2$);
\item $\half [C_{2m}\times O] \cong (m\times T)\pt 2$  (split if $m= 1$);
\item $\frac{1}{3}[C_{3m}\times T] \cong (m\times 2^2)\pt 3 $ (split if $m = 1$);
\item $\half [D_{4m}\times D_{4n}] \cong (D_{2m}\times D_{2n})^\bullet 2$ ($m,n\ge 2$);
\item $\half [C_{2m}\times D_{4n}] \cong (m\times D_{2n})^\bullet  2$  ($n\ge 2$);
\item $\half [C_{2m}\times D_{2n}] \cong (m\times n):2 \cong m\pt D_{2n}$.

\end{itemize}
\end{theorem}

All other finite subgroups of $\Aut(\bbP^1\times \bbP^1)$  are conjugate to a group
$G^0\pt 2$, where the quotient  $2$ is represented by an automorphism which interchanges the two rulings of $\bfF_0$. It is equal to $\tau\circ g$, where $\tau$ is the switch $(x,y)\mapsto (y,x)$ and $g\in \PGL(2)\times \PGL(2)$. 

\bigskip
\subsection{Finite groups of automorphisms of $\bfF_n, n \ne 0$}
 Let $S$ be a minimal ruled surface $\bfF_n, n \ne 0$. If $n = 1$, the group $\Aut(\bfF_1)$ is isomorphic to a subgroup of $\Aut(\bbP^2)$ leaving one point fixed. We will not be interested in such subgroups so we assume that $n \ge 2$.

\begin{theorem}\label{thh} Let $S = \bfF_n, n \ne 0$. We have
$$\Aut(\bfF_n) \cong \bbC^{n+1}:(\GL(2)/\mu_n),$$
where $\GL(2)/\mu_n$ acts on $\bbC^{n+1}$ by means of its natural linear representation in the space of binary forms of degree $n$.  Moreover,
$$\GL(2)/\mu_n
 \cong \begin{cases}\bbC^*:\PSL(2),& \text{if $n$ is even},\\
 \bbC^*:\SL(2),& \text{if $n$ is odd}.
\end{cases}
 $$
\end{theorem}

\begin{proof} This is of course well-known. We identify $\bbF_n$ with the weighted projective plane $\bbP(1,1,n)$.
An automorphism is given by  the formula
$$(t_0,t_1,t_2) \mapsto (at_0+bt_1,ct_0+dt_1,et_2+f_n(t_0,t_1)),$$ where $f_n$ is a homogeneous polynomial of degree $n$. The vector space   $\bbC^{n+1}$ is identified with the normal subgroup of transformations $[t_0,t_1,t_2+f_n(t_0,t_1)]$. The quotient by this subgroup is isomorphic to the subgroup of transformations $[at_0+bt_1,ct_0+dt_1,et_2]$ modulo  transformations of the form $[\lambda t_0,\lambda t_1,\lambda^nt_2]$.  This  group  is obviously isomorphic to $\GL(2)/\mu_n$. Consider the natural projection $\GL(2)/\mu_n\to \PGL(2) \cong \PSL(2)$. Define a homomorphism
$\SL(2)\to \GL(2)/\mu_n$ by assigning to a matrix $A$ the coset of $A$ modulo $\mu_n$. If $n$ is even, the kernel of this homomorphism is $\la -I_2\ra$, so we have a splitting $\GL(2) \cong \bbC^*:\PGL(2)$. If $n$ is odd, the homomorphism is injective and defines a splitting  $\GL(2)/\mu_n \cong \bbC^*:\SL(2)$.
\end{proof}

Suppose $G$ is a finite subgroup of $\Aut(\bbF_n)$. Obviously, $G$ is contained in the subgroup
$\GL(2)/\mu_n$.

Suppose $G\cap \bbC^* = \{1\}$. Then $G$ is isomorphic to a subgroup of $\PGL(2)$ (resp. $\SL(2))$ over which the extension splits. Note that the kernel $\bbC^*$ of the projections $\GL(2)/\mu_n\to \PGL(2)$ or $\GL(2)/\mu_n\to \SL(2)$ is the center. Thus each finite subgroup $H$ of $\PGL(2)$ (resp. $\SL(2))$ defines $k$  conjugacy classes of subgroups isomorphic to $H$, where $k = \#\Hom(H,\bbC^*) = \#G/[G,G]$.

If $G\cap \bbC^* \cong \mu_m$ is non-trivial, the group is a central extension $m\pt H$, where $H$ is a polyhedral group, if $n$ is even, and a binary polyhedral group otherwise. We can apply Lemma \ref{coh} to find some cases when the extension must split. In other cases the structure of the group is determined by using
Theorem  \ref{clas}. We leave this to  the reader.

\section{Automorphisms of conic bundles}\label{sec4}
\subsection{Geometry of conic bundles}
Let $\phi:S\to \bbP^1$ be a conic  bundle with singular fibres over points in a finite set $\Sigma \subset \bbP^1$. We assume that $k = \#\Sigma > 0$.
  Recall that each singular fibre $F_x,x\in \Sigma,$  is the bouquet of two $\bbP^1$'s.

 Let $E$ be a section of the conic bundle fibration $\phi$. The Picard group of $S$ is freely generated by the divisor classes of $E$, the class $F$ of a fibre, and the classes of $k$  components of singular fibres, no two in the same fibre. The next lemma follows easily from the intersection theory on $S$.

 \begin{lemma}\label{ineqq} Let $E$ and $E'$ be two sections with negative self-intersection $-n$. Let $r$ be the number of components of singular fibres which intersect both $E$ and $E'$. Then $k-r$ is even and
$$2E\cdot E' =k-2n-r.$$
In particular,
$$k \ge 2n+r .$$
\end{lemma}

Since a conic bundle $S$  is isomorphic to a blowup of
a minimal ruled surface, it  always contains a section $E$ with
negative self-intersection $-n$. If $n \ge 2$, we obviously get $k\ge 4$. If $n = 1$, since $(S,G)$ is minimal, there exists $g\in G$ such that $g(E)\ne  E$ and $E\cap g(E) \ne \emptyset$. Applying  the previous lemma we get
$$k \ge 4.$$

\subsection{Exceptional conic bundles}\label{exc} We give three different constructions of the same conic bundle, which we will call an \emph{exceptional conic bundle}.

\smallskip\noindent
\textit{First construction.}

 Choose a ruling $p:\bfF_0\to \bbP^1$ on $\bfF_0$ and fix two points on the base, say $0$ and $\infty$. Let $F_0$ and $F_\infty$ be the corresponding fibres. Take $g+1$ points $a_1,\ldots,a_{g+1}$ on $F_0$ and $g+1$ points $a_{g+2},\ldots,a_{2g+2}$ on $F_\infty$ such that no two lie in the same fibre of the second ruling $q:\bfF_0\to \bbP^1$. Let $\sigma:S\to \bfF_0$ be the blowup of the points $a_1,\ldots,a_{2g+2}$. The composition $\pi=q\circ \sigma: S\to \bbP^1$ is a conic bundle with $2g+2$ singular fibres $R_i+R_i'$ over the points $x_i = q(a_i), i = 1,\ldots,2g+2$. For $i=1,\ldots,g+1$,
$R_i =\sigma^{-1}(a_i)$ and $R_{n+i}$ is  the proper
transform of the fibre $q^{-1}(a_i)$. Similarly, for $i=1,\ldots,n$,
$R_i' $ is  the proper transform of the fibre $q^{-1}(a_i)$
and $R_{g+1+i}' =\sigma^{-1}(a_{g+1+i}).$

Let $E_0,E_\infty$ be the proper transforms of
$F_0,F_\infty$ on $S$. Each is  a section of the conic bundle
$\pi$. The section $E_0$ intersects $R_1,\ldots,R_{2g+2},$ and the
section $E_\infty$ intersects $R_1',\ldots,R_{2g+2}'.$

Let
$$D_0 = 2E_0+\sum_{i=1}^{2g+2}R_i,  \quad D_\infty =   2E_\infty+\sum_{i=1}^{2g+2}R_i'.$$
It is easy to check that $D_0\sim D_\infty$.  Consider the pencil
$\calP$  spanned by the curves $D_0$ and $D_\infty$. It has $2g+2$
simple  base points $p_i = R_i\cap R_i'$. Its general member is a
nonsingular curve $C$. In fact, a standard  formula for computing the Euler characteristic of a fibred surface in terms of the Euler characteristics of fibres shows
that all members except $D_0$ and $D_\infty$ are nonsingular curves.
Let $F$ be a fibre of the conic bundle. Since $C\cdot F = 2$, the
linear system $|F|$ cuts out a $g_2^1$ on $C$, so it is a hyperelliptic
curve or the genus $g$ of $C$ is 0 or 1. The  points $p_i$ are
obviously the ramification points   of the $g_2^1$. Computing the
genus of $C$ we find that it is equal  to $g$, thus
$p_1,\ldots,p_{2g+2}$ is the set of ramification points. Obviously all
nonsingular members  are isomorphic curves. Let $\sigma:S'\to S$ be
the blowup the base points $p_1,\ldots,p_{2g+2}$ and let $\overline{D}$
denote the proper transform of a curve on $S$. We have
$$2\bar{E}_0+2\bar{E}_\infty+\sum_{i=1}^{2g+2}(\bar{R}_i+\bar{R}_i' + 2\sigma^{-1}(p_i))\sim 2\sigma^*(C).$$
This shows that there exists a double over $X'\to S'$ branched along the divisor  $\sum_{i=1}^{2g+2}(\bar{R}_i+\bar{R}_i' )$. Since $\bar{R_i}^2 = \bar{R}_i'{}^2 = -2$, the ramification divisor on $X'$ consists of $4g+4$ $(-1)$-curves. Blowing them down we obtain a surface $X$ isomorphic to the product $C\times \bbP^1$. This gives us the following.

\medskip\noindent
\textit{Second construction. } A pair $(C,h)$ consisting of a nonsingular curve and an involution $h\in \Aut(C)$ with quotient $\bbP^1$ will be called a \emph{hyperelliptic curve}. If $C$ is of genus $g\ge 2$, then $C$ is a hyperelliptic curve and $h$ is its involution defined by the unique $g_2^1$ on $C$. Let $\delta$ be an involution of
$\bbP^1$ defined by $(t_0,t_1)\mapsto (t_0,-t_1)$. Consider the
involution $\tau = h\times \delta$ of the product $X= C\times
\bbP^1$. Its fixed points are $4g+4$ points  $c_i\times\{0\}$ and
$c_i\times \{\infty\}$, where $X^{\la h\ra} = \{c_1,\ldots,c_{2g+2}\}$.
Let $X'$ be a minimal resolution of $X/(\tau)$. It is easy
to see that the images of the curves $\{c_i\}\times \bbP^1$ are
$(-1)$-curves on $X'$. Blowing them down we obtain our exceptional
conic bundle.

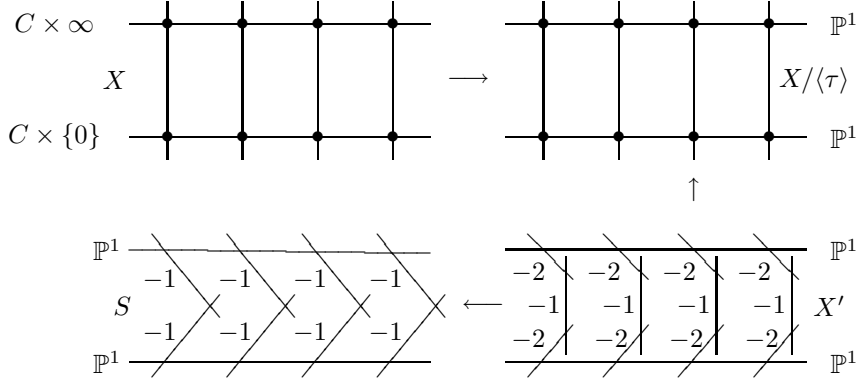
\begin{figure}[h]
\xy @={(0,0),(10,0),(20,0),(30,0),(50,0),(60,0),(70,0),(80,0)}@@{*{\bullet}};
@={(0,-15),(10,-15),(20,-15),(30,-15),(50,-15),(60,-15),(70,-15),(80,-15)}@@{*{\bullet}};
(-5,0)*{};(35,0)**{}**\dir{-};
(45,0)*{};(85,0)**{}**\dir{-};
(-5,-15)*{};(35,-15)**{}**\dir{-};
(45,-15)*{};(85,-15)**{}**\dir{-};
(0,3)*{};(0,-18)**\dir{-};
(10,3)*{};(10,-18)**\dir{-};
(20,3)*{};(20,-18)**\dir{-};
(30,3)*{};(30,-18)**\dir{-};
(50,3)*{};(50,-18)**\dir{-};
(60,3)*{};(60,-18)**\dir{-};
(70,3)*{};(70,-18)**\dir{-};
(80,3)*{};(80,-18)**\dir{-};
(-5,-30)*{};(35,-30.5)**{}**\dir{-};
(45,-30)*{};(85,-30)**{}**\dir{-};
(-5,-45)*{};(35,-45)**{}**\dir{-};
(45,-45)*{};(85,-45)**{}**\dir{-};
(-2,-28)*{};(7,-39)**{}**\dir{-};
(-2,-47)*{};(7,-36)**{}**\dir{-};
(8,-28)*{};(17,-39)**{}**\dir{-};
(8,-47)*{};(17,-36)**{}**\dir{-};
(18,-28)*{};(27,-39)**{}**\dir{-};
(18,-47)*{};(27,-36)**{}**\dir{-};
(28,-28)*{};(37,-39)**{}**\dir{-};
(28,-47)*{};(37,-36)**{}**\dir{-};
(48,-28)*{};(54,-34)**{}**\dir{-};
(48,-47)*{};(54,-40)**{}**\dir{-};
(53,-31)*{};(53,-44)**{}**\dir{-};
(58,-28)*{};(64,-34)**{}**\dir{-};
(58,-47)*{};(64,-40)**{}**\dir{-};
(63,-31)*{};(63,-44)**{}**\dir{-};
(68,-28)*{};(74,-34)**{}**\dir{-};
(68,-47)*{};(74,-40)**{}**\dir{-};
(73,-31)*{};(73,-44)**{}**\dir{-};
(78,-28)*{};(84,-34)**{}**\dir{-};
(78,-47)*{};(84,-40)**{}**\dir{-};
(83,-31)*{};(83,-44)**{}**\dir{-};
(-15,0)*{C\times \infty};
(-15,-15)*{C\times \{0\}};
(-7,-7.5)*{X};
(90,0)*{\bbP^1};
(90,-15)*{\bbP^1};
(-8,-30)*{\bbP^1};
(-6,-37.5)*{S};
(-8,-45)*{\bbP^1};
(90,-30)*{\bbP^1};
(90,-45)*{\bbP^1};
(88,-37.5)*{X'};
(86,-7.5)*{X/\la \tau\ra};
(-1,-34)*{-1};(-1,-41)*{-1};
(9,-34)*{-1};(9,-41)*{-1};
(19,-34)*{-1};(19,-41)*{-1};
(29,-34)*{-1};(29,-41)*{-1};
(48,-33)*{-2};(48,-42)*{-2};
(58,-33)*{-2};(59,-42)*{-2};
(68,-33)*{-2};(69,-42)*{-2};
(78,-33)*{-2};(79,-42)*{-2};
(50,-37.5)*{-1};
(60,-37.5)*{-1};(70,-37.5)*{-1};
(80,-37.5)*{-1};
(40,-7.5)*{\longrightarrow};
(42,-37.5)*{\longleftarrow};
(70,-22)*{\uparrow};
\endxy
\end{figure}

\textit{Third construction.}

Let us consider  a quasi-smooth hypersurface $Y$  of degree $2g+2$ in weighted projective space $\bbP = \bbP(1,1,g+1,g+1)$ given by an equation
\begin{equation}\label{hypeq}
F_{2g+2}(T_0,T_1)+T_2T_3 = 0,
\end{equation}
where $F_{2g+2}(T_0,T_1)$ is a homogeneous polynomial of degree
$2g+2$ without multiple roots. The surface is a double cover of
$\bbP(1,1,g+1)$ (the cone over a Veronese curve of degree $g+1$)
branched over the curve  $F_{2g+2}(T_0,T_1)+T_2^2 = 0$. The
preimages of the singular point of $\bbP(1,1,g+1)$ with coordinates 
$(0,0,1)$ is a pair of  singular points of  $Y$ with coordinates
$(0,0,1,0)$ and $(0,0,0,1)$. The singularities are locally
isomorphic to the singular points of a cone of the Veronese surface
of degree $g+1$. Let $S$ be a minimal resolution of $Y$. The
preimages of the singular points are disjoint smooth rational
curves $E$ and $E'$ with self-intersection $-(g+1)$. The projection
$\bbP(1,1,g+1,g+1)\to \bbP^1, (t_0,t_1,t_2,t_3)\mapsto (t_0,t_1)$
lifts to a conic bundle on $S$ with sections $E,E'$. The pencil
$\lambda T_2+\mu T_3 = 0$ cuts out a pencil of curves on $Y$ which
lifts to a pencil of bi-sections of the conic bundle $S$ with $2g+2$
base points $(t_0,t_1,0,0)$, where $F_{2g+2}(t_0,t_1) = 0$.

It is easy to see that this is a general example of an exceptional
conic bundle. In Construction 2, we blow down the sections
$E_0,E_\infty$ to singular points. Then consider an involution
$g_0$ of the surface which is a descent of the automorphism of the
product $C\times \bbP^1$ given by $\id_C\times \psi$, where
$\psi:(t_0,t_1)\mapsto (t_1,t_0)$. The quotient by $(g_0)$
gives $\bbP(1,1,g+1)$ and the ramification divisor is the image on
$S$  of the curve $C\times (1,1)$ or $C\times (1,-1)$. On one of
these curves $g_0$ acts identically, on the other one it acts as the
involution defined by the $g_2^1$.

\begin{proposition}\label{cdp} Let $\phi:S\to \bbP^1$ be a minimal conic $G$-bundle with $k\le 5$ singular fibres. Then $S$ is a Del Pezzo surface, unless $k = 4$ and $S$ is an exceptional conic bundle.
\end{proposition}

\begin{proof} Since $k\le 5$ we have $K_S^2 = 8-k\ge 3$. By Riemann-Roch, $|-K_S|\ne \emptyset$. Suppose $S$ is not a Del Pezzo surface. Then there exists an irreducible  curve $C$ such that $-K_S\cdot C \le 0$. Suppose, the equality takes place. By Hodge's Index Theorem, $C^2 < 0$, and by the adjunction formula, $C^2 = -2$ and $C\cong \bbP^1$. If   strict inequality takes place, then  $C$ is a component of a divisor $D\in |-K_S|$, hence $|-K_S-C|\ne \emptyset$ and $|K_S+C|=\emptyset$. Moreover, since $K_S^2 > 0$, we have $C\not\in |-K_S|$. Applying Riemann-Roch to the divisor $K_S+C$ we easily obtain that $C$ is of arithmetic genus 0, and hence $C\cong \bbP^1$. By adjunction, $C^2 \le -2$.  In both cases we have a smooth rational curve with $C^2\le -2$.

If $k = 4$ and $S$ is an exceptional conic bundle, then $S$ is not  a Del Pezzo surface since it has  sections with self-intersection $-2$.  Assume this is not the case.  Let $C$ be the union of smooth rational curves with self-intersection $< -2$. It is obviously a $G$-invariant curve, so we can write $C\sim -aK_S-bf$, where $f$ is the divisor class of a fibre of $\phi$. Intersecting with $f$ we get $a> 0$. Intersecting with $K_S$, we get $2b > ad$, where $d = 8-k \ge 3$. It follows from Lemma \ref{ineqq}, that $S$ contains a section $E$ with self-intersection $-2$ or $-1$. Intersecting $C$ with $E$ we get $0\le C\cdot E = a(-K_S\cdot E)-b\le a -b$.  This contradicts the previous inequality. Now let us take $C$ to be the union of $(-2)$-curves. Similarly, we get $2b = ad$ and
$C^2 = -aK_S\cdot C-bC\cdot f = -bC\cdot f = -2ab$. Let $r$ be the number of irreducible components of $C$. We have $2a = C\cdot f \ge r$ and
$-2r\le C^2 = -2ab  \le -br$. If $b = 2$, we have the equality everywhere, hence $C$ consists of $r = 2a$ disjoint sections, and $8 = rd$. Since $d\ge 3$, the only solution is $d = 4, r = 2$, and this leads to the exceptional conic bundle. Assume $b = 1$. Since $C^2= -2a$ is even, $a$ is a positive integer, and we get
$2= ad$. Since $d\ge 3$, this is impossible.
\end{proof}

\smallskip
\subsection{Automorphisms of an exceptional conic bundle} Let us describe the automorphism group of an exceptional conic bundle. The easiest way to do it to use Construction 3. We denote by $Y_g$ an exceptional conic bundle given by equation \eqref{hypeq}. Since we are interested only in minimal groups we assume that $g\ge 1$.

 Since $K_{Y_g} = \calO_{\bbP}(-2)$, any automorphism $\sigma$ of $Y_g$ is a restriction of an automorphism of $\bbP$. Let $G_1$ be the subgroup of $\SL(2)$ of transformations preserving the zero divisor of $F_{2g+2}(T_0,T_1)$  and $\chi_1:G_1\to \bbC^*$ be the multiplicative character of $G_1$ defined by $\sigma_1^*(F_{2g+2}) = \chi_1(\sigma_1)F_{2g+2}$.  Similarly, let $G_2$ be the subgroup of  $\GL(2)$ of matrices preserving the zeroes of $T_2T_3$ and let $\chi_2:G_2\to \bbC^*$ be the character defined by $\sigma^*(T_2T_3) = \chi_2(\sigma_2)T_2T_3$.
Let
$$(G_1\times G_2)^0 = \{(\sigma_1,\sigma_2)\in G_1\times G_2:\chi_1(\sigma_1) = \chi_2(\sigma_2)\}.$$
In the notation of the diagonal products,
$$(G_1\times G_2)^0 = \frac{1}{m}[G_1\times G_2],$$
where $\chi_1(G_1) = \mu_m \subset \bbC^*$. The subgroup
$$K = \la (-I_2,(-1)^{g+1}I_2)\ra $$
acts identically on $Y_g$ and the quotient group is isomorphic to $\Aut(Y_g)$.

Let $\Aut(Y_g)  \cong (G_1\times G_2)^0/K \to \PGL(2)$ be the homomorphism induced by the projection of $G_1$ to $\PGL(2)$. Its image is a finite subgroup $P$ of $\PGL(2)$.  Its kernel $H$ consists of cosets modulo $K$ of pairs $(\pm I_2,\sigma_2)$, where $\chi_2(\sigma_2) = 1$. Clearly, $H\cong \Ker(\chi_2)$.

It is easy to see that $\Ker(\chi_2) \cong \bbC^*:2$ is generated by diagonal matrices with determinant 1 and the matrix which switches the coordinates. Inside of $\GL(2)$ it is conjugate to the normalizer $N$ of the maximal torus in $\SL(2)$.  So we obtain an isomorphism
\begin{equation}\label{extt}
\text{Aut$(Y_g) \cong N \pt P$}.
\end{equation}
Suppose there exists a homomorphism $\eta:G_1\to \bbC^*$ such that $\eta(-1) = (-1)^{g+1}$.  Then the homomorphism
$$G_1\to (G_1\times G_2)^0/K,\  \sigma_1 \mapsto (\sigma_1,\eta(\sigma_1)I_2) \mod K$$
factors through a homomorphism $P\to (G_1\times G_2)^0/K$ and defines a splitting of the extension \eqref{extt}. Since elements of the form $(\sigma_1,\eta(\sigma_1)I_2)$ commute with elements of $N$, we see that the extension is trivial when it splits. It is easy to see that the converse is also true. Since the trivial $\eta$ works when $g$ is odd, we obtain that the extension always splits in this case. Assume $g$ is even and $G_1$ admits  a 1-dimensional representation $\eta$ with $\eta(-I_2) = -1$.  If its kernel is trivial, $G_1$ is isomorphic to a subgroup of $\bbC^*$, hence cyclic. Otherwise, the kernel is a subgroup of $\SL(2)$ not containing the center. It must be a cyclic subgroup of odd order. The image is a cyclic group. Thus $G_1$ is either cyclic, or a binary dihedral group $D_{2n}$ with $n$ odd.

To summarize we have proved the following.

\begin{proposition}\label{autoex} The group of automorphisms of an exceptional conic bundle \eqref{hypeq} is isomorphic to an extension $N\pt P$, where $P$ is the subgroup of $\PGL(2)$ leaving the set of zeroes of $F_{2g+2}(T_0,T_1)$ invariant and $N \cong \bbC^*:2$ is a group of matrices with determinant $\pm 1$ leaving $T_2T_3$ invariant.  Moreover, the extension splits and defines an isomorphism
$$\Aut(Y_g) \cong N\times P$$
if and only if $g$ is odd, or $g$ is even and $P$ is either a cyclic group or a dihedral group $D_{4k+2}$.
\end{proposition}

 Now let $G$ be a finite minimal subgroup of $\Aut(Y_g)$. Assume first that $\Aut(Y_g) \cong N\times P$. Let $N'$ be the projection of $G$ to $N$ and $P'$ be the projection to $P$. Since $G$ is minimal, $N'$ contains an element which switches $V(T_2)$ and $V(T_3)$. Thus $N'$ is isomorphic to a dihedral group $D_{2n}$. Applying Goursat's Lemma we obtain that
 $$G\cong N'\trian_DQ,$$
where $D$ is a common quotient of $N'$ and $Q$. If $N' $ is a dihedral group, then $D$  is either dihedral group  or  a cyclic of order 2. Using Goursat's Lemma  it is easy to list all possible subgroups. We leave it as an exercise to the reader.

If $\Aut(Y_g)$ is not isomorphic to the direct product $N\times P$, we can only say that
$$\text{$G \cong H\pt Q$},$$
where $H$ is a subgroup of $D_{2n}$ or $Q_{4n}$, and $Q$ is a polyhedral group. Note that we can write these extensions in the form
$ n\pt (2\pt Q)$ or $n\pt (2^2\pt Q)$.

\begin{example}\label{exept} Let $\phi:S\to \bbP^1$ be an exceptional conic bundle with $g = 1$. It has 4 singular fibres. According to the first construction, the blow up $S'$ of $S$ at the four singular points of the singular fibres  admits an elliptic fibration  $f:S'\to \bbP^1$ with two singular fibres of type $I_0^*$ in Kodaira's notation. The $j$-invariant of the fibration is zero, and after the degree 2 base change $\bbP^1\to \bbP^1$ ramified at 2 points, the surface becomes isomorphic to the product $E\times \bbP^1$, where $E$ is an elliptic curve. Its $j$-invariant corresponds to the cross ratio of the 4 points, where a section of $\pi$ with self-intersection $-2$ intersects the singular fibres. Conversely, starting from the product, we can divide it by an elliptic involution, to get the conic bundle. This is our second construction.

According to the third construction, the surface can be given by an equation
$$F_4(T_0,T_1) +T_2T_3 = 0$$
in the weighted projective space $\bbP(1,1,2,2)$. The projection to $(T_0,T_1)$ is a rational map undefined at the four points $P_i = (a,b,0,0)$, where $F_4(a,b) = 0$. After we blow them up we get the conic bundle. The projection to $(T_2,T_3$ is a rational map undefined at the two singular points
$(0,0,1,0)$ and $(0,0,0,1)$. After we blow them up, we get the elliptic fibration. We have two obvious commuting involutions  $\sigma_1 = [t_0,t_1,-t_2,-t_3]$ and $\sigma_2 = [t_0,t_1,t_3,t_2]$. The locus of  fixed points of each of them is an elliptic curve with equation $T_2 = T_3$ or $T_2 = -T_3$. The group
$\la \sigma_1,\sigma_2\ra \cong 2^2$ permutes these two curves.

The groups of automorphisms of $S$ is easy to describe. It follows from Proposition \ref{autoex} that $G$ is a finite subgroup of $P\times K$, where $P$ is a subgroup of  $\PGL(2)$ leaving the zeros of $F_4$ invariant and $K$ is a subgroup of $\GL(2)$ leaving the zeroes of $t_2t_3$ invariant.  First we choose coordinates $T_0,T_1$ to write $F_4$ in the form $T_0^4+T_1^4+aT_0^2T_1^2, a^2 \ne 4$. It is always possible  if $F_4$ has 4 distinct roots (true in our case). Let $P$ be the subgroup of $\PGL(2)$ leaving the set of zeroes of $F_4$ invariant. It is one of the following groups $1, 2, 4, 2^2, D_8, A_4$. If $a^2 \ne 0, -12, 36$, then $P$ is  a subgroup of $2^2$. If $a^2 = 0, 36$, then $P$ is a subgroup of $D_8$. If $a^2 = -12$, then $P$ is a subgroup of $A_4$.

Suppose $a$ is not exceptional. Let $\bar{P}$ be the corresponding binary group. Then it leaves $F_4$ invariant, so $G$ is a subgroup of $K'\times P$, where $K'$ consists of matrices leaving $t_2t_3$ invariant. It is generated by diagonal matrices with determinant 1 and the transformations $[t_1,t_0]$.  Thus
$$G \subset D_{2n}\times 2^2,$$
where we use that a a finite subgroup of $K'$ is either cyclic or binary dihedral.

Suppose $a^2 = 0,36$ and $G$ contains an element of order 4. The form can be transformed to the form $T_0^4+T_1^4$. The value of the character at an element $\tau$ of order 4 is equal to $-1$.   We obtain
$$G \subset \half[Q\times  P] \cong (Q'\times P')^\bullet 2 \cong Q'{}_\bullet P$$
where $Q$ is a finite subgroup of $K$ and  the diagonal product is taken with respect to the subgroup $Q' = Q\cap K' $ of $Q$ and the subgroup $P\cap 2^2$ of index 2 of $P$. The group $Q'$ is cyclic or dihedral.

Suppose $a = -2\sqrt{3}i$ and $G$ contains an element $g$ of order 3 given by
$[\frac{1-i}{2}t_0+\frac{1-i}{2}t_1,-\frac{1+i}{2}t_0-\frac{1+i}{2}t_1]$. Its character is defined by $\chi(g) = \epsilon_3$. We obtain
$$G \subset \frac{1}{3}[Q\times P], $$
where the diagonal product is taken with respect to the subgroup $Q = \chi_2^{-1}(\mu_3)$ of $K$ and the subgroup $P' = P\cap 2^2$ of index 3 of $P$. Again $G$ is a subgroup of $D_{2n}\pt P$.
\end{example}

\subsection{Minimal conic bundles $G$-surfaces}
Now assume $(S,G)$ is a minimal $G$-surface such that $S$ admits a conic bundle map $\phi:S\to \bbP^1$.  As we had noticed before the number of singular fibres $k$ is greater or equal to $4$. Thus

\begin{equation}  K_S^2 = 8-k\le 4.
\end{equation}

Let $(S,G)$ be a  rational $G$-surface and
\begin{equation}\label{inj}
a:G \to \Or(\Pic(X)), \quad g\mapsto (g^*)^{-1}\end{equation}
be the natural representation of $G$ in the orthogonal group of the Picard group. We denote by $G_0$ the kernel of this representation. Since $k > 2$ and $G_0$ fixes  any component of a singular fibre, it acts identically on the base of the conic bundle fibration. Since $G_0$ fixes the divisor class of a section, and sections with negative self-intersection do not move in a linear system, we see that $G_0$ fixes pointwise any section with negative self-intersection. If we consider a section as a  point of degree 1 on the generic fibre, we see that $G_0$ must be is a cyclic group.

\begin{proposition}\label{first} Assume $G_0\ne \{1\}$.  Then  $S$ is an exceptional conic bundle.
\end{proposition}

\begin{proof} Let $g_0$ be  a non-trivial element from $G_0$. Let $E$ be a section with $E^2 = -n < 0$. Take an element $g\in G$ such that $E'= g(E)\ne E$.  Since  $g_0$ has two fixed points on each component of a singular fibre we obtain that $E$ and $E'$ do not intersect the same component. By Lemma \ref{ineqq}, we obtain that $k = 2n$. Now we blow down $n$ components in $n$ fibres  intersecting $E$ and $n$ components in the remaining $n$ fibres intersecting $E'$ to get a minimal ruled surface with two disjoint sections with self-intersection $0$. It must be isomorphic to $\bfF_0$.  So, we see that $S$ is an exceptional conic bundle (Construction 1)  with $n= g+1$.
\end{proof}

\smallskip

\emph{From now on in this section, we assume that $G_0= \{1\}$}.

\smallskip

Let $S_\eta$ be the general fibre of $\phi$. By Tsen's theorem it is
isomorphic to $\bbP_K^1$, where $K= \bbC(t)$ is the field of
rational functions of the base. Consider $S_\eta$ as a scheme over
$\bbC$. Then
$$\Aut_\bbC(S_\eta)  \cong \Aut_K(S_\eta):\PGL(2)\cong  \dJJ(2),$$
where $\dJJ(2)$ is a de Jonqui\`eres subgroup of $\Cr(2)$ and $\Aut_K(S_\eta)\cong \PGL(2,K)$.  A finite minimal group $G$ of automorphisms of a conic bundle is isomorphic to a subgroup of $\Aut_\bbC(S_\eta)$. Let $G_K = G\cap \Aut_K(S_\eta)$ and $G_B \cong G/G_K$ be the image of $G$ in $\PGL(2)$. We have an extension of groups
\begin{equation}\label{exten}
1\to G_K \to G \to G_B\to 1
\end{equation}

Let $\calR$ be the subgroup of $\Pic(S)$ spanned by the divisor classes of $R_i-R_i', i = 1,\ldots,k.$ It is obviously $G$-invariant and  $\calR_\bbQ$ is equal to the orthogonal complement of $\Pic(S)_\bbQ^G$ in $\Pic(S)_\bbQ$. The  orthogonal group of the quadratic lattice $\calR$  is isomorphic  to the wreath product $2\wr S_k$. The normal subgroup $2^k$ consists of transformations which switch some of the  $R_i$'s  with $R_i'$. A subgroup isomorphic to $S_k$  permutes  the classes  $R_i-R_i'$.

 \begin{lemma}\label{fixlocus} Let $G$ be a  minimal group of automorphisms of $S$. There exists an element $g\in G_K$ of  order 2 which  switches the components of some singular  fibre.  \end{lemma}

\begin{proof} Since $G$ is minimal, the $G$-orbit of any $R_i$ cannot consist of disjoint components of fibres. Thus it contains a pair $R_j, R_j' $ and hence there exists an element $g\in G$ such that $g(R_j) = R_j'$. If $g$ is of odd order $2k+1$, then $g^{2k}$ and $g^{2k+1}$ fix $R_j$, hence $g$ fixes $R_j$. This contradiction shows that $g$ is of even order $2m$. Replacing $g$ by an odd power, we may assume that $g$ is of order $m = 2^a$.

Assume $a = 1$. Obviously the singular point $p = R_j\cap R_j'$ of the fibre belongs to the fixed locus  $S^g$ of $g$. Suppose $p$ is an isolated fixed point. Then we can choose local coordinates at $p$  such that  $g$ acts   by $(z_1,z_2)\mapsto (-z_1,-z_2)$, and hence acts identically on the tangent directions. So it cannot switch the components. Thus  $S^g$ contains a curve  not contained in fibres which passes through $p$. This implies that $g\in G_K$.

Suppose $a > 1$. Replacing $g$ by $g'= g^{m/2}$ we get an automorphism of order 2 which fixes the point $x_j$ and the components $R_j,R_j'$. Suppose $S^{g'}$ contains one of the components, say $R_j$. Take a general point $y\in R_j$. We have
$g'(g(y)) = g(g'(y)) = g(y).$
This shows that $g'$ fixes $R_j'$ pointwise. Since $S^{g'}$ is smooth, this is impossible. Thus  $g'$ has 3 fixed points $y,y',p$ on $F_j$, two on each component. Suppose $y$ is an isolated fixed point lying on $R_j$. Let $\pi:S\to S'$ be the blowing down of $R_j$. The element $g'$ descends to an automorphism of order 2 of $S'$ which has an isolated fixed point at $q= \pi(R_j)$.
Then it acts identically on the tangent directions at $q$, hence on $R_j$. This contradiction shows that $S^{g'}$ contains a curve intersecting $F_j$ at $y$ or at $p$, and hence $g'\in G_K$.  If $g'$ does not switch components of any fibre, then it acts identically on $\Pic(S)$. By our assumption it is impossible. Thus $g'$ has to switch components of some fibre. 
\end{proof}

The restriction of the homomorphism $G\to \Or(\calR) \cong 2^k:S_k$ to $G_K$ defines a surjective homomorphism
$$\rho: G_K\to 2^s, \ s\le k.$$
 An element from $\Ker(\rho)$ acts identically on $\calR$ and hence on $\Pic(S)$. By Lemma \ref{fixlocus}, $G_K$ is not trivial and $s > 0$.
A finite subgroup of $\PGL(2,K)$ does not admit a surjective homomorphism to $2^s$ for $s > 2$. Thus $s = 1$ or $ 2$.

\medskip\noindent

Case 1: $s = 1$.

Let $G_K = \la h\ra$.  The element $h$   fixes two  points on each
nonsingular fibre. The closure of these points is a one-dimensional
component $C$ of $S^h$. It is a smooth bisection of the fibration. If $C$ is reducible, then it consists of two disjoint sections. Each fibre meets one of the components at a nonsingular point. This implies that $h$ cannot switch components of any fibre and hence acts identically on $\Pic(S)$. The latter case has been excluded. So, $C$ is a smooth irreducible curve and the projection $S\to \bbP^1$ defines a double cover $f:C\to \bbP^1$. Let $x$ be a ramification point of the double cover. Since no smooth $h$-invariant curve can be tangent to $C$ (its image under the projection to the quotient $S/\la h\ra$ locally splits in the cover at the intersection point  with the branch divisor), the point $x$ must be the  singular point of the fibre passing through $x$. The components of this  fibre must be switched by $h$ because otherwise we would have three invariant tangent directions at $x$. Conversely, if the components of a fibre are switched by $h$, the singular point of the fibre must be a ramification point of $f:C\to \bbP^1$. Thus we obtain that the number of ramification points is equal to the number $m$ of fibres whose components are switched by $h$, and by Hurwitz's formula, the genus of $C$ is equal to $(m-2)/2$.

\medskip\noindent
Case 2: $s = 2$.

Let $g_1,g_2$ be two elements from $G_K$ which are mapped to
generators of the image of $G_K$ in $2^k$. Let $C_1$ and $C_2$ be
one-dimensional components of the sets $S^{g_1}$ and $S^{g_2}$. As
in the previous case we show that $C_1$ and $C_2$ are smooth
irreducible curves of genera $g(C_1)$ and $g(C_2)$. Let $\Sigma_1$
and $\Sigma_2$ be the sets of branch points of the corresponding double
covers.  For any point $x\in \Sigma_1\cap\Sigma_2$ the
transformation $g_3 = g_1g_2$ fixes the components of the fibre
$F_x$. For any point $x\in \Sigma_1\setminus \Sigma_2$, $g_3$
switches the components of $F_x$. Let $C_3$ be the one-dimensional
component of $S^{g_3}$ and $\Sigma_3$ be the set of branch points of
$C_3$.  We see that $\Sigma_i = \Sigma_j+\Sigma_k$ for distinct
$i,j,k$, where $ \Sigma_j+\Sigma_k=
(\Sigma_j\cup\Sigma_k)\setminus(\Sigma_j\cap\Sigma_k)$. This implies
that there exist three  binary forms
$p_1(t_0,t_1),p_2(t_0,t_1),p_3(t_0,t_1)$, no two of which have a
common root, such that  $\Sigma_1 = V(p_2p_3), \Sigma_2 = V(p_1p_3),
\Sigma_3 = V(p_1p_2)$. Setting $d_i = \deg p_i$, we get
$$2g(C_i)+2  = d_j+d_k.$$

Let us summarize what we have learnt \footnote{The statement of this theorem in the published version of the paper contains  a wrong assertion that the number $m$ is equal to the number of singular fibres. The mistake was in using a wrong argument in the proof of Lemma 5.6. We thank J. Blanc for pointing out this mistake in his paper \cite{Bl}.}.

\begin{theorem}\label{itog} Let $G$ be a minimal finite group of automorphisms of a conic bundle $\phi:S\to \bbP^1$ with a set $\Sigma$ of singular fibres.  Assume $G_0=\{1\}$. Then $k = \#\Sigma > 2$ and one of the following cases occurs.
\begin{itemize}
\item[(1)] $G=2P$, where the central involution $h$ fixes pointwise an irreducible smooth bisection $C$ of $\pi$ and switches the components in $m$ fibres over the branch points of the $g_2^1$ on $C$ defined by the projection $\pi$. The curve $C$
is a curve of genus $g = (m-2)/2$.  The group $P$ is isomorphic to a group of automorphisms of $C$ modulo the involution defined by the $g_2^1$.
\item[(2)] $G \cong 2^2\pt P$, each nontrivial  element $g_i$ of the subgroup $2^2$ fixes pointwise an irreducible smooth bisection $C_i$. The set $\Sigma$ contains a subset $\Sigma'$ which is partitioned in 3 subsets $\Sigma_1,\Sigma_2,\Sigma_3$ such that the projection $\phi:C_i\to \bbP^1$ ramifies over  $\Sigma_j + \Sigma_k , i\ne j\ne k$.
The group $P$  is  a  subgroup of $\Aut(\bbP^1)$ leaving the set $\Sigma$ and its  partition into 3 subsets $\Sigma_i$ invariant.
\end{itemize}
\end{theorem}

It follows from Lemma \ref{coh} that in Case 1  the non-split extension is isomorphic to a binary polyhedral group, unless $G = O$ or $D_{2n}$, where $n$ is even.

\begin{remark} It follows from the previous description that any abelian group $G$ of automorphisms of a conic bundle must be a subgroup of some extension $Q\pt P$, where $Q$ is a dihedral, binary dihedral or cyclic group, and $P$ is a polyhedral group. This implies that $G$ is either a cyclic group, or a group $2\times m$, or $2^2\times m$,  or $2^4$. All these groups occur (see  Example \ref{blanc} and  \cite{Blanc}).
\end{remark}

\subsection{Automorphisms of hyperelliptic curves}\label{5.5}
We consider a  curve of genus  g equipped with a linear series $g_2^1$ as a curve $C$ of degree $2g+2$ in $\bbP(1,1,g+1)$ given by an equation
$$T_2^2+F_{2g+2}(T_0,T_1) = 0.$$
An automorphism $\sigma$ of $C$  is defined by a transformation
$$(t_0,t_1,t_2)\mapsto (at_1+bt_0,ct_1+dt_0,\alpha t_2),$$
where $\left(\begin{smallmatrix}a&b\\c&d\end{smallmatrix}\right)\in \SL(2)$ and
$F(aT_0+bT_1,cT_0+dT_1) = \alpha^2F(T_0,T_1)$. So to find the group of automorphisms of $C$ we need to know  relative invariants  $\Phi(T_0,T_1)$ for finite subgroups $\bar{P}$ of $\SL(2,\bbC)$ (see \cite{Springer}). The set of relative invariants is a finitely generated $\bbC$-algebra. Its generators are called Gr\"undformen.
We will list the Gr\"undformen (see \cite{Springer}). We  will  use them  later for the description of automorphism groups of Del Pezzo surfaces of degree 1.

\medskip\noindent
\begin{itemize}
\item \emph{$\bar{P}$ is a cyclic group of  order $n$}.
\end{itemize}
A generator is given by the matrix
$$g = \begin{pmatrix}\epsilon_n&0\\
0&\epsilon_n^{-1}\end{pmatrix}.$$
The Gr\"undformen are $t_0$ and $t_1$ with characters determined by
$$\chi_1(g) = \epsilon_n, \quad \chi_2(g) = \epsilon_n^{-1}.$$

\medskip\noindent
\begin{itemize}
\item \emph{$\bar{P} \cong Q_{4n}$ is a binary dihedral group of order  $4n$}.
\end{itemize}

Its generators are given by the matrices
$$g_1 = \begin{pmatrix}\epsilon_{2n}&0\\
0&\epsilon_{2n}^{-1}\end{pmatrix}, \quad g_2 = \begin{pmatrix}0&i\\
i&0\end{pmatrix}.$$
 The Gr\"undformen are
\begin{equation}\label{cyclic}
\Phi_1 = T_0^n+T_1^n, \quad \Phi_2 = T_0^n-T_1^n, \quad \Phi_3 = T_0T_1.
\end{equation}
The generators $g_1$ and $g_2$ act on the Gr\"undformen with characters
$$\chi_1(g_1) = \chi_2(g_1) =  -1, \chi_1(g_2) = \chi_2(g_2) = i^{n},$$
$$\chi_3(g_1) = 1, \chi_3(g_2) = -1.$$

\medskip\noindent
\begin{itemize}
\item \emph{$\bar{P}$ is a binary tetrahedral group of order  $24$}.
\end{itemize}

Its generators are given by the matrices
$$g_1 = \begin{pmatrix}\epsilon_{4}&0\\
0&\epsilon_{4}^{-1}\end{pmatrix}, \quad g_2 = \begin{pmatrix}0&i\\
i&0\end{pmatrix}, \quad g_3 = \frac{1}{i-1}\begin{pmatrix}i&i\\
1&-1\end{pmatrix}.$$
The Gr\"undformen are
$$\Phi_1 = T_0T_1(T_0^4-T_1^4),\  \Phi_2, \ \Phi_3  = T_0^4\pm 2\sqrt{-3}T_0^2T_1^2+T_1^4.$$
The generators $g_1,g_2, g_3$ act on the Gr\"undformen with characters
$$\chi_1(g_1) = \chi_1(g_2) =  \chi_1(g_3) = 1,$$
$$\chi_{2}(g_1) = \chi_2(g_2) = 1,\  \chi_2(g_3) = \epsilon_3,$$
$$\chi_{3}(g_1) = \chi_3(g_2) = 1, \ \chi_3(g_3) = \epsilon_3^2.$$

\medskip\noindent
\begin{itemize}
\item \emph{ $\bar{P}$ is a binary octahedral group of order 48}.
\end{itemize}

Its generators are
$$g_1 = \begin{pmatrix}\epsilon_8&0\\
0&\epsilon_8^{-1}\end{pmatrix}, \quad g_2 =  \begin{pmatrix}0&i\\
i&0\end{pmatrix}, \quad g_3 =  \frac{1}{i-1}\begin{pmatrix}i&i\\
1&-1\end{pmatrix}.$$
The  Gr\"undformen are
$$\Phi_1 = T_0T_1(T_0^4-T_1^4),\ \Phi_2 = T_0^8+14T_0^4T_1^4+T_1^8,\ \Phi_3 = (T_0^4+T_1^4)((T_0^4+T_1^4)^2-36T_0^4T_1^4).$$
The generators $g_1,g_2, g_3$ act on the Gr\"undformen with characters
$$\chi_1(g_1) = -1, \chi_1(g_2) =  \chi_1(g_3) = 1,$$
$$\chi_{2}(g_1) = \chi_2(g_2) = \chi_{2}(g_3) = 1,$$
$$\chi_3(g_1) = -1, \quad \chi_3(g_2) = \chi(g_3) = 1.$$

\medskip\noindent
\begin{itemize}
\item \emph{$\bar{P}$ is a binary icosahedral group of order 120}.
\end{itemize}

Its generators are
$$g_1 = \begin{pmatrix}\epsilon_{10}&0\\
0&\epsilon_{10}^{-1}\end{pmatrix}, \quad g_2 = \begin{pmatrix}0&i\\
i&0\end{pmatrix}, \quad g_3 = \frac{1}{\sqrt{5}}\begin{pmatrix}\epsilon_5-\epsilon_5^4&\epsilon_5^2-\epsilon_5^3\\
\epsilon_5^2-\epsilon_5^3&-\epsilon_5+\epsilon_5^4\end{pmatrix}.$$
The Gr\"undformen are
$$\Phi_1 = T_0^{30}+T_1^{30}+522(T_0^{25}T_1^5-T_0^{5}T_1^{25})-10005(T_0^{20}T_1^{10}+T_0^{10}T_1^{20}),$$
$$\Phi_2 = -(T_0^{20}+T_1^{20})+228(T_0^{15}T_1^5-T_0^{5}T_1^{15})-494T_0^{10}T_1^{10},$$
$$\Phi_3 = T_0T_1(T_0^{10}+11T_0^5T_1^5-T_1^{10}).$$
Since $P/(\pm 1) \cong A_5$ is a simple group and all Gr\"undformen are of even degree, we easily see that the characters are trivial.

\subsection{Commuting de Jonqui\`eres involutions}
Recall that a de Jonqui\`eres involution $IH_{g+2}$ is regularized
by an automorphism of the surface $S$ which is obtained from
$\bfF_1$ by blowing up $2g+2$ points. Their images on $\bbP^2$ are
the $2g+2$ fixed points of the  involution of  $H_{g+2}$. Let $\pi:S\to X = S/IH_{g+2}$. Since the
fixed locus of the involution is a smooth hyperelliptic curve of
genus $g$, the quotient surface $X$ is a nonsingular surface. Since
the components of singular fibres of the conic bundle on $S$ are
switched by $IH_{g+2}$, their images on $X$ are isomorphic to
$\bbP^1$. Thus $X$ is a minimal ruled surface $\bfF_e$. What is $e$?

Let $\bar{C} = \pi(C)$ and $\bar{E} = \pi(E)$, where $E$ is the exceptional section on $S$. The curve $\bar{E}$ is a section on $X$ whose preimage in the cover splits. It is either tangent to   $\bar{C}$ at any of its of $g$ intersection points (since $IH_{g+2}(E)\cdot E = g$) or is  disjoint from $\bar{C}$ if $g = 0$.  Let $s$ be the divisor class of a section on $\bfF_e$ with self-intersection $-e$ and $f$ be  the class  of a fibre. It is easy to see that
$$\bar{C} = (g+1+e)f+2s, \quad \bar{E} = \frac{g+e-1}{2}f+s.$$
Let $\bar{R}$ be  a section with the divisor class $s$. Suppose $\bar{R} =
\bar{E}$, then $\bar{R}\cdot \bar{C} = g+1-e = 2g$ implies $g =
1-e$, so $(g,e) = (1,0)$ or $(0,1)$. In the first case, we  get an
elliptic curve on $\bfF_0$ with divisor class $2f+2s$ and $S$ is
non-exceptional conic bundle with $k=4$. In the second case $S$ is
the conic bundle (nonminimal) with $k=2$.

Assume that $(g,e) \ne (1,0)$. Let $R = \pi^{-1}(\bar{R})$ be the
preimage of $\bar{R}$. We have $R^2 = -2e$. If it splits into two
sections $R_1+R_2$, then $R_1\cdot R_2 =  \bar{C}\cdot \bar{R} =
g+1-e$, hence $-2e = 2(g+1-e)+2R_1^2$ gives $R_1^2 = -g-1$. Applying
Lemma \eqref{ineqq}, we get $R_1\cdot R_2 = g+1-e = g-1+(2g+2-a)/2 =
-a/2$, where $a\ge 0$. This gives $e=g+1$, but intersecting
$\bar{E}$ with $\bar{R}$ we get $e\le g-1$. This contradiction shows
that $\bar{R}$ does not split, and hence $R$ is an irreducible
bisection of the conic bundle with $R^2 = -2e$. We have $R\cdot E=
(g-e-1)/2, R\cdot R_i = R\cdot R_i' = 1$, where $R_i+R_i'$ are
reducible fibres of the conic fibration.

This shows that the image of  $R$ in the plane is a hyperelliptic  curve $H_{g'+2}'$ of degree $d = (g-e+3)/2 $ and genus $g' = d-2 = (g-e-1)/2$ with the  point
 $q$ of multiplicity $g'$.  It also passes through the points $p_1,\ldots,p_{2g+2}$. Its Weierstrass points $p_1',\ldots,p_{2g'+2}'$ lie on $H_{g+2}$.  Here we use the notation from  \ref{dejon}. Also the curve $H_{g'+2}'$ is  invariant with respect to the de Jonqui\`eres involution.

 Write the  equation of $H_{g'+2}'$ in the form
\begin{equation}\label{hyp'}
A_{g'}(T_0,T_1)T_2^2+2A_{g'+1}(T_0,T_1)T_2+A_{g'+2}(T_0,T_1) =
0.\end{equation} It follows from the geometric definition of the de
Jonqui\`eres involution that we have the following relation between
the equations of $H'_{g'+2}$ and $H_{g+2}$ (cf. \cite{Coble}, p.126)
\begin{equation}\label{rel1}
F_gA_{g'+2}-2F_{g+1}A_{g'+1}+F_{g+2}A_{g'} = 0.\end{equation}
Consider this as a system of linear equations with coefficients of $A_{g'+2},A_{g'+1},A_{g'}$ considered as the unknowns. The number of the unknowns is equal to $(3g-3e+9)/2$. The number of the equations is $(3g-e+5)/2$. So, for a general $H_{g+2}$ we can solve these equations only if $g= 2k+1, e = 0, d = k+2$ or $g = 2k, e = 1, d = k+1$.
Moreover, in the first case we get a pencil of curves $R$ satisfying these properties, and in the second case we have a unique such curve (as expected). Also the first case  covers our exceptional case $(g,e) = (1,0)$.

For example, if we take $g = 2$ we obtain that the six Weierstrass points $p_1,\ldots,p_6$ of $H_{g+2}$ must be on a conic. Or, if $g = 3$, the eight Weierstrass points together with the point $q$ must be the base points of a pencil of cubics. All these properties are of course not expected for a general set of 6 or 8 points in the plane.

To sum up, we have proved the following.

\begin{theorem}\label{itog2} Let $H_{g+2}$ be a  hyperelliptic curve of degree $g+2$  and genus $g$ defining a de Jonqui\`eres involution $IH_{g+2}$. View this involution as an automorphism $\tau$ of order 2 of the surface $S$ obtained by blowing up the singular point $q$ of $H_{g+2}$ and its $2g+2$ Weierstrass points $p_1,\ldots,p_{2g+2}$. Then
\begin{itemize}
\item[(i)]  the quotient surface $X = S/(\tau)$ is isomorphic to $\bfF_e$ and the ramification curve is $C = S^\tau$;
\item[(ii)] if $H_{g+2}$ is a general hyperelliptic curve then  $e=0$  if $g$ is odd and $e=1$ if $g$ is even;
\item [(iii)] the branch curve $\bar{C}$ of the double cover $S\to \bfF_e$ is a curve from the divisor class $(g+1+e)f+2s$;
\item [(iv)] there exists a section from the divisor class $\frac{g+e-1}{2}f+s$ which is  tangent to $\bar{C}$ at each $g$ intersection points unless $g= 0, e = 1$ in which case it is disjoint from $\bar{C}$;
\item [(v)] the reducible fibres of the conic bundle on $S$ are the preimages of the $2g+2$ fibres from the pencil $|f|$ which are tangent to $\bar{C}$;
\item [(vi)] the preimage of a section from the divisor class $s$ either splits if $(g,e) = (1,0)$ or a curve of genus $g = 0$, or a hyperelliptic curve $C'$ of genus $g' = (g-e-1)/2\ge 1$ which is invariant with respect to $\tau$. It intersects the hyperelliptic curve $C$ at its $2g'+2$ Weierstrass points.
\item [(vii)] the curve $C'$ is uniquely defined if $e > 0$ and varies in a pencil if $e = 0$.
\end{itemize}
\end{theorem}

Let $IH_{g'+2}'$ be the de Jonqui\`eres involution defined by the curve $H_{g'+2}'$ from equation \eqref{hyp'}. Then it can be given in affine  coordinates by formulas \eqref{dj2}, where $F_i$ is replaced with $A_i$. Thus we have two involutions defined by the formulas
\begin{eqnarray}
IH_{g+2}:\
(x',y') &= &\bigl(x, \frac{-yP_{g+1}(x)-P_{g+2}(x)}{P_g(x)y+P_{g+1}(x)}\bigr),\\ \notag
IH_{g'+2}':(x',y') &= &\bigl(x, \frac{-yQ_{g'+1}(x)-Q_{g'+2}(x)}{Q_g(x)y+Q_{g'+1}(x)}\bigr),\\ \notag
\end{eqnarray}
where $P_i$ are the dehomogenizations of the $F_i$'s and $Q_i$ are the dehomogenizations of the $A_i$'s. Composing them in both ways we see that  the relation \eqref{rel1} is satisfied if and only if  the two involutions commute. Thus a de Jonqui\`eres involution can be always included in a group of de Jonqui\`eres transformations isomorphic to $2^2$. In fact, for a general $IH_{g+2}$ there exists a unique such group if $g$ is even and there is a $\infty^1$ such groups when $g$ is odd. It is easy to check   that the involution $IH_{g+2}\circ  H_{g'+2}'$ is the de Jonqui\`eres involution defined by the third hyperelliptic curve with equation
\begin{equation}\label{hyp2}
\det\begin{pmatrix}F_g&F_{g+1}&F_{g+2}\\
A_{g'}&A_{g'+1}&A_{g''+2}\\
1&-T_2&T_2^2\end{pmatrix}
= B_{g''}T_2^2+2B_{g''-1}T_2+B_{g''+2} = 0,\end{equation}
(cf. \cite{Coble}, p.126).

If we blow up the Weierstrass point of the curve $C'$ (the proper transform of $H_{g'+2}'$ in $S$), then we get a conic bundle surface $S'$ from case (2) of Theorem \ref{itog}.

\subsection{A question on extensions} It  still remains to decide  which extensions
\begin{equation}\label{extension}
1\to G_K\to G\to G_B\to 1
\end{equation} describe minimal groups of automorphisms of conic bundles. We do not have the full answer and only make a few remarks and examples. Lemma \ref{coh} helps to decide on splitting in the case when $G_K$ is abelian and central.

\begin{example}\label{extension1} Consider a de Jonqui\`eres transformation
$$\djj_P:(x,y) \mapsto (x, P(x)/y),$$
where $P(T_1/T_0) = T_0^{-2g}F_{2g+2}(T_0,T_1)$ is a dehomogenization of a homogeneous
polynomial $F_{2g+2}(T_0,T_1)$ of degree $2g+2$ defining a
hyperelliptic curve of genus $g$. Choose $F_{2g+2}$ to be a relative
invariant of a binary polyhedral group $\bar{P}$ with
character $\chi:\bar{P}\to \bbC^*$. We assume that $\chi = \alpha^2$
for some character $\alpha:\bar{P}\to \bbC^*$. For any $g =
\left(\begin{smallmatrix}a&b\\c&d\end{smallmatrix}\right)\in
\bar{P}$ define the transformation
$$g:(x,y)\mapsto (\frac{ax+b}{cx+d},\alpha(g)(cx+d)^{-g-1}y).$$
We have
$$P(\frac{ax+b}{cx+d}) = \alpha^2(g)(cx+d)^{-2g-2}P(x).$$
It is immediate to check that $g$ and $\djj_p$ commute. The matrix $ -I_2$ defines the transformation
$g_0:(x,y)\mapsto (x,\alpha(-I_2)(-1)^{g+1}y)$. So, if
$$\alpha(-I_2) = (-1)^{g+1},$$
the action of $\bar{P}$ factors through $P$ and together with
$\djj_P$ generate the group $2\times P$. On the other hand, if
$\alpha(-I_2) =  (-1)^{g}$, we get the group $G = 2\times
\bar{P}$. In this case the group $G$ is regularized on an
exceptional conic bundle with $G_0 \cong 2$. The generator
corresponds to the transformation $g_0$.
\end{example}

Our first general observation is that the extension $G= 2P$ always splits  if $g$ is even, and of course, if $P$ is a cyclic group of odd order.  In fact, suppose $G$ does not split. We can always find an element $g\in G$ which is mapped to an element $\bar{g}$ in $P$ of even order $2d$ such that $g^{2d} = g_0\in G_K$. Now $g_1 = g^d$ defines an automorphism of order 2 of the hyperelliptic curve $C = S^{g_0}$ with  fixed points lying over two fixed points of $\bar{g}$ in $\bbP^1$. None of these points belong to $\Sigma$, since otherwise $g_0$, being a square of $g_1$, cannot switch the components of the corresponding fibre. Since $g_1$ has two fixed points on the invariant fibre and both of them must lie on $C$, we see that  $g_1$ has 4 fixed points. However this contradicts the Hurwitz formula.

\smallskip
Recall that a double cover $f:X\to Y$ of nonsingular varieties with branch divisor $W\subset Y$ is given by an invertible sheaf $\calL$ together with a section $s_W\in \Gamma(Y,\calL^2)$ with zero divisor  $W$. Suppose a group $G$ acts on $Y$ leaving invariant $W$. A \emph{lift} of $G$ is  a group
$\widetilde{G}$ of automorphisms of $X$  such that it commutes with the covering involution $\tau$ of $X$ and the corresponding homomorphism $\widetilde{G}\to \Aut(Y)$ is an isomorphism onto the group $G$.

The following lemma is well-known and is left to the reader.

\begin{lemma} A subgroup $G\subset \Aut(Y)$ admits a lift if and only if $\calL$ admits a $G$-linearization and in the corresponding representation of $G$ in $\Gamma(Y,\calL^2)$ the section $s_W$ is  $G$-invariant.
\end{lemma}

\begin{example}\label{blanc} Let $p_{i}(t_0,t_1), i =0,1,2,$ be binary forms of degree $d$. Consider a curve $C$ in $\bfF_0 \cong \bbP^1\times \bbP^1$ given by an equation
$$F= p_{0}(t_0,t_1)x_0^2+2p_{1}(t_0,t_1)x_0x_1+p_2(t_0,t_1)x_1^2 = 0.
$$
Assume that the binary form $D = p_1^2-p_0p_2$  does not have multiple roots. Then $C$ is a nonsingular hyperelliptic curve of genus $d-1$. Suppose $d = 2a$ is even so that the genus of the curve is odd. Let $P$ be a polyhedral group not isomorphic to a cyclic group of  odd order. Let $V = \Gamma(\bbP^1,\calO_{\bbP^1})$ and $\rho:P\to \GL(S^{2a}V\otimes S^2V)$ be its natural representation, the tensor product of the two natural representations of $P$ in the space of binary forms of even degree. Suppose that  $F\in S^{2a}V\otimes S^2V$ is an invariant. Consider the double cover $S\to \bfF_0$ defined by the section $F$ and the invertible sheaf $\calL = \calO_{\bfF_0}(a,1)$. Now assume additionally that  $P$ does not have a linear representation in $S^aV\otimes V$ whose tensor square is equal to $\rho$. Thus $\calL$ does not admit a $P$-linearization and  we cannot lift $P$ to a group of automorphisms of the double cover. However, the  binary polyhedral group $\bar{P}$ lifts.
Its central involution acts identically on $\bfF_0$, hence lifts to the covering involution of $S$. It follows from the discussion in the previous subsection that  $S$ is a non-exceptional conic bundle, and  the group $\bar{P}$ is a minimal group of automorphisms of $S$ with $G_K \cong 2$ and $G_B \cong P$.

Here is a concrete example. Take
$$p_0 = t_0t_1(t_0^2+t_1^2), \quad p_1 =t_0^4+t_1^4, \quad p_2=t_0t_1(t_0^2-t_1^2).$$
Let $\bar{P}\subset \SL(2)$ be a cyclic group of order 4  that acts  on the variables $t_0,t_1$  via the transformation $[it_0,-it_1]$ and on the variables $x_0,x_1$ via the transformation $[ix_0,-ix_1]$. Then $\bar{P}$ acts on $S^2V\otimes V$ via $[-1,1,-1]\otimes [i,-i]$. The matrix $-I_2$ acts as $1\otimes -1$ and hence $P= \bar{P}/(\pm I_2)$ does not act on $S^2V\otimes V$. This realizes the cyclic group $C_4$ as a minimal group of automorphisms of a conic bundle with $k = 2g+2 = 8$.
\end{example}

The previous example shows that for any $g\equiv 1 \mod 4$ one can realize a binary polyhedral group $\bar{P}= 2.P$ as a minimal group of automorphisms of a conic bundle with $2g+2$ singular fibres. We do not know whether the same is true for $g\equiv 3 \mod 4.$

\begin{example}\label{example} Let $p_i(t_0,t_1), i = 1,2,3, $ be three binary forms of even degree $d$ with no multiple roots. Assume no two have common zeroes. Consider a surface $S$ in $\bbP^1\times \bbP^2$ given by a bihomogeneous form of degree $(d,2)$
\begin{equation}
p_1(t_0,t_1)z_0^2+p_2(t_0,t_1)z_1^2+p_3(t_0,t_1)z_2^2 = 0,
\end{equation}
The surface is nonsingular. The projection to $\bbP^1$ defines a conic bundle structure on $S$ with
singular fibres over the zeroes of the polynomials $p_i$. The curves $C_i$ equal to the preimages of the lines $z_i = 0$ under the second projection are hyperelliptic curves of genus $g = d-1$.  The automorphisms $\sigma_1,\sigma_2$ defined by the negation of one of the first two coordinates $z_0,z_1,z_2$ form a subgroup of $\Aut(S)$ isomorphic to $2^2$. Let $P$ be a finite subgroup of $\SL(2,\bbC)$  and $g\mapsto g^*$ be its natural action on the space of binary forms. Assume that
$p_1,p_2,p_3$ are relative invariants of $P$ with characters  $\chi_1,\chi_2,\chi_3$ such that we can write them  in the form $\eta_i^2$ for some characters $\eta_1,\eta_2,\eta_3$ of $P$. Then $P$ acts on $S$ by the formula
$$g((t_0,t_1),(z_0,z_1,z_2)) = ((g^*(t_0),g^*(t_1)),(\eta_1(g)^{-1}z_0,\eta_2(g)^{-1}z_1,\eta_3(g)^{-1}z_2)).$$

For example, let $P = \la g \ra$ be a cyclic group of order 4. We take $p_1 = t_0^2+t_1^2, p_2 = t_0^2-t_1^2, p_3 = t_0t_1$. It acts on $S$ by the formula
$$g:((t_0,t_1),(z_0,z_1,z_2)) \mapsto ((it_1,it_0),(iz_0,z_1,iz_2)).$$
Thus $g^2$ acts identically on $t_0,t_1, z_1$ and multiplies $z_0,z_2$ by $-1$. We see that  $G_K = \la g^2\ra$ and  the  extension $1\to G_K\to G\to G_B\to 1$ does not split. If we add to the group the transformation $(t_0,t_1, z_0,z_1,z_2)\mapsto (t_0,t_1, z_0,-z_1,z_2)$ we get a non-split extension
$2^{2+1}$.

On the other hand, let  us replace $p_2$ with $t_0^2+t_1^2+t_0t_1$. Define $g_1$ as acting only on $t_0,t_1$ by $[it_1,it_0]$, $g_2$ acts only on $z_0$ by $z_0\mapsto -z_0$ and $g_3$ acts only on $z_1$ by $z_1\mapsto -z_1$.  We  get the groups $\la g_1,g_2\ra = 2^2$ and $\la g_1,g_2,g_3\ra = 2^3$.

In another example we take $P$ to be the dihedral group $D_8$. We take $p_1 = t_0^2+t_1^2, p_2 = t_0^2-t_1^2, p_3 = t_0t_1$. It acts on $S$ by the formula
$$\left(\begin{smallmatrix}i&0\\
0&-i\end{smallmatrix}\right):((t_0,t_1),(z_0,z_1,z_2)) \mapsto ((it_0,-it_1),(iz_0,iz_1,z_2)),$$
$$\left(\begin{smallmatrix}0&i\\
i&0\end{smallmatrix}\right):((t_0,t_1),(z_0,z_1,z_2)) \mapsto ((it_1,it_0),(z_0,iz_1,z_2)),$$
The scalar matrix $c = -I_2$ belongs to $G_K\cong 2^2$ and the quotient $P/(c)\cong 2^2$ acts faithfully  on the base. This gives a non-split extension $2^{2+2}$.

Finally, let us take
$$p_1 = t_0^4+t_1^4,\  p_2 = t_0^4+t_1^4+t_0^2t_1^2, \ p_3 = t_0^4+t_1^4-t_0^2t_1^2.$$
These are invariants  for the group $D_4$ acting via $g_1:(t_0,t_1)\mapsto (t_0,-t_1), \ g_2:(t_0,t_1)\mapsto (t_1,t_0)$. Together with transformations $\sigma_1,\sigma_2$ this generates the group $2^4$
(see another realization of this group in \cite{Beauville}).
\end{example}

\section{Automorphisms of Del Pezzo surfaces}
\subsection{The Weyl group} Let $S$ be a Del Pezzo surface of degree $d$ not isomorphic to $\bbP^2$ or $\bfF_0$. It is  isomorphic to the blowup of $N = 9-d \le 8$ points in $\bbP^2$ satisfying the conditions of generality from section \ref{minimal}. The blowup of one or 2 points is obviously nonminimal (since the exceptional curve in the first case and the proper transform of the line through the two points is $G$-invariant). So we may assume that $S$ is a Del Pezzo surface of degree $d\le 6$.

Let $\pi:S\to \bbP^2$ be the blowing-up map. Consider the
factorization \eqref{decom} of $\pi$ into a composition of blowups
of $N = 9-d$ points. Because of the generality condition, we may assume
that none of the points $p_1,\ldots, p_N$ is infinitely near, or,
equivalently, all exceptional curves  $\calE_i$ are
irreducible curves. We identify them with curves $E_i =
\pi^{-1}(p_i)$.  The divisor classes $e_0 = [\pi^*(\text{line}], e_i
= [E_i], i = 1,\ldots, N,$ form a basis of $\Pic(S)$. It is called a
\emph{geometric basis}.

Let
$$\alpha_1 = e_0-e_1-e_2-e_3, \alpha_2 = e_1-e_2,\ldots,\alpha_N = e_{N-1}-e_N.$$
For any $i = 1,\ldots,N$ define a \emph{reflection} isometry $s_i$ of the abelian group $\Pic(S)$
$$s_i:x\mapsto  x+(x\cdot \alpha_i)\alpha_i.$$
Obviously, $s_i^2 = 1$ and $s_i$ acts identically on the orthogonal complement of $\alpha_i$. Let $W_S$ be the group of  automorphisms of $\Pic(S)$ generated by
the transformations $s_1,\ldots,s_N$. It is called the \emph{Weyl
group of $S$}. Using the basis $(e_0,\ldots,e_N)$ we identify $W_S$
with a group of isometries of the odd unimodular quadratic form
$q:\bbZ^{N+1}\to \bbZ$ of signature $(1,N)$ defined by
$$q_N(m_0,\ldots,m_N) = m_0^2-m_1^2-\cdots-m_N^2.$$
 Since
$K_S = -3e_0+e_1+\cdots+e_N$ is orthogonal to all $\alpha_i$'s, the
image of $W_S$ in $\Or(q_N)$ fixes the vector $k_N =
(-3,1,\ldots,1)$. The subgroup of $\Or(q_N)$ fixing $k_N$ is denoted
by $W_N$ and is called the \emph{Weyl group of type $E_N$}. The orthogonal complement
$\calR_N$ of $k_N$ equipped with the restricted inner-product, is called the \emph{root lattice} of $W_N$.

We denote by $\calR_S$ the sublattice of $\Pic(S)$ equal to the
orthogonal complement of the vector $K_S$. The vectors
$\alpha_1,\ldots,\alpha_N$ form a $\bbZ$-basis of $\calR_S$. By
restriction the Weyl group  $W_S$ is isomorphic to a subgroup of
$\Or(\calR_S)$. A choice of a geometric basis $\alpha_1,\ldots,\alpha_N$
defines an isomorphism from $\calR_S$ to the root lattice $Q$ of a
finite root system  of type $E_N (N= 6,7,8), D_5 (N= 5), A_4 (N =
4)$ and $A_2+A_1 (N = 3)$. The group $W_S$ becomes isomorphic to the
corresponding Weyl group $W(E_N)$.

The next lemma is  well-known and its proof goes back to Kantor \cite{Kantor} and Du Val \cite{DuVal}.  We refer for modern proofs to \cite{Alberich} or  \cite{Dolgachev}.

\begin{lemma} Let $(e_0',e'_1,\ldots,e_N')$ be another geometric basis in $\Pic(S)$ defined by a birational morphism $\pi':S\to \bbP^2$ and a choice of a factorization of $\pi'$ into a composition of blowups of  points.  Then the transition matrix is an element of $W_N$. Conversely, any element of $W_N$ is a transition matrix of two geometric bases in $\Pic(S)$.
\end{lemma}

The next lemma is also well-known and is left to the reader.

\begin{lemma} If $d\le 5$, then  the natural  homomorphism
$$\rho:\Aut(S)\to W_S$$
is injective.
\end{lemma}

We will use the known classification of conjugacy classes  in the Weyl groups. According to \cite{Carter} they are indexed by certain graphs. We call them {\it Carter graphs}. One writes each element $w\in W$ as the product of two involutions $w_1w_2$, where each involution is the product of reflections with respect to orthogonal roots. Let $\calR_1, \calR_2$ be the corresponding sets of such roots. Then the graph has vertices identified with elements of the set $\calR_1\cup \calR_2$ and two vertices $\alpha,\beta$ are joined by an edge if and only if  $(\alpha,\beta) \ne 0$. A Carter graph with no cycles is a Dynkin diagram. The subscript in the notation of a Carter graph indicates the number of vertices. It is also equal to the difference between the rank of the root lattice $Q$ and the rank of its fixed sublattice $Q^{(w)}$.

\begin{table}[h]
\begin{center}
\begin{tabular}{| l |r| r | r }
\hline
Graph&Order &Characteristic polynomial\\ \hline
$A_k$&$k+1$&$t^k+t^{k-1}+\cdots+1$\\ \hline
$D_k$&$2k-2$&$(t^{k-1}+1)(t+1)$\\ \hline
$D_k(a_1)$&l.c.m$(2k-4,4)$&$(t^{k-2}+1)(t^2+1)$\\ \hline
$D_k(a_2)$&l.c.m$(2k-6,6)$&$(t^{k-3}+1)(t^3+1)$\\ \hline
\vdots&\vdots&\vdots\\ \hline
$D_k(a_{\frac{k}{2}-1})$&even\  $k$ &$(t^{\frac{k}{2}}+1)^2$\\ \hline
$E_6$&12&$(t^4-t^2+1)(t^2+t+1)$\\ \hline
$E_6(a_1)$&9&$t^6+t^3+1$\\ \hline
$E_6(a_2)$&6&$(t^2-t+1)^2(t^2+t+1)$\\ \hline
$E_7$&18&$(t^6-t^3+1)(t+1)$\\ \hline
$E_7(a_1)$&14&$t^7+1$\\ \hline
$E_7(a_2)$&12&$(t^4-t^2+1)(t^3+1)$\\ \hline
$E_7(a_3)$&30&$(t^5+1)(t^2-t+1)$\\ \hline
$E_7(a_4)$&6&$(t^2-t+1)^2(t^3+1)$\\ \hline
$E_8$&30&$t^8+t^7-t^5-t^4-t^3+t+1$\\ \hline
$E_8(a_1)$&24&$t^8-t^4+1$\\ \hline
$E_8(a_2)$&20&$t^8-t^6+t^4-t^2+1$\\ \hline
$E_8(a_3)$&12&$(t^4-t^2+1)^2$\\ \hline
$E_8(a_4)$&18&$(t^6-t^3+1)(t^2-t+1)$\\ \hline
$E_8(a_5)$&15&$t^8-t^7+t^5-t^4+t^3-t+1$\\ \hline
$E_8(a_6)$&10&$(t^4-t^3+t^2-t+1)^2$\\ \hline
$E_8(a_7)$&12&$(t^4-t^2+1)(t^2-t+1)^2$\\ \hline
$E_8(a_8)$&6&$(t^2-t+1)^4$\\ \hline

\end{tabular}
\end{center}
\label{tablecarter}
\caption{Carter graphs and characteristic polynomials}\label{tab1}
\end{table}
Note that the same conjugacy classes may correspond to different graphs (e.g. $D_3$ and $A_3$, or $2A_3+A_1$ and $D_4(a_1)+3A_1$).

The Carter graph  determines the characteristic polynomial of $w$. In particular, it  gives the trace $\Tr_2(g)$ of $g^*$ on the cohomology space $H^2(S,\bbC) \cong \Pic(S)\otimes \bbC$. The latter should be compared with the Euler-Poincar\`e characteristic of the fixed locus $S^g$ of $g$ by applying the Lefschetz fixed-point formula.

\begin{equation}\label{trace}
\Tr_2(g) = s-2+\sum_{i\in I}(2-2g_i),
\end{equation}
where $S^g$ the  disjoint union of  smooth curves $R_i, i\in I,$ of genus $g_i$ and $s$ isolated fixed points.

To determine whether a finite subgroup $G$ of $\Aut(S)$ is minimal, we use the well-known formula from the character theory of finite groups

$$\rank~\Pic(S)^G = \frac{1}{\#G}\sum_{g\in G}\Tr_2(g).$$

The tables for conjugacy classes of elements from the Weyl group $W_S$ give the values of the trace on the lattice $\calR_S = K_S^\perp$. Thus the group is minimal if and only if the sum of the traces add up to 0.

\subsection{Del Pezzo surfaces of degree $6$}

Let $S$ be a Del Pezzo surface of degree $6$. We fix a geometric
basis  $e_0,e_1,e_2,e_3$ which is defined by a birational morphism $\pi:S\to \bbP^2$
with indeterminacy points $p_1 = (1,0,0), p_2 = (0,1,0)$ and $p_3 = (0,0,1)$. The vectors
$$\alpha_1 = e_0-e_1-e_2-e_3, \alpha_2 =e_1-e_2,\alpha_3= e_2-e_3)$$
form a basis of the lattice $\calR_S$ with Dynkin diagram of type $A_2+A_1$.  The Weyl group
$$W_S = \la s_{1}\ra\times \la s_{2},s_{3}\ra\cong 2\times S_3.$$
The representation $\rho:\Aut(S) \to W_S$ is surjective. The reflection $s_{1}$ is realized by the lift of the standard quadratic transformation $\tau_1$.  The reflection $s_2$( resp. $s_{3}$) is realized by the projective transformations $[x_1,x_0,x_2]$ (resp. $[x_0,x_2,x_1])$. The kernel of $\rho$ is the maximal torus $T$ of $\PGL(3)$, the quotient of $(\bbC^*)^3$ by the diagonal subgroup $\bbC^*$. Thus
$$\Aut(S) \cong T:(S_3\times 2) \cong N(T):2,$$
where $N(T)$ is  the normalizer of $T$ in $\PGL(2)$. It is easy to check that $s_1$ acts on $T$ as the inversion automorphism.

Let $G$ be a minimal finite subgroup of $\Aut(S)$. Obviously, $\rho(G)$ contains $s_1$ and $s_2s_3$ since  otherwise $G$ leaves invariant $\alpha_1$ or  one of the vectors $2\alpha_1+\alpha_2$, or $\alpha_1+2\alpha_2$. This shows that $G\cap N(T)$ is an imprimitive subgroup of $\PGL(3)$. This gives

\begin{theorem} Let $G$ be a minimal subgroup of a Del Pezzo surface of degree $6$. Then
$$G = H_\bullet  \la s_1\ra,$$
where $H$ is  an imprimitive finite subgroup of $\PGL(3)$.
\end{theorem}

Note that one of the groups from the theorem is the group  $2^2:S_3\cong S_4$. Its action on  $S$ given by the equation
$$x_0y_0z_0-x_1y_1z_1 = 0$$
in $(\bbP^1)^3$ is given in \cite{BT}.

\subsection{Del Pezzo surfaces of degree $d = 5$.}

In this case  $S$ is isomorphic to the blowup of the reference
points $p_1 = (1,0,0),\  p_2 = (0,1,0), \ p_3 = (0,0,1), p_4 =
(1,1,1)$. The lattice  $\calR_S$ is of type $A_4$ and $W_S \cong
S_5$. It is known that the homomorphism $\rho:\Aut(S) \to W_S$ is an
isomorphism. We already know that it is injective. To see the
surjectivity one can argue, for example, as follows.

 Let $\tau$ be the standard  quadratic transformation with base points  $p_1,p_2,p_3$. It follows from its formula that the point $p_4$ is a fixed point. We know that $\tau$ can be regularized  on the Del Pezzo surface $S'$  of degree 6 obtained by the blowup of the first three points. Since the preimage of $p_4$ in $S'$ is a fixed point,  $\tau$ lifts to an automorphism of  $S$. Now let $\phi$ be a projective transformation such that $\phi(p_1) = p_1, \phi(p_2) = p_2, \phi(p_4) = p_3$. For example, we take
 $A = [t_0-t_2,t_1,-t_2,-t_2]$. Then the quadratic transformation $\phi^{-1} \tau \phi$ is not defined at the points $p_1,p_2,p_4$ and fixes the point $p_3$. As above, it can be lifted to an involution of $S$. Proceeding in this way we define 4 involutions $\tau = \tau_1,\ldots,\tau_4$ of $S$ each fixes one of the exceptional curves. One checks that their images in the Weyl group $W_S$ generate the group.

Another way to see the isomorphism $\Aut(S)\cong S_5$ is to use a well-known isomorphism between $S$ and the moduli space $\overline{\calM_{0,5}} \cong (\bbP^1)^5/\!/\SL(2)$. The group $S_5$ acts by permuting the factors.

\begin{theorem} Let $(S,G)$ be a minimal Del Pezzo surface of degree $d = 5$. Then $G = S_5, A_5, 5:4, 5:2$, or $C_5$.
\end{theorem}

\begin{proof} As we have just shown $\Aut(S) \cong W_4 \cong S_5$. The group $S_5$ acts on $\calR_S \cong \bbZ^4$ by means of its standard irreducible 4-dimensional representation (view $\bbZ^4$ as a subgroup of $\bbZ^5$ of vectors with coordinates added up to zero and consider the representation of $S_5$ by switching the coordinates). It is known that a maximal proper subgroup of $S_5$ is equal (up to a conjugation) to one of three subgroups $S_4, S_3\times 2, A_5,  5:4$. A maximal subgroup of $A_5$ is either $5\times 2$ or $S_3$ or $D_{10} = 5:2$. It is easy to see that the groups $S_4$ and $S_3\times 2$ have invariant elements in the lattice $Q_4$. It is known that an element of order 5 in $S_5$ is a cyclic permutation, and hence has no invariant vectors. Thus any subgroup $G$ of $S_5$ containing an element of order 5 defines a minimal surface $(S,G)$. So, if $(S,G)$ is minimal, $G$ must be equal to one of the groups from the assertion of the theorem.
\end{proof}

\subsection{Automorphisms of a Del Pezzo surface of degree $d = 4$}

 In this case $\calR$ is of type $D_5$ and $W_S \cong 2^4:S_5$. We use the following well-known classical result.

 \begin{lemma} Let $S$ be a Del Pezzo surface of degree 4. Then $S$ is isomorphic to a nonsingular surface of degree 4 in $\bbP^4$ given by equations
 \begin{equation}\label{eqn}
F_1 = \sum_{i=0}^4 T_i^2 = 0, \quad F_2 = \sum_{i=0}^4 a_iT_i^2 = 0,\end{equation}
where all $a_i$'s are distinct.
\end{lemma}

\begin{proof} It is known that a Del Pezzo surface in its anti-canonical embedding is projectively normal. Using Riemann-Roch, one obtains that $S$ is a complete intersection $Q_1\cap Q_2$ of two quadrics. Let $\calP = \lambda Q_1+\mu Q_2$ be the pencil spanned by these quadrics. The locus of singular quadrics in the pencil is a homogeneous equation of degree 5 in the coordinates $\lambda,\mu$. Since $S$ is nonsingular, it is not hard to see that the equation has no multiple roots (otherwise $\calP$ contains a reducible quadric or there exists a quadric in the pencil with singular point at $S$, in both cases  $S$ is singular).  Let $p_1,\ldots,p_5$ be the singular points of singular quadrics from the pencil. Suppose they  are contained in a hyperplane $H$. Since no quadrics in the pencil contains $H$, the restriction $\calP|H$ of  the pencil of quadrics to $H$ contains $\ge 5$ singular members. This implies that all  quadrics in $\calP|H$ are singular. By Bertini's theorem, all quadrics are singular at some point  $p\in H$. This implies that all quadrics in $\calP$ are tangent to $H$ at $p$. One of the quadrics must be singular at $p$, and hence $S$ is singular at $p$. This contradiction shows that $p_1,\ldots,p_5$ span $\bbP^4$.  Choose coordinates in $\bbP^4$ such that the singular points of singular quadrics from $\calP$ are the points $(1,0,0,0,0), (0,1,0,0,0),$ and so on. Then each  hyperplane $V(t_i)=(t_i=0)$ is a common tangent hyperplane of quadrics from $\calP$ at the point $p_i$. This easily implies that the equations of quadrics are given by \eqref{eqn}.
\end{proof}

Let $Q_i = V(a_iF_1-F_2), i = 0,\ldots, 4,$ be one of the singular
quadrics in the pencil $\calP$. It is a cone over
a nonsingular quadric in $\bbP^3$, hence it contains two families of
planes. The intersection of a plane with any other quadric in the
pencil is a conic contained in $S$.  Thus each $Q_i$ defines a pair
of pencils of conics $|C_i|$ and $|C_i'|$, and it is easy to see
that $|C_i+C_i'| = |-K_S|$.

\begin{proposition} Let $S$ be a Del Pezzo surface given by equations \eqref{eqn}. The divisor classes
$c_i = [C_i]$ together with $K_S$ form a basis of  $\Pic(S)\otimes \bbQ$. The Weyl group $W_S$ acts on this basis by permuting the  $c_i$'s and sending $c_i$ to $c_i' = [C_i'] = -K_S-c_i$.
\end{proposition}

\begin{proof}  If we choose a geometric  basis $(e_0,e_1,\ldots,e_5)$ in $\Pic(S)$, then the 5 pairs of pencils of conics are defined by the classes $e_0-e_i, 2e_0-e_1-\cdots -e_5+e_i$. It is easy to check that the classes $[C_i]$'s and $K_S$ form a basis in $\Pic(S)\otimes \bbQ$. The  group $W_S$ contains a subgroup  isomorphic to $S_5$ generated by the reflections in vectors $e_1-e_2,\ldots,e_4-e_5.$, It acts by permuting $e_1,\ldots,e_5$,  hence permuting the pencils $|C_i|$. It is equal to the semi-direct product of $S_5$ and  the subgroup isomorphic to $2^4$ which is generated by the conjugates of the product $s$ of two commuting reflections with respect to the vectors $e_0-e_1-e_2-e_3$ and $e_4-e_5$. It is easy to see that
$s([C_4]) = [C_4'], s([C_5]) = [C_5']$ and $s([C_i]) = [C_i]$ for $i\ne 4,5$.
This easily implies that $W_S$ acts by permuting the classes $[C_i]$ and switching even number of them to $[C_i']$.
\end{proof}

\begin{corollary}\label{act} Let $W(D_5)$ act in $\bbC^5$ by permuting the coordinates and switching the signs of even number of coordinates. This linear representation of $W(D_5)$ is isomorphic to the representation of $W(D_5)$ on   $\calR_S\otimes \bbC$.
\end{corollary}

The group of projective automorphisms generated by the
transformations which switch $x_i$ to $-x_i$ generates  a subgroup
$H$ of $\Aut(S)$ isomorphic to $2^4$. We identify the group $H$ with the linear space of  subsets of  even cardinality of the set
 $J = \{0,1,2,3,4\}$, or, equivalently, with the subspace of $\bbF_2^J$ of functions with support at a subset of even cardinality. We equip $H$ with the symmetric bilinear form defined by the dot-product in  $\bbF_2^J$, or, equivalently, by  $(A,B) = \#A\cap B \mod 2$. We denote elements of $H$ by $i_A$, where $i_A$ is the characteristic function of $A\subset J$.

There are two kinds of involutions $i_A$.   An involution of the first
kind corresponds to a subset $A$  of  4 elements. The
set of fixed points of such an involution is a hyperplane section of $S$, an
elliptic curve. The trace formula \eqref{trace} gives that the  the trace of $i_A$
in $\Pic(S)$  is equal to $-2$. The corresponding conjugacy class in
$W_5$ is of type $4A_1$. There are 5 involutions of the first kind.
The quotient surface $S/\la i_A\ra = Q$ is isomorphic to a nonsingular
quadric. The map $S\to Q$ coincides with the map $S\to
\bbP^1\times\bbP^1$ that is given by the pencils $|C_i|$ and
$|C_i^\prime|$.

Involutions of the second type correspond to  subsets $A$ of cardinality $2$. The fixed-point set of such involution consists of 4 isolated points. This gives that the trace is equal to 2, and the conjugacy class is of type $2A_1$. The quotient $S/(i_A)$ is isomorphic to the double cover of $\bbP^2$ branched along the union of two conics.

The subgroup of the Weyl group $W(D_5)$ generated by involutions from the conjugacy class of type $2A_1$ is the normal subgroup $2^4$ in the decomposition $W(D_5) \cong 2^4:S_5$. The product of two commuting involutions from this conjugacy class is an involution of type $4A_1$. Thus the image of $H$ in $W_S$ is a normal subgroup isomorphic to $2^4$.

Let $G \cong 2^a$ be a subgroup of $2^4$. All cyclic groups $G$ are not minimal.

There are three kinds of subgroups H of order 4 in $2^4$. A subgroup of the first kind does not contain an involution of the first kind. An example is the group generated by $i_{01},i_{12}$. The trace of its nonzero elements equal to 1. So this group is not minimal.

A subgroup of the second type contains only one involution of the first kind. An example is the group
generated by $i_{01},i_{23}.$ The trace formula gives $\rank~\Pic(S)^H = 2$. So it is also  nonminimal.

A subgroup of the third kind contains two involutions of the first kind. For, example a group generated by $i_{1234},i_{0234}$. It contains 2 elements with trace $-3$ and one element with trace $1$. Adding up the traces we see the group is a minimal group. It is easy to see that
$S^H$ consists of 4 isolated points.

Now let us consider subgroups of order 8 of $2^4$. They are parametrized by the same sets which parametrize involutions.  A subgroup $H_A$ corresponding to a subset $A$ consists of involutions $i_B$ such that $\#A\cap B$ is even. The subsets $A$ correspond to linear functions on $2^4$.  If $\#A = 2$, say $A = \{0,1\}$, we see that $H_A$ contains the
involutions $i_{01}, i_{01ab}, i_{cd}, c,d\ne 0,1$. Adding up the traces we obtain that  these subgroups are minimal.

If $\#A = 4$, say $A = \{1,2,3,4\}$, the subgroup $H_A$ consists of $i_{1234}$ and $i_{ab}$, where $a, b \ne 0$. Adding up the traces we obtain that  $H_A$ is not minimal.

Since $2^4$ contains a minimal subgroup, it is minimal itself.

Now suppose that  the image $G'$ of  $G$ in $S_5$ is non-trivial. The subgroup $S_5$ of $\Aut(S)$ can be realized as the stabilizer of a set of 5 skew lines on $S$ (corresponding to the basis vectors $e_1,\ldots,e_5$). Thus any subgroup $H$ of $S_5$ realized as a group of automorphisms of $S$ is isomorphic to a group of projective transformations of $\bbP^2$ leaving invariant a set of 5 points. Since there is a unique conic through these points, the group is isomorphic to a finite group of $\PGL(2)$ leaving invariant a set of 5 distinct points. In section \ref{sec3}, we listed all possible subgroups of $\GL(2)$ and in section \ref{sec4} we described their relative invariants. It follows that a subgroup leaves invariant a set of 5 distinct points if and only if it is one of the following groups
$C_2,C_3,C_4, C_5, S_3, D_{10}$.  The corresponding binary forms of degree 5 are projectively equivalent to the following binary forms:
\begin{itemize}
\item $C_2: t_0(t_0^2-t_1^2)(t_0^2+at_1^2), \ a\ne -1,0,1$;
\item $C_4: t_0(t_0^2-t_1^2)(t_0^2+t_1^2)$;
\item $C_3,S_3: t_0t_1(t_0-t_1)(t_0-\epsilon_3 t_1)(t_0-\epsilon_3^2t_1)$;
\item $C_5,D_{10}: (t_0-t_1)(t_0-\epsilon_5 t_1)(t_0-\epsilon_5^2t_1)(t_0-\epsilon_5^3t_1)(t_0-\epsilon_5^4t_1)$.
\end{itemize}
The corresponding surfaces are projectively equivalent to the following surfaces
\begin{eqnarray}\label{dp4}
C_2&: &x_0^2+x_1^2+x_2^2+x_3^2+x_4^2 =
x_0^2+ax_1^2-x_2^2-ax_3^2= 0, \ a\ne -1,0,1,\\
 C_4&:&x_0^2+x_1^2+x_2^2+x_3^2+x_4^2 = x_0^2+ix_1^2-x_2^2-ix_3^2 = 0,\\
S_3&:&x_0^2+\epsilon_3 x_1^2+\epsilon_3^2x_2^2+x_3^2 =
x_0^2+\epsilon_3^2x_1^2+\epsilon_3 x_2^2+x_4^2 = 0,\\ 
D_{10}&: &\sum_{i=0}^4\epsilon_5^ix_i^2=\sum_{i=0}^4
\epsilon_5^{4-i}x_i^2 = 0 \label{X}
\end{eqnarray}

\begin{remark} \label{rmknew} Note that equations (6.4),  (6.5) and (6.6)  are  specializations of equation (6.3) . It is obvious for equation (6.4) where we have to take $a = i$. Equation (6.3)  specializes to equation (6.5) when we take $a= \pm \frac{1}{\sqrt{-3}}$ (use that the Moebius transformation of order 3 $x\mapsto \frac{ax+1}{x+a}$ permutes cyclically $\infty,a,-a$ and fixes $1,-1$). Equation (6.3) specializes to equation (6.6) if we take $a = -2\pm \sqrt{5}$ (use that the Moebius transformation $x\mapsto \frac{x+2a-1}{x+1}$ permutes cyclically $(\infty,1,a,-a,-1)$). We thank J. Blanc for this observation.
\end{remark}

Since the subgroup $S_5$ leaves the class $e_0$ invariant,  it remains to consider subgroups $G$ of $2^4:S_5$ which are not contained in $2^4$ and not conjugate to a subgroup of  $S_5$. We use the following  facts.

1) Suppose $G$ contains a minimal subgroup of $2^4$. Then $G$ is minimal.

2) Let $\overline{G}$ be the image of $G$ in $S_5$. Then it is a subgroup of one of the groups listed above.

3) The group $W(D_5)$ is isomorphic to the group of transformations
of $\bbR^5$ which consists of permutations of coordinates and
changing even number of signs of the coordinates. Each element $w\in
W(D_5)$ defines a permutation of the coordinate lines which can be
written as a composition of cycles $(i_1\ldots i_k)$. If $w$ changes
signs of even number of the  coordinates $x_{i_1},\ldots,x_{i_k}$, the cycle
is called positive. Otherwise it is called a negative cycle. The
conjugacy class of $w$ is determined by the lengths of positive  and
negative cycles, except when all cycles of even length and positive
in which case there are two conjugacy classes.  The latter case does
not occur in the case when $n$ is odd. Assign to a positive cycle of
length $k$ the Carter graph $A_{k-1}$. Assign to a pair of negative
cycles of lengths $i\ge j$ the Carter graph of type $D_{i+1}$ if $j
= 1$ and $D_{i+j}(a_{j-1})$ if $j > 1$. Each conjugacy class is
defined by the sum of the graphs.  We identify $D
 _2$ with $2A_1$, and $D_3$ with $A_3$. In Table 2 below we give the conjugacy classes of elements in $W(D_5),$ their characteristic polynomials and the traces in the root lattice of type $D_5$.

\bigskip

\begin{table}[ht]
\centering
{\small \begin{tabular}{||r |r |r |r |r ||}\hline
Order&Notation&Characteristic  polynomial&Trace&Representatives\\ \hline
2&$A_1$&$t+1$&3&$(ab)$\\ \hline
2&$2A_1$&$(t+1)^2$&1&$(ab)(cd), (ab)(cd)i_{abcd}$\\ \hline
2&$2A_1^*$&$(t+1)^2$&1&$(ab)(cd), (ab)(cd)i_{ab}$\\ \hline
2&$3A_1$&$(t+1)^3$&-1&$(ab)i_{cd}$\\ \hline
2&$4A_1$&$(t+1)^4$&-3&$i_{abcd}$\\ \hline
3&$A_2$&$t^2+t+1$&2&$(abc), (abc)i_{ab}$\\ \hline
4&$A_3$&$t^3+t^2+t+1$&1&$(abcd), (abcd)i_{ab}, (abcd)i_{abcd}$\\ \hline
4&$A_1+A_3$&$(t^3+t^2+t+1)(t+1)$&-1&$(ab)(cd)i_{ae}$\\ \hline
4&$D_4(a_1)$&$(t^2+1)^2$&1&$(ab)(cd)i_{ac}$\\ \hline
5&$A_4$&$(t^4+t^3+t^2+t+1)$&0&$(abcde),(abcde)i_{A}$\\ \hline
6&$A_2+A_1$&$(t^2+t+1)(t+1)$&0&$(ab)(cde)$\\ \hline
6&$A_2+2A_1$&$(t^2+t+1)(t+1)^2$&-2&$(abc)i_{abde}, (abc)i_{de}$\\ \hline
6&$D_4$&$(t^3+1)(t+1)$&0&$(abc)i_{abce}$\\ \hline
8&$D_5$&$(t^4+1)(t+1)$&-1&$(abcd)i_{abce}, (abcd)i_{de}$\\ \hline
12&$D_5(a_1)$&$(t^3+1)(t^2+1)$&0&$(abc)(de)i_{ac}$\\ \hline
\end{tabular}}
\bigskip
\caption{Conjugacy classes in $W(D_5)$}\label{tab2}
\end{table}

In the following $\overline{G}$ denotes the image of $G$ in $K = W(D_5)/2^4\cong S_5$.
\smallskip
Case 1. $\overline{G} \cong C_2$.

It follows from the description  of the image of
$\Aut(S)$ in $W(D_5)$ given in Corollary \ref{act}, that $\overline{G}$ is generated by the permutation  $s
= (02)(13).$ Let $g\not\in G\cap 2^4$. Then  $g = s$ or $g = s
i_A$ for some $A$. It follows from Table \ref{tab2} that $g$ is either of type $2A_1$,  or of type $A_1+A_3$, or of type $D_4(a_1)$.  Let $K = G\cap 2^4 \cong 2^a$. If $a = 0$ or $1$, the group is not minimal.

\smallskip
 $a = 2$.

Suppose first that $s$ acts identically on $K$. Then the group is commutative isomorphic to $2^3$ if it does not contain elements of order 4 and $2\times 4$ otherwise. In the first case
\begin{equation}\label{setF}
K =  F:= \la i_{02},i_{13}\ra.
\end{equation}
is the subspace of fixed points of $s$ in $2^4$. Since $(si_A)^2 = i_{A+s(A)}$, we see that $G\setminus K = \{si_A, A\in K\}$. Consulting Table \ref{tab2} we compute the traces of all elements from $G$ to conclude that the total sum is equal to 8. Thus the group is not minimal.

In the second case $G$ contains 4 elements of order 4 of the form  $si_A$, where $A\not \in K$. Suppose $A+s(A) = \{0,1,2,3\}$. Since $K$ is a subspace of the second type, it contributes 4 to the total sum of the traces. Thus the sum of the traces of the elements of order 4 must be equal to $-4$. In other words they have to be elements with trace $-1$ of the form $si_A$, where $\#A = 2, 4\in A$. This gives the unique conjugacy class of a minimal group isomorphic to $2\times 4$. It is represented by the group
$G = \la K,si_{04}\ra$.

Assume now that $G$ is  non-abelian, obviously isomorphic to $D_8$. The subspace $K$ contains one element from the set $F$.  The nontrivial coset contains 2 elements of order 4  and two elements of order 2. Suppose $K$ is of the second type with the sum of the traces of its elements equal to 4.  Two elements of order 2 in $G\setminus K$ have the trace equal to 1. Elements of order 4 have the trace equal to 1 or $-1$. So the group cannot be minimal. Thus $K$ must be of the third type, the minimal one. This gives us the minimal group conjugate to the subgroup
$G = \la i_{1234},i_{02},si_{04}\ra$ isomorphic to $D_8$.

\smallskip
 $a = 3$.

There are three $s$ invariant subspaces of $2^4$ of dimension 3. Their orthogonal complements are spanned by the one of the vectors in the set \eqref{setF}. As we saw earlier, if $K^\perp =\la i_A\ra,$ where $\#A = 2$, the subspace $K$ is a minimal group. Otherwise the total sum of the traces of elements from $K$ is equal to $8$. In the first case  we may assume that  $K = \la i_{14},i_{34},i_{02}\ra$. All elements of order 2 in  the nontrivial coset have the trace equal to 1. Thus we must have elements of order 4 in the coset with trace $-1$. Let $si_A$ be such an element, where we may assume that $A = \{0,4\}$. Thus $G = \la K,si_{04}\ra$. Its nontrivial coset has 4 elements of order 4 with trace $1$ and 4 elements of order 4 with trace $-1$. The group is minimal. It is a non-split extension $(2^3)^\bullet 2$. Its center is isomorphic to $2^2$. The classification of groups of order 16 from Table \ref{order8} shows that this is group is isomorphic to $L_{16}$.

\smallskip
$a = 4$.

In this case $G = 2^4:2$, where the extension is defined by the action of $s$ in $2^4$. The group has 2-dimensional center with the quotient isomorphic to $2^3$.

\smallskip
Case 2. $\overline{G} \cong C_3$.

We may assume that $\overline{G} =\la s\ra $, where $s = (012)$.  Applying Lemma \ref{zas}, we obtain that $G$ is a split extension
$K:3$, where $K= G\cap 2^4 \cong 2^a$.  Since there are no minimal elements of order 3, we must have $a> 0$. If $a = 1$, the group is $2:3 \cong 6$. There are no minimal elements of order 6, so we may assume that $a > 1$.

Assume $a = 2$. The group is abelian $2^2\times 3$ or non-abelian  $2^2:3 \cong A_4$. In the first case, $K = \la i_{0123},i_{0124}\ra$ is the subspace of the third type, the minimal one. Thus the total number of elements in the nontrivial cosets is equal to $0$. An element of order 3 has the trace equal to $1$. An element of order 6 has trace equal to $-2$ of $0$. So we must have an element of order 6 with trace $-2$. It must be equal to  $si_{34}$. Its cube is $i_{34}\in K$. So we get one conjugacy class of a minimal group isomorphic to $2^2\times 3$. It is equal to $\la K,si_{34}\ra$.

If $G \cong A_4$, the subspace $K$ is not minimal. The group does not contain elements of order 6. So the traces of all elements not from $K$ are positive. This shows that the group cannot be minimal.

Assume $a = 3$. The subspace $K$ is minimal if and only if its orthogonal complement is generated by $i_{34}$. Again the group must contain $si_{34}$  with trace $-2$ and hence equal to $\la K,si_{34}\ra = \la K,s\ra$.
 The group $G$ is minimal and is isomorphic to
$2\times (2^2:3) \cong  2\times A_4.$

Finally assume that $K = \la i_{0123}\ra^\perp$ is not minimal. We have computed earlier the sum of the traces of its elements. It is equal to 8. Again it must contain an element of order 6 equals  $si_{0123}$. Since $K$ contains $i_{0123}$, the group contains $s$. Now we can add all the traces and conclude that the group is not minimal.

Of course we should not forget the minimal group $2^4:3$.

\smallskip
Case 3. $\overline{G} \cong  S_3$.

The group  $\overline{G}$ is generated by the permutations of coordinates $(012)$ and $(12)(34)$. It is immediately checked that $H = G\cap 2^4$ is not trivial for minimal $G$. The only subgroup of $H$ invariant with respect to the conjugation action of $\overline{G}$ on $H$ is $H$ itself. This gives a minimal group isomorphic to $2^4:S_3$. The extension is defined by the restriction to $S_3$ of the natural action of $S_4 = W(A_4)$ on its root lattice modulo 2.

\smallskip
Case 4. $\overline{G} \cong  C_4$.

The group $2^4:4$ contains $2^4:2$, so all minimal groups of the
latter group are minimal subgroups of $2^4:4$. Without loss of
generality, we may assume that the group  $\overline{G}$ is generated by
the permutations of coordinates $s = (0123)$.   The only proper
subgroup of $2^4$ invariant with respect to the conjugation action
of $\overline{G}$ on $K = G\cap 2^4$ is either  $\la i_{0123}\ra$ or its orthogonal complement. In the first case $G \cong  2\pt 4 \cong 2\times 4$ or $8$. In the first case the group is not minimal. In the second case  $G = \la (0123)i_{0123}\ra $ is minimal.

Assume $K = \la i_{0123}\ra^\perp$.
If $s\not \in G$, then $si_{A}\in G$, where $A\not \in K$. The Table \ref{tab2} shows that all such elements are minimal of order 8.  This gives a minimal group $G \cong 2^3:4 = 2^2:8$.

Next we have to consider the case when $s \in G$ so that $G = \la K,s\ra$. The total number of traces of elements from $K$ is equal to 8. Consulting the Table \ref{tab2} we obtain that the elements in the cosets $sK, s^2K, s^3K$ have the trace equal to 1. So the group is not minimal.

Our last minimal group in this case is $2^4:4$.

\smallskip
Case 5. $\overline{G} = C_5$ or $D_{10}$.

Again, we check using the table of conjugacy classes that no group isomorphic to $D_{10}$ is minimal. Also no proper subgroup of $H$ is invariant with respect to conjugation by a permutation of order 5, or by a subgroup of $S_5$ generated by $(012)$ and $(12)(34)$.  Thus we get two minimal groups isomorphic to $2^4:5$ or $2^4:D_{10}$.

The following theorem summarizes what we have found.

\begin{theorem}\label{clasdp4} Let $(S,G)$ be a minimal Del Pezzo surface of degree $d = 4$. Then $G$ is isomorphic to one of the following groups
\begin{enumerate}
\item $\Aut(S) \cong  2^4$.
$$2^4, \ 2^3, \ 2^2.$$
\item $\Aut(S) \cong 2^4:2$
$$ 2\times 4, \ D_8,\ L_{16},\ 2^4:2,
$$
and from the previous case.
\item $\Aut(S) \cong 2^4:4$
$$  \ 8, \ 2^2:8,\  2^4:4,$$
and from the previous two cases.
\item $\Aut(S) \cong  2^4:S_3$.
$$ 2^2\times 3,\  2\times A_4, \quad 2^4:3, \ 2^4:S_3,$$
and from Cases 1) and 2).
\item $\Aut(S) \cong  2^4:D_{10}$
$$ 2^4:D_{10},\  2^4:5,
$$
and from Cases 1) and 2).
\end{enumerate}
\end{theorem}

\subsection{Cubic surfaces}

The following theorem gives the classification of cyclic subgroups of $\Aut(S)$ and identifies the conjugacy classes of their generators.

\begin{theorem}\label{clas3} Let $S$ be a nonsingular cubic surface admitting a non-trivial automorphism $g$ of order $n$. Then $S$ is equivariantly isomorphic to one of the following surfaces $V(F)$ with
\begin{equation}\label{dact}
g = [t_0,\epsilon_n^at_1,\epsilon_n^bt_2,\epsilon_n^ct_3].\end{equation}
\begin{itemize}
\item $4A_1$ $(n = 2)$,  $(a,b,c)= (0,0,1)$,
$$F = T_3^2L_1(T_0,T_1,T_2)+T_0^3+T_1^3+T_2^3+\alpha T_0T_1T_2.$$
\item $2A_1$ $(n = 2)$, $(a,b,c)= (0,1,1),$
$$F = T_0T_2(T_2+\alpha T_3)+T_1T_3(T_2+\beta T_3)+T_0^3+T_1^3.$$
\item $3A_2$ $(n = 3)$, $(a,b,c)= (0,0,1),$
$$F = T_0^3+T_1^3+T_2^3+T_3^3+\alpha T_0T_1T_2.$$
\item $A_2$ $(n = 3)$, $(a,b,c)= (0,1,1),$
$$F = T_0^3+T_1^3+T_2^3+T_3^3.$$
\item $2A_2$ $(n = 3)$, $(a,b,c)= (0,1,2),$
$$F = T_0^3+T_1^3+T_2T_3(T_0+aT_1)+T_2^3+T_3^3.$$
\item $D_4(a_1)$ $(n = 4)$, $(a,b,c)= (0,2,1),$
$$F = T_3^2T_2+L_3(T_0,T_1)+T_2^2(T_0+\alpha T_1).$$
\item $A_3+A_1$ $(n = 4)$, $(a,b,c)= (2,1,3),$
$$F = T_0^3+T_0T_1^2+T_1T_3^2+T_1T_2^2.$$
\item $A_4$ $(n = 5)$, $(a,b,c)= (4,1,2)$,
$$F = T_0^2T_1+T_1^2T_2+T_2^2T_3+T_3^2T_0.$$
\item $E_6(a_2)$ $(n = 6)$, $(a,b,c)= (0,3,2),$
$$F = T_0^3+T_1^3+T_3^3+T_2^2(\alpha T_0+T_1).$$
\item $D_4$ $(n = 6)$, $(a,b,c)= (0,2,5),$
$$F = L_3(T_0,T_1)+T_3^2T_2+T_2^3.$$
\item $A_5+A_1$ $(n = 6)$, $(a,b,c)= (4,2, 1),$
$$F = T_3^2T_1+T_0^3+T_1^3+T_2^3+\lambda T_0T_1T_2.$$

\item $2A_1+A_2$ $(n = 6)$, $(a,b,c)= (4,1,3),$
$$F = T_0^3+\beta T_0T_3^2+T_2^2T_1+T_1^3.$$
\item $D_5$ $(n = 8)$, $(a,b,c)= (4,3,2),$
$$F = T_3^2T_1+T_2^2T_3+ T_0T_1^2+T_0^3.$$
\item $E_6(a_1)$ $(n = 9)$, $(a,b,c)= (4,1,7),$
$$F = T_3^2T_1+T_1^2T_2+ T_2^2T_3+T_0^3.$$
\item $E_6$ $(n = 12)$, $(a,b,c)= (4,1,10),$
$$F = T_3^2T_1+T_2^2T_3+T_0^3+T_1^3.$$
\end{itemize}
\end{theorem}

We only sketch a proof, referring for the details to \cite{Topics}. Let $g$ be a nontrivial projective automorphism of $S = V(F)$ of order $n$. All  possible values of $n$ can be obtained   from the classification of conjugacy classes of $W(E_6)$. Choose coordinates to assume that $g$ acts as in \eqref{dact}.  Then $F$ is a sum of monomials which belong to the same eigensubspace of $g$ in its action in the space of cubic polynomials. We list all possible eigensubspaces. Since $V(F)$ is nonsingular, the square or the cube of each variable divides some monomial entering in $F$. This allows one to list all possible nonsingular $V(F)$ admitting an automorphism $g$.  Some obvious linear change of variables allows one to find normal forms. Finally, we determine the conjugacy class by using the trace formula \eqref{trace} applied to the locus of fixed points of $g$ and its powers.

The conjugacy class labeled by the  Carter graph with $6$ vertices defines a minimal cyclic group.

\begin{corollary} The following conjugacy classes define  minimal cyclic groups of automorphisms of a cubic surface $S$.
\begin{itemize}
\item $3A_2$ of order 3,
\item $E_6(a_2)$ of order 6,
\item $A_5+A_1$ of order 6,
\item $E_6(a_1)$ of order 9,
\item $E_6$ of order 12.
\end{itemize}
\end{corollary}

Next we find all possible automorphism groups of nonsingular cubic surfaces. Using a normal form of a cubic admitting a cyclic group of automorphisms from given conjugacy class, we determine all other possible symmetries of the equation. We refer for the details to \cite{Topics}.  The list of possible automorphism groups of cubic surfaces is given in Table \ref{tab3}.

{\small \begin{table}[ht]
\begin{center}
\begin{tabular}{| l|r | r | r | r|}
\hline
\footnotesize{Type}&\footnotesize{Order}&\footnotesize{Structure} &\footnotesize{$F(T_0,T_1,T_2,T_3)$}&\footnotesize{Parameters}\\ \hline
I&648 &$3^3:S_4$&\footnotesize{$T_0^3+T_1^3+T_2^3+T_3^3$}&\\ \hline
II&120 &$S_5$&\footnotesize{$T_0^2T_1+T_0T_2^2+T_2T_3^2+T_3T_1^2$}&\\ \hline
III&108&$H_3(3):4$&\footnotesize{$T_0^3+T_1^3+T_2^3+T_3^3+6aT_1T_2T_3$}&\footnotesize{ $20a^3+8a^6 =1$}\\ \hline
IV&54&$H_3(3):2$&\footnotesize{$T_0^3+T_1^3+T_2^3+T_3^3+6aT_1T_2T_3$}&$a-a^4\ne 0,$\\
&&&&$8a^3\ne -1$,\\
&&&&$20a^3+8a^6 \ne 1$\\ \hline
V&24&$S_4$&\footnotesize{$T_0^3+T_0(T_1^2+T_2^2+T_3^2)$}&$9a^3\ne 8a$\\
&&&\footnotesize{$+aT_1T_2T_3$}&$8a^3\ne -1,$\\
\hline
VI&12&$S_3\times 2$&\footnotesize{$T_2^3+T_3^3 +aT_2T_3(T_0+T_1)+T_0^3+T_1^3$}&$a\ne 0$\\ \hline
VII&8&8&\footnotesize{$T_3^2T_2+T_2^2T_1+T_0^3+T_0T_1^2$}&\\ \hline
VIII&6&$S_3$&\footnotesize{$T_2^3+T_3^3 +aT_2T_3(T_0+bT_1)+T_0^3+T_1^3$}&$a^3\ne -1$\\ \hline
IX&4&$4$&\footnotesize{$T_3^2T_2+T_2^2T_1+T_0^3+T_0T_1^2+aT_1^3$}&$a\ne 0$\\ \hline
X&4&$2^2$&\footnotesize{$T_0^2(T_1+T_2+aT_3)+T_1^3+T_2^3$} \\
&&&\footnotesize{$+T_3^3+6bT_1T_2T_3$}&$8b^3\ne -1$\\
\hline
XI&2&$2$&\footnotesize{$T_1^3+T_2^3+T_3^3+6aT_1T_2T_3$}&$b^3,c^3\ne 1$\\
&&&\footnotesize{$+T_0^2(T_1+bT_2+cT_3)$}&$b^3\ne c^3$\\
&&&&$8a^3\ne -1,$\\
\hline

\end{tabular}
\end{center}
\caption{Groups of automorphisms of cubic surfaces}\label{table1}\label{tab3}
\end{table}}

\begin{remark} Note that there are various  ways to write the equation of cubic surfaces from the table. For example, using the identity
$$(x+y+z)^3+\epsilon_3(x+\epsilon_3y+\epsilon_3^2z)^3+\epsilon_3^2(x+\epsilon_3^2y+\epsilon_3z)^3 = 9(x^2z+y^2x+z^2y)$$
we see that the Fermat cubic can be given by the equation
$$T_0^3+T_1^2T_3+T_3^2T_2+T_2^2T_1 = 0.$$
Using Theorem \ref{clas3} this exhibits a symmetry of order 9 of the surface, whose existence is not obvious in the original equation.

Using the Hesse form of an equation of a nonsingular plane cubic curve we see that a surface  with equation
$$T_0^3+F_3(T_1,T_2,T_3) = 0$$
is projectively equivalent to a surface with the equation
$$T_0^3+T_1^3+T_2^3+T_3^3+6aT_0T_1T_2 = 0.$$
The special values of the parameters $a =0, 1,\epsilon_3,\epsilon_3^2$ give the Fermat cubic. The values $a$ satisfying $20a^3+8a^6 = 1$ give a plane cubic with an  automorphism of order 4 (a harmonic cubic). Since a harmonic cubic is isomorphic to the cubic with equation
$T_1^3+T_1T_2^2+T_3^3 = 0$, using Theorem \ref{clas3} we see symmetries of order 6 from the conjugacy class $E_6(a_2)$  for surfaces of type $III,IV$ and of order 12 for the surface
$$T_3^2T_1+T_2^2T_3+T_0^3+T_1^3 = 0$$
of type $III.$
\end{remark}

It remains to classify minimal groups $G$.  Note that if $G$ is realized as a group of projective (or weighted projective) automorphisms  of a family  of surfaces $(S_t)$, then $G$ is  a subgroup of the group of projective automorphisms of any surface  $S_{t_0}$ corresponding to a special  value $t_0$ of the parameters. We indicate this by writing $S'\to S$.  The types of $S'$ when it happens are
$$IV \to  III, \ IV \to I, \quad   VI,VIII,IX \to I, \quad XI\to X.$$
So it suffices to consider the surfaces of types I, II, III, V, VII, X.

We will be using the following lemma, kindly communicated to us by R. Griess. For completeness sake, we provide its proof.

\begin{lemma}\label{bob} Let $S_{n+1}$ act naturally on its root lattice
$\calR_n = \{(a_1,\ldots,a_{n+1})\in \bbZ^{n+1}:a_1+\cdots+a_{n+1} = 0\}$. Let  $\calR_n(p)\cong \bbF_p^n$  be the reduction of $\calR_n$ modulo a prime number $p > 2$ not dividing $n+1$. Then the set of conjugacy classes of subgroups dividing a splitting of $\bbF_p^n:S_{n+1}$ is bijective to the set $\bbF_p$ if $p\mid n+1$ and consists of one element otherwise.
\end{lemma}

\begin{proof} It is easy to see that, fixing a splitting, there is a natural bijection between conjugacy classes of splitting subgroups in $A:B$ and the cohomology set  $H^1(B,A)$, where $B$ acts on $A$ via the homomorphism $B\to \Aut(A)$ defining the semi-direct product. So, it suffices to prove that $H^1(S_{n+1},\calR_n(p)) \cong \bbF_p $ if $p\mid n+1$ and zero otherwise. Consider the permutation representation of $S_{n+1}$ on $M = \bbF_p^{n+1}$. We have an exact sequence of $S_{n+1}$-modules
$$0\to \calR_n(p) \to M \to \bbF_p \to 0$$
defined by the map $(a_1,\ldots,a_{n+1} )\to a_1+\cdots+a_{n+1}$.  The module $M$ is the induced module of the trivial representation of the subgroup $S_n$ of $S_{n+1}$. By Eckmann-Shapiro's Lemma,
$$H^1(S_{n+1},M) = H^1(S_n,\bbF_p) = \Hom(S_n,\bbF_p).$$
Suppose $p\nmid n+1$, then the exact sequence splits, and we get
$$0 = H^1(S_{n+1},M) = H^1(S_{n+1},\calR_n(p))\oplus H^1(S_{n+1},\bbF_p). $$
Since $H^1(S_{n+1},\bbF_p) = \Hom(S_{n+1},\bbF_p) = 0$, we get the result. If $p\mid n+1$, then
$H^0(S_{n+1},M) = 0, H^0(S_{n+1},\bbF_p) = \bbF_p$ and the exact sequence of cohomology gives the desired result.
\end{proof}

\medskip\noindent
Type I.

Let us first classify  $\bbF_3$-subspaces of the group $K = 3^3$. We view $3^3$ as the $S_4$-module $\calR_3(3) \cong \bbF_3^3$ from the previous lemma.   We denote the image of a vector $(a,b,c,d)$ in $K$ by $[a,b,c,d]$.  In our old notations
$$[a,b,c,d] = [\epsilon^at_0,\epsilon^bt_1,\epsilon^ct_2,\epsilon^dt_3].$$
There are 13 ($=\#\bbP^2(\bbF_3)$) one-dimensional subspaces in $3^3$. The group $S_4$ acts on this set with 3 orbits. They are represented by vectors
$[1,2,0,0], [1,1,1,0], [1,1,2,2]$ with respective stabilizer subgroups $2^2, S_3$ and $D_8$.  We call them lines of type I,II,III, respectively. As subgroups they are cyclic groups of order 3 of the following types.
$$\begin{cases}
      2A_2& \text{Type I}, \\
      3A_2& \text{Type II},\\
      A_2&\text{Type III}.
\end{cases}
$$

The conjugacy class of a  2-dimensional subspace  $K$ is determined by its orthogonal complement in $3^3$ with respect to the dot-product pairing on $\bbF_3^4$. Thus we have 3 types of 2-dimensional subspaces of types determined by the type of its orthogonal complement.

 An easy computation gives the following table.
\begin{equation}\label{subspaces}
\begin{matrix}\text{Type}&3A_2&2A_2&A_2&\text{Trace}\\
\text{I}&4&2&2&0\\
\text{II}&2&6&0&0\\
\text{III}&0&4&4&2
\end{matrix}
\end{equation}
Here we list the types of the nontrivial elements in the subspace and the last column is the sum $\frac{1}{9}\sum_{g\in L}\Tr(g|\calR_S)$. This gives us two conjugacy classes of minimal subgroups isomorphic to $3^2$.

Let $G$ be a subgroup of $\Aut(S)$, $\overline{G}$ be its image in $S_4$ and $K = G\cap 3^3$. Let
$k = \dim_{\bbF_3} K$.

\smallskip
\emph{Case 1}: $k = 0$.

In this case $G$ defines a  splitting of the projection $3^3:S_4\to S_4$. Assume $G \cong S_4$. It follows from Lemma \ref{bob} that there are 3 conjugacy classes of subgroups isomorphic to $S_4$ which define a splitting.  Let us show each of them is minimal.

We start with the standard $S_4$ generated by permutations of the coordinates. It contains 6 elements of type $4A_1$, 8 elements of type $2A_2$, 3 elements of type $2A_1$ and 6 elements of type $A_3+A_1$. Adding up the traces we obtain that the group is minimal.

Suppose $G$ is another subgroup isomorphic but not conjugate to the previous $S_4$. It corresponds to a 1-cocylce $\phi:S_4\to \bbF_3^3$ defined by a vector $v = (a_1,a_2,a_3,a_4)$ with $a_1+a_2+a_3+a_4$
$$\phi_v(\sigma) = (a_1,a_2,a_3,a_4)-(a_{\sigma(1)},a_{\sigma(2)},a_{\sigma(3)},a_{\sigma(4)}) = (a,b,c,d).$$
The cohomology class of this cocycle depends only on the sum of the coordinates of the vector $v$. Without loss of generality we may choose $v = (1,0,0,0)$ and drop the subscript in $\phi$.  Thus a new $S_4$ consists of transformation $\sigma\circ \phi(\sigma)$.  We check that the type of an element $\sigma\phi(\sigma)$ is equal to the type of $\sigma$ for each $\sigma\in S_4$. This shows that the representation of a new $S_4$ in $\calR_S$ has the same character as that of the old $S_4$. This shows that all three $S_4$'s are minimal.

Let $G$ be mapped isomorphically to a subgroup $G'$ of $S_4$. If $G'$  is a 2-group, then it is contained in a 2-Sylow subgroup of one of the $S_4$'s. By the above its representation  in $\calR_S$ is the same as the restriction of the representation of $S_4$. A 2-Sylow subgroup of $S_4$ contains 3 elements of type $2A_1$, two element of type $4A_1$ and 2 elements of type $A_3+A_1$. Adding up the traces we get 8. Thus the group is not minimal.

 If $G'$ is not a 2-group, then $G'$ is either $S_3$ or $3$. A lift of a permutation $(i23)$ is given by a matrix
 $$\begin{pmatrix}0&0&\epsilon_3^a\\
 \epsilon_3^b&0&0\\
 0&\epsilon_3^c&0\end{pmatrix},$$
 where $a+b+c\equiv 0 \mod 3$. It is immediately checked that all such matrices define an element of type $2A_2$.  Adding up the traces we see that a lift of $S_3$ is minimal.

\smallskip
\emph{Case 2}:  $k = 1$.

Let $G'$ be the image of $G$ in $S_4$ and $K = G\cap 3^3$. Clearly, the subspace $K$ must be invariant with respect to the restriction of the homomorphism $S_4\to \Aut(3^3)$ to $G'$.

Assume $K$ is of type $I$, say generated by $\bfn_K = [1,2,0,0]$. Then its stabilizer in $S_4$ is generated by $(12), (34), (123)$ isomorphic to $2^2$. The conjugation by $(12)$ sends $\bfn_K$ to $-\bfn_K$, and the conjugation by $(34)$ fixes $\bfn_K$. Thus the group $K:2^2$ is isomorphic $S_3\times 2$.
It is easy to check that the product $g = (34)\bfn_K$ is of order 6, of minimal type $A_5+A_1$. Thus the  groups  $K:2^2$ and its cyclic subgroup of order 6 are minimal. Also its subgroup of order 6 is minimal. Its subgroup $\la \bfn_K, (12)\ra \cong S_3$ contains $3$ elements of  type $4A_1$ and   2 elements of type $2A_2$.  Adding up the traces we see that $S_3$ is also minimal. It is obviously not conjugate to the $S_3$ from the previous case.

Assume $K$ is of type $II$, say generated by $\bfn_K = [1,1,1,0]$.  Since $\bfn_K$ is minimal, any subgroup in this case is  minimal. The stabilizer of $K$ in $S_4$ is generated by $(123), (12)$ and is isomorphic to $S_3$. Our group $G$ is a subgroup of $3\pt S_3$ with $K$ contained in the center. There are three non-abelian groups of order 18: $D_{18}, 3\times S_3, 3^2:2$. The extension in the last group is defined by the automorphism of $3^2$ equal to the minus identity. In our case the image of $G$ in $S_3$ acts identically on $K$. Since the center of $D_{18}$ or  $3^2:2$ is trivial, this implies that either  $G$ is a cyclic  subgroup of $D_{18}$ of order 9 or 3, or a subgroup of
$S_3\times  3$ in which case $G \cong 3, 3^2, 6, 3\times S_3$.  Note that the group $3^2$ is not conjugate to a subgroup of $3^3$. To realize  a cyclic subgroup of order 9 is enough to take $g = \bfn_K(234)$. Note that the Sylow subgroup of $3\pt S_3$ is of order 9, so all  3-subgroups of order 9 are conjugate.

Assume $K$ is of type $III$, say generated by $\bfn_K = [1,1,2,2]$. The stabilizer group is generated by $(12), (34), (13)(24)$ and is isomorphic to $D_8$. Our group $G$ is a subgroup of $3:D_8$. The split extension is defined by the homomorphism $D_8\to 2$ with kernel $\la (12), (34)\ra \cong 2^2$. The subgroup isomorphic to $D_8$ is contained in a  nonminimal $S_4$, hence is not minimal. Let $H = \la \bfn_K,(12),(34)\ra \cong 6\times 2 $ so that $3:D_8 \cong H:2$. The subgroup $H$ contains 4 elements of type $D_4$, 4 elements of type $A_2$, 2 elements of type $A_1$ and one element of type $2A_1$. Adding up the traces of elements from $H$ we get the sum equal to  $24$. The nontrivial coset contains $8$ elements of type $A_3+A_1$, 4 elements of type $4A_1$ and one element of type $2A_1$. Adding up the traces we get 0. Thus the group is not minimal. So this case does not reveal any minimal groups.

\smallskip
\emph{Case 3}: $k = 2$.

The image of $G$ in $S_4$ is contained in the stabilizer of the orthogonal vector $\bfn_K$. Thus $G$ is a subgroup of one of the following  groups
$$
G = \begin{cases}
K:2^2&\text{if $K$ is of type I},\\
\text{$K\pt S_3$}&\  \text{if $K$ is of type II},\\
K:D_8&\text{if $K$ is of type III}.
\end{cases}
$$
 Since $K$ of type I or II contains a minimal element of order 3, the subgroups of  $K:2^2$ and $K\pt S_3$ are minimal.  Recall that they all contain $K$.

 Assume $K$ is of type I. Recall that the $S_4$-module $3^3$ is isomorphic to the root lattice $\calR_3(3)$ of $A_3$ modulo 3. A permutation  $\sigma$ of order 2 represented by the transposition $(12)$ decomposes the module into the sum of one-dimensional subspaces with eigenvalues $-1,1,1$. So, if $\sigma$ fixes $\bfn_K$ it acts on $K$ with eigenvalues $-1,1$. Otherwise it acts identically on $K$. The product $(12)(34)$ acts with eigenvalues $(-1,-1,1)$, so if it fixes $\bfn_K$ then it acts as the minus identity on $K$. In our case $\bfn_K = [1,2,0,0]$ and $(12), (12)(34)\in 2^2$ fix $\bfn_K$. Accordingly, $(12)$ acts as $(-1,1)$ giving a subgroup $K:2 \cong 3\times S_3$, $(12)(34)$ gives the subgroup $K:2\cong 3^2:2\not\cong 3^2:2$. Finally, $(34)$ acts identically on $K$ giving the subgroup $3^2\times 2$. So we obtain 3 subgroups of index 2 of $K:2^2$ isomorphic to   $3^2:2, 3^2\times 2, S_3\times 3$. The remaining subgroups are $K\cong 3^2$ and $K:2^2\cong 3^2:2^2$.

 Assume $K$ is of type II. Again we have to find all subgroups $H$ of $G = K \pt S_3$ containing $K$.
 Elements of order 2 are transpositions in $S_3$. They  fix  $\bfn_K$. Arguing as above we see that $K:2 \cong 3\times S_3$.  An element of order 3 in $S_3$ fixes $\bfn_K$. Hence it acts in the orthogonal complement as an element of type $A_2$ in the root lattice of type $A_2$ modulo 3. This defines a unique non-abelian group of order 27 isomorphic to the Heisenberg group $H_3(3)$. The third group is $K\pt S_3 \cong H_3(2):2 \cong 3(3^2:2)$.

Assume $K$ is of type III. This time $K$ is not minimal. Each subgroup of order 2 of type $4A_1$  of $D_8$ defines a subgroup $K:2 \cong 3\times S_3$ of $K:D_8$. It contains 3 elements of order 2, of type $4A_1$, 6 elements of order 6, of type $D_4$, and the elements from $K$. Adding up the traces we get the sum equal to 18. So the subgroup is not minimal.   An element of order 2 of type $2A_1$ defines a subgroup $3^2:2$  not isomorphic to $3\times S_3$.  It contains 9 elements of type $2A_1$, and the elements from $K$. Adding up the traces we get the sum equal to $36$. So the group is not minimal.

Assume $G \cong K:4$. It contains the previous group $3^2:2$ as a subgroup of index 2. It has 9 2-Sylow subgroups of order 4. Thus the nontrivial coset consists of 18 elements of order 4. Each element of order 4 has the trace equal tom 0. This shows that the group is not minimal.

Finally it remains to investigate the group $K:D_8$. It contains the previous group $K:4$ as a subgroup of index 2. The sum of the traces of its elements is equal to 36. The nontrivial coset consists of the union of 4 subsets, each is equal to the set of nontrivial elements in the group of type $K:2 \cong 3\times S_3$. The sum of traces of elements in each subset is equal to $12$. So the total sum is 72 and the group is not minimal.

\smallskip
\emph{Case 3}: $k = 3$.

This gives the groups $3^3:H$, where $H$ is a subgroup of $S_4$ which acts on $3^3$ via the restriction of the homomorphism  $S_4\to \Aut(3^3)$ describing the action of $W(A_3)$ on its root lattice modulo 3. We have the groups
$$3^3:S_4, \ 3^3:D_8,\ 3^3:S_3, \ 3^3:2^2 (2), \ 3^3:3, \ 3^3:4, \ 3^3:2 (2).$$

\medskip\noindent
Type II.

The surface is isomorphic to the Clebsch diagonal cubic surface in $\bbP^4$ given by the equations
$$\sum_{i=0}^4T_i^3 = \sum_{i=0}^4T_i = 0.$$
The group $S_5$ acts by permuting the coordinates. The orbit of the line $x_0 = x_1+x_2 = x_3+x_4 = 0$ consists of 15 lines. It is easy to see that the remaining 12 lines form a double-six. The lines in the double-six are described as follows.

Let $\omega$ be a primitive 5th root of unity. Let $\sigma = (a_1,\ldots,a_5)$ be a permutation of $\{0,1,2,3,4\}$.  Each line
$\ell_\sigma$ is the span by a pair of  points  $(\omega^{a_1},\ldots,\omega^{a_5})$ and
$(\omega^{-a_1},\ldots,\omega^{-a_5})$. This gives $12$ different lines. One  immediately checks that $\ell_\sigma\cap \ell_{\sigma'} \ne \emptyset $ if and only if
$\sigma' = \sigma\circ \tau$ for some odd permutation $\tau$.  Thus the orbit of the alternating subgroup $A_5$ of any line defines a set of 6 skew lines (a sixer) and therefore $A_5$ is not  minimal. Let $\ell_1,\ldots,\ell_6$ be a sixer. It is known that the divisor classes $\ell_i,K_S$ span $\Pic(S)\otimes \bbQ$. This immediately implies that $\Pic(S)^{A_5}$ is spanned (over $\bbQ$) by  $K_S$ and the sum  $\sum \ell_i$. Since $S_5$ does not leave this sum invariant, we see that  $S_5$ is a minimal group.

A maximal subgroup of $S_5$ not contained in $A_5$ is isomorphic to  $S_4$, or $5:4$, or $2\times S_3$. The subgroups isomorphic to $S_4$ are conjugate so we may assume that it consists of permutations of $1,2,3,4$. The group has 6 elements of type $4A_1$ conjugate to $(12)$, 3 elements of type $2A_1$ conjugate to $(12)(34)$, 8 elements of  type $2A_2$ conjugate to $(123)$ and 6 elements of type $A_3+A_1$ conjugate to $(1234)$. The total sum of the traces is equal to 0. So the group is minimal. This gives another, non-geometric proof of the minimality of $S_5$.

Consider a 2-Sylow subgroup $G$ of $S_4$  isomorphic to $D_8$. It consists of 5 elements of order 2, two of type $4A_1$ (from the conjugacy class of $(12)$) and $3$ of type $2A_1$ (from the conjugacy class of
$(12)(34)$). Its cyclic subgroup of order 4 is generated by an element of type $A_3+A_1$. Adding up the traces we see that the subgroup is not minimal. Thus $S_4$ has no minimal proper subgroups.

A subgroup isomorphic to $5:4$ is conjugate to a subgroup generated by two cycles $(01234)$ and
$(0123)$. Computing the traces, we find that the group is not minimal.  The subgroup  isomorphic to $2\times S_3$ is conjugate to a subgroup generated by $(012), (01), (34)$. Its element of order 6 belongs to the conjugacy class $D_4$. So this group is different from the isomorphic group in the previous case. Computing the traces we find that it is not minimal.

\medskip\noindent
Type III.

The surface is a specialization of  a surface of type $IV$. Recall that each nonsingular plane cubic curve is isomorphic to a member of the Hesse pencil
$$T_1^3+T_2^3+T_3^3+6aT_1T_2T_3 = 0.$$
The group of projective automorphisms leaving the pencil invariant is the Hesse group $G_{216}$ of order 216. It is isomorphic to $3^2:\overline{T}$. The stabilizer of a general member of the pencil is isomorphic to a non-abelian extension $3^2:2$. It is generated by
$$g_1 = [t_1,\epsilon_3 t_2, \epsilon_3^2t_3], \ g_2 = [t_2, t_3,t_1], \ g_3 = [t_2,t_1,t_3].$$
The pencil contains 6 members isomorphic to a harmonic cubic. They correspond to the values of the parameters satisfying the equation $8a^6+20a^3-1 = 0$. The stabilizer of a harmonic member is the group $3^2:4$. The additional generator is given by the matrix
$$g_4 = \frac{1}{\sqrt{3}}\begin{pmatrix}1&1&1\\
1&\epsilon_3&\epsilon_3^2\\
1&\epsilon_3^2&\epsilon_3\end{pmatrix}
$$
The pencil also contains 4 anharmonic cubics isomorphic to the Fermat cubic. They correspond to the parameters $a$ satisfying the equation $a^4-a = 0$. The stabilizer of an anharmonic member is isomorphic to $3^2:6$. The additional generator is given by $g_5 = [x_1,\epsilon_3x_2,\epsilon_3x_3]$.

All curves from the Hesse pencil have 9 common  inflection points. If we fix one of them, say $(1,-1,0)$, all nonsingular members acquire a group law. The group of automorphisms generated by $g_1,g_2$ correspond to translations by 3-torsion points. The automorphism $g_3$ is the negation automorphism.
The automorphism $g_4$ is the complex multiplication by $\sqrt{-1}$. The automorphism $g_5$ is the complex multiplication by $e^{2\pi i/3}$.

The Hesse group admits a central  extension $3G_{216} \cong H_3(3):\overline{T}$ realized as the complex reflection group in $\bbC^3$. It acts  linearly on the variables  $T_1,T_2,T_3$ leaving the polynomial  $T_1^3+T_2^3+T_3^3+6aT_1T_2T_3$ unchanged. We denote by $\tilde{g}_i$ the automorphism of the cubic surface obtained from the automorphism $g_i$ of the Hesse pencil by acting identically on the variable $T_0$. The center of the group $3G_{216}$ is generated by $c =
[\tilde{g}_1,\tilde{g}_2] =  [1,\epsilon_3,\epsilon_3,\epsilon_3]$. This is an element of order 3, of minimal type $3A_2$.

Now we have  a complete description of  the automorphism group of a surface of type $IV$. Any minimal subgroup of $H_3(3):2$ can be found among minimal subgroups of surfaces of type I.  However, we have 2 non-conjugate subgroups of type $S_3$ equal to $\la\tilde{g}_1,\tilde{g}_3\ra$ and $\la\tilde{g}_2,\tilde{g}_3\ra$, and two non-conjugate subgroups in $S_3\times 3$ obtained from the previous groups by adding the central element $c$.

Surfaces of type III acquire additional minimal subgroups of the form $A:4$, where $A$ is a subgroup of $H_3(3)$. The element $\tilde{g}_4$ acts by conjugation on the subgroup $H_3(3)$ via
  $(\tilde{g}_1,\tilde{g}_2) \mapsto  (\tilde{g}_2^2, \tilde{g}_1).$
Using $\tilde{g}_4$, we can conjugate the subgroups isomorphic to $S_3, 3\times S_3, 3^2$. Also we get two new   minimal groups  $H_3(3):4$ and $ 12$.

\medskip\noindent
Type V.

  The group $S_4\cong 2^2:S_3$ acts by permuting the coordinates $T_1,T_2,T_3$ and multiplying them by $-1$ leaving the monomial $T_1T_2T_3$ unchanged.  To make the action explicit, we identify $2^2$ with the subspace of $\bbF_2^3$ of vectors whose coordinates add up to 0. The semi-direct product corresponds to the natural action of $S_3$ by permuting the coordinates. Thus
  $g = ((a,b,c),\sigma)\in 2^2:S_3$ acts as the transformation $[t_0,(-1)^at_{\sigma(1)},(-1)^bt_{\sigma(2)},(-1)^ct_{\sigma(3)}]  $.
It is easy to compute the types of elements of $S_4$ in their action on $S$. The group contains 3 elements of type $2A_1$, 6 elements of type $4A_1$, 8 elements of type $2A_2$ and 6 elements of type $A_3+A_1$. Adding up the traces we see that the group is minimal. The subgroup $S_3$ is minimal. No other subgroup is minimal.

\medskip\noindent
Type VII.

The automorphism group of the surface of type VII  is a nonminimal cyclic group of order 8.

\medskip\noindent
Type X.

The automorphism group of the surface of type X consists of the identity, two involutions of type $4A_1$ and one involution of type $2A_1$. Adding up the traces, we get that the group is not minimal.

Let us summarize our result in the following.

\begin{theorem}\label{clasdp3} Let $G$ be a minimal subgroup of automorphisms of a nonsingular cubic surface. Then $G$ is isomorphic to one of the following groups.
\begin{enumerate}
\item $G$ is a subgroup of automorphisms of a surface of type I.
$$S_4 (3), \  S_3\  (2), \ S_3\times 2, \
 S_3\times 3\  (2),\ 3^2:2 \ (2), \ 3^2:2^2,$$
 $$H_3(3):2, \ H_3(3),  \  3^3:2\  (2),\  3^3:2^2\  (2), \ 3^3:3,\ 3^3:S_3,\  3^3:D_8, \ 3^3:S_4,\  3^3:4,$$
$$3^3, \ 3^2 \ (3), \  3^2\times 2,\  9, \ 6\  (2), \ 3.$$
\item $G$ is a subgroup of automorphisms of a surface of type II.
$$S_5, \quad S_4.$$
\item $G$ is a subgroup of automorphisms of a surface of type III.
$$H_3(3):4,\  H_3(3):2, \ H_3(3), \  S_3\times 3, \  S_3,\ 3^2,\ 12, \  6,\ 3.$$

\item $G$ is a subgroup of automorphisms of a surface of type IV.
$$ H_3(3):2, \ H_3(3), \  S_3 \ (2), \ 3\times S_3 \ (2),\ 3^2 \ (2),\  6,\ 3.$$
\item $G$ is a subgroup of automorphisms of a surface of type V.
$$S_4,\ S_3.$$
\item  $G$ is a subgroup of automorphisms of a surface of type VI.
$$6,\quad S_3\times 2, \quad S_3.$$
\item  $G$ is a subgroup of automorphisms of a surface of type VIII.
$$S_3.$$
\end{enumerate}
\end{theorem}

\subsection{Automorphisms of Del Pezzo surfaces of degree 2}
Recall that the center of the Weyl group $W(E_7)$ is generated by an element $w_0$ which acts on the root lattice as the negative of the identity. Its conjugacy class is of type $A_1^7$. The quotient group $W(E_7)' = W(E_7)/\la w_0\ra $ is isomorphic to the simple group $\Sp(6,\bbF_2)$. The extension $2.\Sp(6,\bbF_2)$ splits by the subgroup $W(E_7)^+$  equal to the kernel of the determinant homomorphism $\det:W(E_7)\to \{\pm 1\}$. Thus we have
$$W(E_7) = W(E_7)^+\times \la w_0\ra.$$

Let $H$ be a subgroup of $W(E_7)'$. Denote by $H^+$ its lift to an isomorphic subgroup  of $W^+$. Any other isomorphic lift of $H$ is defined by a nontrivial homomorphism $\alpha:H\to \la w_0\ra \cong 2$. Its elements are the products $h\alpha(h), h\in H^+$. We denote such a lift by $H_\alpha$. Thus all lifts are parametrized by the group $\Hom(H,\la w_0\ra)$and $H^+$ corresponds to the trivial homomorphism.  Note that $wH_\alpha w^{-1} = (w'Hw'{}^{-1})_\alpha$, where $w'$ is the image of $w$ in $W(E_7)'$. In particular,  two lifts of the same group are never conjugate.

Now we apply this to our geometric situation. Let $S$ be a Del Pezzo  surface of degree 2.  Recall that the map $S\to \bbP^2$ defined by $|-K_S|$ is a degree 2 cover. Its branch curve is a nonsingular curve of degree 4. It is convenient to view a Del Pezzo surface of degree 2 as a hypersurface in the weighted projective space $\bbP(1,1,1,2)$ given by an equation of degree 4
\begin{equation}\label{dp2}
T_3^2+F_4(T_0,T_1,T_2) = 0.
\end{equation}
 The  automorphism of the cover  $\gamma = [t_0,t_1,t_2,-t_3]$  defines  the conjugacy class of  a Geiser
involution of $\bbP^2$. For any divisor class $D$ on $S$ we have $D+\gamma_0^*(D)\in |-mK_S|$ for some integer $m$. This easily implies that $\gamma^*$ acts as the minus identity in $\calR_S$. Its image in the Weyl group $W(E_7)$ is the generator $w_0$ of its center. Thus the Geiser involution is the geometric realization of $w_0$.

Let $\rho:\Aut(S) \to W(E_7)$ be the natural injective homomorphism corresponding to a choice of a geometric basis in $\Pic(S)$. Denote by $\Aut(S)^+$ the full preimage of $W(E_7)^+$. Since $W(E_7)^+$ is a normal subgroup,  this definition is independent of a choice of a geometric basis. Under the restriction homomorphism $\Aut(S) \to \Aut(B)$ the group $\Aut(S)^+$ is mapped isomorphically to $\Aut(B)$ and we obtain
$$\Aut(S)^+ \cong \Aut(S)/\la \gamma\ra \cong \Aut(B).$$
>From now on we will identify any subgroup $G$ of $\Aut(B)$ with a subgroup of $\Aut(S)$ which we call
the \emph{even lift} of $G$. Under the homomorphism $\rho:\Aut(S)\to W(E_7)$ all elements of $G$ define \emph{even conjugacy classes}, i.e. the conjugacy classes of elements from $W(E_7)^+$. It is immediate to see that a conjugacy class is even if and only if the sum of the subscripts in its Carter graph  is even. An isomorphic  lift of a subgroup $G$ to a subgroup of $\Aut(S)$ corresponding to some nontrivial homomorphism $G\to \la \gamma\ra$ (or, equivalently to a subgroup of index 2 of $G$) will be called an \emph{odd lift} of $G$.

The odd and even lifts of the same group are never conjugate, two minimal lifts are conjugate in $\Aut(S)$ if and only if the groups are conjugate in $\Aut(B)$. Two odd lifts of $G$ are conjugate if and only if they correspond to conjugate subgroups of index 2 (inside of the normalizer of $G$ in $\Aut(B))$.

The following simple lemma will be heavily used.

\begin{lemma}\label{xxx} Let $G$ be a subgroup of $\Aut(B)$ and $H$ be its subgroup of index 2. Assume $H$ is a minimal subgroup of $\Aut(S)$ (i.e. its even lift is such a subgroup). Then $G$ is  minimal in its even lift and its odd lift corresponding to $H$. Conversely, if $G$ is minimal in both lifts, then $H$ is a minimal subgroup.
\end{lemma}

\begin{proof} Let $\Tr(g)_\pm $ be the trace of $g\in G$ in the representation of $G$ in $\calR_S$ corresponding to the minimal (resp. odd) lift.
Suppose $G$ is minimal in both lifts. Then
$$\sum_{g\in G}\Tr_+(g) = \sum_{g\in H}\Tr_+(g) +\sum_{g\not \in H}\Tr_+(g) =  0,$$
$$\sum_{g\in G}\Tr_-(g) = \sum_{g\in H}\Tr_-(g) +\sum_{g\not \in H}\Tr_-(g) = $$
$$= \sum_{g\in H}\Tr_+(g) -\sum_{g\not \in H}\Tr_+(g) = 0.$$
This implies that $\sum_{g\in H}\Tr_+(g) = 0$, i.e. $H$ is minimal. The converse is obviously true. \
\end{proof}

Since $\gamma$ generates  a minimal subgroup of automorphisms of $S$, any group containing $\gamma$ is minimal. So, we classify first subgroups of $\Aut(B)$ which admit minimal lifts. These will be all minimal subgroups of $\Aut(S)$ which do not contain the Geiser involution $\gamma$. The remaining minimal groups will be of the form $\la \gamma\ra \times \widetilde{G}$, where $\widetilde{G}$ is any lift of a subgroup $G$ of $\Aut(B)$. Obviously, the product does not depend on the parity of the lift.

As in the case of cubic surfaces we first classify cyclic subgroups.

\begin{lemma}\label{finorder} Let $g$ be an automorphism of order $n > 1$ of a nonsingular plane quartic $C = V(F)$. Then one can choose coordinates in such a way that
$g = [t_0,\epsilon_n^at_1,\epsilon_n^bt_2]$ and $F$ is given in the following list.
\begin{itemize}
\item[(i)]  $(n = 2), (a,b) = (0,1),$
$$F = T_2^4+T_2^2L_2(T_0,T_1)+L_4(T_0,T_1). $$
\item[(ii)] $(n = 3), (a,b) = (0,1)$,
$$F = T_2^3L_1(T_0,T_1)+L_4(T_0,T_1). $$
\item[(iii)] $(n = 3), (a,b) = (1,2)$,
$$F = T_0^4+\alpha T_0^2T_1T_2+T_0T_1^3+T_0T_2^3+ \beta T_1^2T_2^2.$$
\item[(iv)] $(n = 4), (a,b) = (0,1),$
$$F = T_2^4+L_4(T_0,T_1). $$
\item[(v)] $(n = 4), (a,b) = (1,2),$
$$F = T_0^4+T_1^4+T_2^4+\alpha T_0^2T_2^2+\beta T_0T_1^2T_2. $$
\item[(vi)]  $(n = 6), (a,b) = (3,2)$,
$$F = T_0^4+T_1^4+\alpha T_0^2T_1^2+T_0T_2^3. $$
\item[(vii)]  $(n = 7), (a,b) = (3,1),$
$$F = T_0^3T_1+T_1^3T_2+T_2^3T_0. $$
\item[(viii)] $(n = 8), (a,b) = (3,7),$
$$F = T_0^4+T_1^3T_2+T_1T_2^3. $$
\item[(ix)] $(n = 9), (a,b) = (3,2),$
$$F = T_0^4+T_0T_1^3+T_2^3T_1. $$
\item[(x)] $(n = 12), (a,b) = (3,4)$,
$$F = T_0^4+T_1^4+T_0T_2^3.$$
\end{itemize}
Here the subscript in polynomial $L_i$ indicates its degree.
\end{lemma}

Also observe  that the diagonal matrix $(t,t,t,t^2)$ acts identically on $S$.

Let $g\in \Aut(B)$ be an element of order $n$ of type $(\ast)$ from the previous Lemma. The following Table identifies the conjugacy class of two lifts $\tilde{g}$ of  $g$   in the Weyl group $W(E_7)$.  If $n$ is even, then $g$ admits two lifts in $\Aut(S)$ of order $n$. If $n$ is odd, then one of the lifts is of order $n$ and another is of order $2n$. We denote by $(\ast)_+$ the conjugacy class of the lift which is represented by an element from $W(E_7)^+$ (of order $n$ if $n$ is odd). The conjugacy class of another lift is denoted by $(\ast)_-$.  The last column of the Table gives the trace of $g$ on $\calR_S$.

\bigskip
\begin{table}[htdp]
\centering
\begin{tabular}{||r |r |r |r ||}\hline
Type &Order&Notation&Trace\\ \hline
$(0)_-$&2&$7A_1$&-7\\
\hline
$(i)_+$&2&$4A_1$&-1\\
$(i)_-$&2&$3A_1$&1\\
\hline
$(ii)_+$&3&$3A_2$&-2\\
$(ii)_-$&6&$E_7(a_4)$&2\\
\hline
$(iii)_+$&3&$2A_2$&1\\
$(iii)_-$&6&$D_6(a_2)+A_1$&-1\\
\hline
$(iv)_+$&4&$D_4(a_1)$&3\\
$(iv)_-$&4&$2A_3+A_1$&-3\\
\hline
$(v)_+$&4&$2A_3$&-1\\
$(v)_-$&4&$D_4(a_1)+A_1$&1\\
\hline
$(vi)_+$&6&$E_6(a_2)$&2\\
$(vi)_-$&6&$A_2+A_5$&-2\\
\hline
$(vii)_+$&7&$A_6$&0\\
$(vii)_-$&14&$E_7(a_1)$&0\\
\hline
$(viii)_+$&8&$D_5+A_1$&-1\\
$(viii)_-$&8&$D_5$&1\\
\hline
$(ix)_+$&9&$E_6(a_1)$&1\\
$(ix)_-$&18&$E_7$&-1\\
\hline
$(x)_+$&12&$E_6$&0\\
$(x)_-$&12&$E_7(a_2)$&0\\
\hline

\end{tabular}
\caption{Conjugacy classes of automorphisms of a Del Pezzo surface of degree 2}\label{tab4}
\end{table}

The following is the list of elements of finite order which generate a minimal cyclic group of automorphisms. To identify the conjugacy class of a minimal lift we use the trace formula \eqref{trace}. If both lifts have the same trace, we distinguish them by computing the traces of their powers.

\begin{enumerate}
\item Order 2 ($A_1^7$) (The Geiser involution) \ $g = [t_0,t_1,t_2,-t_3]$
$$F = T_3^2+F_4(T_0,T_1,T_2). $$
\item Order 4 ($2A_3+A_1$)\ $g = [t_0,t_1,it_2,t_3]$
$$F = T_3^2+T_2^4+L_4(T_0,T_1).$$
\item Order 6 ($E_7(a_4)$)\ $g = [t_0,t_1,\epsilon_3t_2,-t_3]$
$$F = T_3^2+T_2^3L_1(T_0,T_1)+L_4(T_0,T_1).$$
\item Order 6 $(A_5+A_2)$\ $g = [t_0,-t_1,\epsilon_3t_2,-t_3]$
$$ F = T_3^2+T_0^4+T_1^4+T_0T_2^3+aT_0^2T_1^2.$$
\item Order 6 $(D_6(a_2)+A_1)$\ $g = [t_0,\epsilon_3x_1,\epsilon_3^2x_2,-x_3]$
$$F = T_3^2+T_0(T_0^3+T_1^3+T_2^3)+ T_1T_2(\alpha T_0^2+\beta T_1T_2).$$
\item Order 12 $(E_7(a_2))$\ $g = [t_0,\epsilon_4t_1,\epsilon_3t_2,t_3]$
$$F = T_3^2+T_0^4+T_1^4+T_0T_2^3, \ (t_0,t_1,t_2,t_3).$$
\item Order 14 ($E_7(a_1))$\ $g = [t_0,\epsilon_4t_1,\epsilon_3t_2,t_3]$
$$F = T_3^2+T_0^3T_1+T_1^3T_2+T_2^3T_0.$$
\item Order 18 $(E_7)$ \ $g = [t_0,\epsilon_3t_1,\epsilon_9^2t_2,-t_3]$
$$F = T_3^2+T_0^4+T_0T_1^3+T_2^3T_1.$$
\end{enumerate}
 Using the information about cyclic groups of automorphisms of plane quartics, it is not hard to get the classification of possible automorphism groups (see \cite{Topics}). It is given in Table 5.

\begin{table}[h]
\begin{center}
\begin{tabular}{| l|r | r | r | r|r|}
\hline
\footnotesize{Type}&\footnotesize{Order}&\footnotesize{Structure} &\footnotesize{Equation}&\footnotesize{Parameters}\\ \hline
I&336 &\small{$2\times L_2(7)$}&\footnotesize{$T_3^2+T_0^3T_1+T_1^3T_2+T_2^3T_0$}&\\ \hline
II&192&\small{$2\times (4^2:S_3)$}&\footnotesize{$T_3^2+T_0^4+T_1^4+T_2^4$}&\\ \hline
III&96&\small{$2\times 4A_4$}&\footnotesize{$T_3^2+T_2^4+T_0^4+aT_0^2T_1^2+T_1^4$}&\footnotesize{$a^2 = -12$}\\ \hline
IV&48&\small{$2\times S_4$}&\footnotesize{$T_3^2+T_2^4+T_1^4+T_0^4+$}&\footnotesize{$a\ne \frac{-1\pm \sqrt{-7}}{2}$}\\
&&&\footnotesize{$+a(T_0^2T_1^2+T_0^2T_2^2+T_1^2T_2^2)$}&\\\hline
V&32&\small{$2\times AS_{16}$}&\footnotesize{$T_3^2+T_2^4+T_0^4+aT_0^2T_1^2+T_1^4$}&\footnotesize{$a^2\ne 0,-12, 4,36$}\\ \hline
VI&18&\small{$18$}&\footnotesize{$T_3^2+T_0^4+T_0T_1^3+T_1T_2^3$}&\\ \hline
VII&16&\small{$2\times D_8$}&\footnotesize{$T_3^2+T_2^4+T_0^4+T_1^4+aT_0^2T_1^2 +bT_2^2T_0T_1$}&\footnotesize{$a,b\ne 0 $}\\ \hline
VIII&12&$2\times 6$&\footnotesize{$T_3^2+T_2^3T_0+T_0^4+T_1^4+aT_0^2T_1^2$}&\\ \hline
IX&12&\small{$2\times S_3$}&\footnotesize{$T_3^2+T_2^4+aT_2^2T_0T_1+T_0(T_2^3+T_0^3)+ bT_0^2T_1^2$}& \\ \hline
X&$8$&$2^3$&\footnotesize{$T_3^2+T_2^4+T_1^4+T_0^4$}&\footnotesize{distinct $a,b,c\ne 0$}\\
&&&\footnotesize{$+aT_2^2T_0^2+bT_1^2T_2^2+cT_0^2T_1^2$}&\\ \hline
XI&6&$6$&\footnotesize{$T_3^2+T_2^3T_0+L_4(T_0,T_1)$}&\\ \hline
XII&4&$2^2$&\footnotesize{$T_3^2+T_2^4+T_2^2L_2(T_0,T_1)+L_4(T_0,T_1)$}&\\ \hline
XIII&2&$2$&\footnotesize{$T_3^2+F_4(T_0,T_1,T_2)$}&\\ \hline

\end{tabular}
\end{center}
\caption{Groups of automorphisms of Del Pezzo surfaces of degree 2}\label{tab5}
\end{table}

Next we find minimal subgroups of automorphisms of a Del Pezzo surface of degree 2.

As in the previous case it is enough to consider surfaces $S'$ which are not specialized to surfaces $S$ of other types. When this happens we write $S'\to S$. We have
$$IX\to IV \to I,II, $$
$$XII\to X\to VII\to V\to II, III,$$
$$  XI\to VIII \to III.$$
All of this is immediate to see, except the degeneration $VIII\to III$. This is achieved by some linear change of variables transforming the form $x^3y+y^4$ to the form $u^4+2\sqrt{3}iu^2v^2+v^4$. So it suffices to consider surfaces of types
I, II, III, VI.

Before we start the classification we advice the reader to go back to the beginning of the section and recall the concepts of odd and even lifts of subgroups of $\Aut(B)$.

\medskip
Type I.

Since $L_2(7)$ has no subgroups of index 2 (in fact, it  is  a simple group), it admits a unique lift to a subgroup of $\Aut(S)$.  It is known that the group $L_2(7)$ is generated by  elements of order 2, 3 and 7. Consulting Table 4, we find that an element of order 2 must be of type $4A_1$, an element of order 3 must be of types $3A_2$ or $2A_2$, and element of order 7 is  of type $A_6$. To decide the type of a generator $g$ of order 3, we use that it acts as a cyclic permutation of the coordinates in the plane, hence has 3 fixed points $(1,1,1), (1,\eta_3,\eta_3^2), (1,\eta_3^2,\eta_3)$. The last two of them lie on the quartic. This easily implies that $g$ has 4 fixed points on $S$, hence its trace in $\Pic(S)$ is equal to 2. This implies that $g$ is of type $2A_2$.  Comparing the traces with the character table of the group $L_2(7)$ we find that the representation of $L_2(7)$  in $(\calR_S)\otimes\bbC$ is  an irreducible  7-dimensional  representation of $L_2(7)$. Thus the group is minimal.

Assume $G$ is a proper subgroup of $L_2(7)$. It is known that maximal
subgroups of $L_2(7)$ are isomorphic to $S_4$ or $7:3$. There are
two conjugacy classes of subgroups isomorphic to $S_4$ (in the realization
$L_2(7) \cong L_3(2)$ they occur as   the stabilizer subgroups of a point or a line in
$\bbP^2(\bbF_2)$). Since $S_4$ contains a unique subgroup of index 2, each subgroup can be lifted in two ways. Consider the even lift of $S_4$  lying in  $L_2(7)$.
To find the restriction of the 7-dimensional representation $V_7 = (\calR_S)_\bbC$ to $G$ we apply  the Frobenius Reciprocity formula. Let $\chi_k$ denote  a $k$-dimensional irreducible  representation   of $L_2(7)$ and $\bar{\chi}_k$ be its restriction to $S_4$. It is known that the induced character  of the trivial representation of $S_4$ is equal to $\chi_1+\chi_6$ (see \cite{Atlas}). Applying the Frobenius Reciprocity formula we get
 $\la \bar{\chi}_1,\bar{\chi}_7\ra = \la \chi_1+\chi_6,\chi_7\ra  = 0.$
This computation shows that the  even lifts of the two conjugacy classes of $S_4$ in $L_2(7)$ are minimal subgroups. It follows from Lemma \ref{xxx} the  the odd lifts are minimal only if the lift of the subgroup $A_4$ of $S_4$ is minimal.  One checks that the induced character of the trivial representation of $A_4$  is equal to
 $\chi_1+\chi_6+\chi_7$. By the Frobenius Reciprocity formula, the restriction of $V_7$ to $A_4$ contains the trivial summand. Thus $A_4$ is not minimal and we conclude that there are only 2 non-conjugate lifts of $S_4$ to a minimal subgroup of $\Aut(S)$.

 Next consider the subgroup $7:3$. It admits a unique lift. The induced representation of its trivial representation has  the character equal to $\chi_1+\chi_7$. Applying  the Frobenius Reciprocity formula, we see that this group is not minimal.

Let $H$ be any subgroup of $L_2(7)$ which admits a minimal lift. Since $\Aut(S)$ does not contain minimal elements of order 3 or 7, $H$ must  be a subgroup of $S_4$. Since $A_4$ does not admit a minimal lift, $H$ is either a cyclic group or isomorphic to either $2^2$ or $D_8$. The only cyclic group which may admit a minimal lift is a cyclic group of order 4. However, the character table for $L_2(7)$ shows that the value of the character $\chi_7$ at an element of order 4 is equal to $-1$, hence it is of type $2A_3$.  It follows from the Table that this element does not admit minimal lifts.

Suppose $G \cong 2^2$. In the even lift, it contains 3 nontrivial elements of type $4A_1$. Adding up the traces we see that this group is not minimal. In the odd lift, it contains one element of type $4A_1$ and two of type $3A_1$. Again, we see that the group is not minimal.

Assume $G\cong D_8$. The group $S_4$ is the normalizer of $D_8$. This shows that there are  two conjugacy classes of subgroups isomorphic to $D_8$.   The group $G$  admits 2 lifts. In the even lift it contains 2 elements of type $2A_3$ and 5 elements of type $4A_1$.  Adding up the traces, we obtain that the lift is minimal. Since the lift of $4$ is not minimal, the odd lift of $D_8$ is not minimal.

\medskip
Type II.

The group $\Aut(B)$ is generated by  the transformations
$$g_1=  [t_0,it_1,t_2,-t_3],\  \tau =  [t_1,t_0,t_2,t_3],\  \sigma= [t_0,t_2,t_1,t_3]$$
of  types $D_4(a_1), 4A_1, 4A_1$, respectively. Let $g_2 = \sigma g_1\sigma^{-1}= [t_0,t_1,it_2,-t_3].$

We have
$$\tau g_1\tau^{-1} = g_1^{-1}g_2^{-1}, \ \tau g_2 \tau^{-1} = g_2.$$
The elements $g_1, g_2, \gamma$ generate a normal subgroup isomorphic to $4^2$. The quotient group is isomorphic to $S_3$. Its generators  of order 2 can be represented by $\tau$ and
$\sigma$. The elements $g_1^2,g_2^2, \tau,\sigma$ generate a subgroup (not normal) isomorphic to $S_4$. Thus
\begin{equation}
\Aut(B) \cong 4^2:S_3
\end{equation}
 and
\begin{equation}
\Aut(S) \cong 2\times (4^2:S_3).
\end{equation}
Consider the natural homomorphism $f:\Aut(B)\to S_3$ with kernel $4^2$. We will consider different cases corresponding to a possible image of a subgroup $G \subset \Aut(B)$ in $S_3$. For the future use we observe that $\Aut(B)$  does not contain elements of order 6 because its square is an element of type $3A_2$ but all our elements of order 3 are of type $2A_2$. Also it   does not contain $2^3$ (this follow from the presentation of the group). We will also use that $\Aut(B)$ contains 2 conjugacy classes of elements of order 4 of types $D_4(a_1)$ (represented by $ g_1$) and $2A_3$ (represented by $g_1g_2$).

 \smallskip
\emph{Case 3}: $f(G) = \{1\}$.

In this case $G$ is a subgroup of $4^2$. The group itself contains 3 elements of  type $4A_1$,  6 elements  of type $D_4(a_1)$ and 6 elements of type $2A_3$. The sum of the traces is equal to 16. Thus the group is not minimal. So no subgroup is minimal in the even lift. An odd lift corresponding to the homomorphism $4^2\to \la \gamma\ra$ sending an element of type $D_4(a_1)$ to $\gamma$ defines an odd lift. There is only one conjugacy class of subgroups of index 2 in $4^2$. It defines an odd lift of $4^2$. We may assume that the subgroup of index 2 is generated by $g_1,g_2^2$. It admits two odd lifts  corresponding to the subgroups $\la g_1^2,g_2^2\ra$ and $\la g_1g_2^2\ra$. Finally a cyclic subgroup $4$  of type $D_4(a_1)$ admits an odd lift. No other subgroup admits a minimal lift.

 \smallskip
\emph{Case 2}: $|f(G)| = 2$.

Replacing the group by a conjugate group, we may assume that $f(G) = \la \tau\ra$. We have
$$G_1 =f^{-1}(\la \tau\ra) = \la \tau,g_1,g_2\ra \cong 4^2:2 \cong 4D_8,$$
where the center is generated by $g_2$.

Let $H = \la \tau, g_1^2g_2\ra$.  One immediately checks that $H$ is normal in $G_1$ and  isomorphic to $D_8$. We have $G_1 \cong D_8:4$.
The subgroup $H$ consists of 5 elements of type $4A_1$ and 2 elements of type $2A_3$. Adding up the traces we obtain that $H$ is minimal in its even lift. Thus $G_1$ is minimal in its even lift.
The subgroup $G_2 = \la \tau,g_1^2,g_2\ra$ is of order 16. It contains $H$ defining a split extension $D_8:2$ with center generated by $g_2$. It is isomorphic to the group $AS_{16}$  (see Table \ref{order8}) is of index 2 in $G_1$. Since it is minimal, the odd lift of $G_1$ corresponding to this subgroup is minimal.

We check that $\tau g_1$ is of order 8 and the normalizer of the cyclic group $\la \tau g_1\ra$ is generated by this subgroup and $g_1^2$. This gives  us another subgroup $G_3$  of index 2 of $G_1$. It is a group of order 16  isomorphic to $M_{16}$. An element of order 8 is of type $D_5(a_1)+A_1$.
Thus the sum of the traces is equal to $8$. Adding up the traces of elements in the nontrivial coset of
$\la \tau g_1\ra$ we get that the sum is equal to $-8$ (all elements have the trace equal to $-1$). This shows that $G_3$ is minimal. Thus the corresponding odd lift of $G_1$ is minimal.

Let $G$ be a subgroup of index 2 of $G_1$ and $g = \tau g_1^ag_2^b\in G$ be the element of largest possible order in $H$.  We verify that $g^2 = g_2^{2b-a}$. If $g$ is of order 8, we check that it generates either $\la \tau g_1\ra $ considered earlier or its conjugate subgroup. Its normalizer is conjugate to the subgroup $G_3$ considered earlier.
If $g$ is of order 4, then $2b-a \equiv  2 \mod 4$. We list all possible cases and find that  all elements of order 4 are conjugate.  Thus we may assume that $G$ contains
$g = \tau g_1^2$. Now we check that the normalizer of this group is our group  $G_2$.

So, all subgroups of index 2 are accounted for. They are two of them isomorphic to $AS_{16}$ and $M_{16}$. They are all minimal in their even lift, and hence define odd lifts of $G_1$.

Let $G$ be a subgroup of index 4 of $f^{-1}(\tau)$. It follows from above argument that $G$ is conjugate to a subgroup of index 2 of $G_2$ or $G_3$. It could be $D_8, 8$, or $2\times 4$.  The first group is minimal, hence $D_8:2$ admits an odd minimal lift. Other two groups are not minimal. The last group admits an odd minimal lift. Note that it is not conjugate to odd $2\times 4$ from Case 1. Finally a cyclic group of type $D_4(a_1)$ admits an odd minimal lift. It is not conjugate to a group from Case 1.

  \smallskip
\emph{Case 3}: $|f(G)| = 3$.

Without loss of generality we may assume that $f(G) = \la \sigma\tau\ra$. By Lemma \ref{zas}, $G$ is a split extension $H:3$, where $H$ is a subgroup of $\la g_1,g_2\ra$. Let $G_1 = f^{-1}(\la \sigma\tau\ra)$.  It is a split extension $4^2:3$. By Sylow's Theorem, all subgroups of order 3 are conjugate. Thus we may assume that $H$ contains $\sigma\tau$.
The possibilities are $G_1$ or
$G_2 = \la g_1^2,g_2^2,\sigma\tau\ra  \cong 2^2:3 \cong A_4$.  The group $A_4$ has 3 elements of type $4A_1$,  4 elements of type $2A_2$ and 4 elements of type $3A_2$. Adding up the traces we see that the group is minimal. Thus $G_1$ is minimal too. The group $G_1$ does not have  subgroups of index 2, so it does not admit  odd lifts. Other groups in this case are conjugate  to the nonminimal group $\la \sigma\tau \ra$.

  \smallskip
\emph{Case 4}:$f(G) = S_3$.

In this case $G\cap f^{-1}(\la \sigma\tau\ra)$ is a subgroup of index 2 equal  to one of  the two groups considered in the previous case. We get $G = \Aut(B)$, or $G \cong 2^2:S_3 \cong S_4$, or $S_3$. Considering the preimage of $\la \tau\ra$, we find that all groups isomorphic to $S_4$ are conjugate and their Sylow 2-subgroup is $D_8$ from the previous case. Thus both  $\Aut(B)$ and $S_4$ admit two minimal lifts. A group isomorphic to $S_3$ contains 2 elements of type $2A_2$ and it is not minimal in any lift.

\medskip
Type III.

We assume that $a = 2\sqrt{3}i$ in the equation of the surface. The group $\Aut(B)$ is isomorphic to $4A_4$. It is generated  (as always in its even lift) by
$$g_1 = [t_1,t_0,t_2,-t_3], \ g_2 = [it_1,-it_0,t_2,-t_3],$$
$$g_3 = [\epsilon_8^7t_0+\epsilon_8^7t_1,\epsilon_8^5t_0+\epsilon_8t_1, \sqrt{2}\epsilon_{12}t_2,2\epsilon_{6}t_3],
c = [t_0,t_1,it_2,-t_3]$$
The ``complicated''  transformation $g_3$ is of order 3 (see our list of Gr\"undformen for binary polyhedral groups). The generators $g_1,g_2$ are of type $4A_1$, the generator $g_3$ is of type $2A_2$ and the generator $c$ is of type $D_4(a_1)$.

The element $c$ generates the center.  We have
$g_1g_2 = g_2g_1c^2$. This shows that the quotient by $\la c\ra$ is isomorphic to $A_4$ and the subgroup $\la c,g_1,g_2\ra \cong 4D_4$ is a group of order 16 isomorphic to the group $AS_{16}$ (see Table \ref{order8}).

Let $f:\Aut(B)\to A_4$ be the natural surjection with kernel $\la c \ra$. Let $G$ be a subgroup of $\Aut(B)$.

  \smallskip
\emph{Case 1}: $G\subset \Ker(f) \cong  4$.

There are no even minimal subgroups. The whole kernel admits a minimal odd lift.

  \smallskip
\emph{Case 2}: $\# f(G)= 2$.

Without loss of generality we may assume that $f(G) = \la g_1\ra$. The subgroup $f^{-1}(\la g_1\ra )$ is generated by $c,g_1$ and is isomorphic to $4\times 2$. It is not minimal in the even lift and minimal in the unique odd lift. Its subgroup $\la cg_1\ra$ of order 4 is of type $2A_3$ and does not admit minimal lifts.

  \smallskip
\emph{Case 3}:$f(G) = \la g_3\ra$.

We have  $G = f^{-1}(\la g_3\ra) = \la c,g_3\ra = \la cg_3\ra \cong 12$. The element $cg_3$ is of type $E_6$, hence not minimal. Its square is an element of type $E_6(a_2)$, also not minimal.  The subgroup $\la (cg_3)^2\ra \cong 6$ defines an odd minimal lift of $G$. The subgroup $\la (cg_3)^4\ra$ is of order 3. It defines an odd minimal lift of $\la (cg_3)^2\ra$. The group $\la g_3\ra$ admits an odd minimal lift.

  \smallskip
\emph{Case 4}: $f(G) = \la g_1,g_2\ra \cong 2^2$.

 The subgroup $H = f^{-1}(\la g_1,g_2\ra)$ is generated by $c,g_1,g_2$. As we observed earlier, it is isomorphic to the group $AS_{16}$ from Table \ref{order8}. A proper subgroup is conjugate in $\Aut(B)$ to either $\la g_1,g_2\ra \cong D_8$ or $\la cg_1,g_2\ra$. All of the subgroups are isomorphic to $D_8$ with center generated by $c^2$. The cyclic subgroup of order 4 is of type $2A_3$, thus the  subgroups are minimal in the even lift (we have done this computation for surfaces of type II).  Thus the  group $H$ is minimal in the even lift and also minimal in two odd lifts corresponding to its two  subgroups of index 2.

  \smallskip
\emph{Case 5}: $f(G) = A_4$.

It is easy to see that $G$ has non-trivial center (the center of its Sylow 2-subgroup). It is equal to $\la c\ra$ or $\la c^2\ra$. In the first case $G = \Aut(B)$. Since it contains minimal subgroups it is minimal.

A subgroup $G$ of index 2 is isomorphic to $2A_4 \cong D_8:3$. Its  Sylow 2-subgroup  is equal to one of the two subgroups isomorphic to $D_8$ from Case 4. Thus $\Aut(B)$ admits two odd lifts. Since $G$ has no subgroups of index 2, the odd lifts of $G$ do not exist.

\medskip
Type VI.

In this case $\Aut(B) \cong 9$ is  not minimal so does not admit minimal lifts.

\smallskip
To summarize our investigation we give two lists. In the first we list all groups which do not contain the Geiser involution $\gamma$. We indicate by $+$ or $-$ the types of their lifts. Also we indicate  the number of conjugacy classes.

All other minimal groups are of the form $\la \gamma\ra \times G$, where $G$ is one of the lifts of a subgroup of $\Aut(B)$. In the second list  we give only  groups  $2\times G$, where $G$ does not admit a minimal lift. All other groups are of the form $2\times G$, where $G$ is given in the previous table.

\begin{table}[h]
\begin{center}
\begin{tabular}{| l|r | r | r | r||}
\hline
\footnotesize{Type of $S$}&Group&Lift&Conjugacy classes\\ \hline
I&$L_2(7)$&$+$&1\\
&$S_4$&+&2\\
&$D_8$&+&2\\ \hline
II&$4^2:S_3$&+,-&2\\
&$S_4$&+,-&2\\
&$4^2:3 $&+&1\\
&$A_4$&+&1\\
&$4^2:2 \cong D_8:4$&+,-,-&3\\
&$M_{16}$&+&1\\
&$AS_{16}$&+,-&2\\
&$D_8$&+&1\\
&$4^2$&-&1\\
&$2\times 4$&-&2\\
&$4$&-&2\\ \hline
III&$4A_4$&+,-&2\\
&$2A_4\cong D_8:3$&+&1\\
&$AS_{16}$&+,-&2\\
&$D_8$&+&1\\
&$12$&-&1\\
&$6$&-&1\\
&$2\times 4$&-&1\\
&$4$&-&1\\ \hline
IV&$S_4$&+&1\\
&$D_8$&+&1\\ \hline
V&$AS_{16}$&+,-&2\\
&$D_8$&+&1\\
&$2\times 4$&-&2\\
&$4$&-&1\\ \hline
VII&$D_8$&+&1\\\hline
VIII&$6$&-&1\\\hline

\end{tabular}
\end{center}
\caption{Minimal groups of automorphisms  not containing $\gamma$}\label{tab6}
\end{table}

\begin{theorem}\label{clasdp2} Let $G$ be a minimal group of automorphisms of a Del Pezzo surface of degree 2. Then $G$ is either equal to a minimal lift of a subgroup from Table \ref{tab6} or equal to $\gamma \times G'$, where $G'$ is either from the table or one of the following groups of automorphisms of the branch quartic curve $B$
\begin{enumerate}
\item Type I:
$7: 3,\ A_4,\ S_3, \ 7, \ 4, \ 3, \ 2.$
\item Type II:
$2^2,\ S_3,\ 8, \ 4,\  3,\ 2.$
\item Type III:
$2^2,\ 4,  \ 2.$
\item Type IV:
$S_3, \ 2^2,\ 3, \ 2.$
\item Type V:
$2^2, \ 2.$
\item Type VI:
$9,\ 3.$
\item Type VII:
$2^2,\ 4,\ 2$
\item Type VIII:
$3$.
\item Type IX:
$S_3,\ 3, \ 2.$
\item Type X:
$2^2,\ 2.$
\item Type XI:
$3.$
\item Type XII:
$\{1\}.$

\end{enumerate}
\end{theorem}

\subsection{Automorphisms of Del Pezzo surfaces of degree 1}\label{clasdp1}
Let $S$ be a Del Pezzo surface of degree 1. The linear
system $|-2K_S|$ defines a finite map of degree 2 onto a quadric
cone $Q$  in $\bbP^3$. Its branch locus is a nonsingular curve $B$ of genus 4
cut out by a cubic surface. Recall that a singular quadric is
isomorphic to the weighted projective space $\bbP(1,1,2)$. A curve
of genus 4 of degree $6$ cut out in $Q$ by a cubic surface is given
by equation $F(T_0,T_1,T_2)$ of degree 6. After change of
coordinates it can be given by an equation
$T_2^3+F_4(T_0,T_1)T_2+F_6(T_0,T_1) = 0$, where $F_4$ and $F_6$ are
binary forms of degree $4$ and $6$.   The double cover of $Q$
branched along such curve is isomorphic to a hypersurface of degree
6 in $\bbP(1,1,2,3)$
\begin{equation}
T_3^2+T_2^3+F_4(T_0,T_1)T_2+F_6(T_0,T_1) = 0.
\end{equation}

The vertex of $Q$ has coordinates $(0,0,1)$ and its preimage in the cover consist of one point $(0,0,1,a)$, where $a^2+1 = 0$ (note that $(0,0,1,a)$ and $(0,0,1,-a)$ represent the same point on $\bbP(1,1,2,3))$. This is the base-point of $|-K_S|$. The members of $|-K_S|$ are isomorphic to genus 1 curves with equations $y^2+x^3+F_4(t_0,t_1)x+F_6(t_0,t_1) = 0$.  The locus of zeros of $\Delta = F_4^3+27F_6^2$ is the set of points in $\bbP^1$ such that the corresponding genus 1 curve is singular. It consists of $a$ simple roots and $b$ double roots. The zeros of $F_4$ are either common zeros  with $F_6$ and $\Delta$, or represent nonsingular elliptic curves isomorphic to an anharmonic plane cubic curve.  The zeros of $F_6$ are either common zeros  with $F_4$ and $\Delta$, or represent nonsingular elliptic curves isomorphic to a harmonic plane cubic curve.

Observe that no common root of $F_4$ and $F_6$ is a multiple root of $F_6$ since otherwise the surface is singular.

Since the ramification curve of the cover $S\to Q$ (identified with the branch curve $B$) is obviously invariant with respect to $\Aut(S)$ we have a natural surjective homomorphism
\begin{equation}\label{rest3}
\Aut(S) \to \Aut(B).
\end{equation}
Its kernel is generated by the deck involution $\beta$ which we call  the \emph{Bertini involution}. It defines the Bertini involution in $\Cr(2)$. The Bertini involution is the analog of the Geiser involution for Del Pezzo surfaces of degree 2. The same argument as above shows that $\beta$ acts in $\calR_S$ as the minus of the identity map. Under the homomorphism $\Aut(S) \to W(E_8)$ defined by a choice of a geometric basis, the image of $\beta$ is the  elements $w_0$ generating the center of $W(E_8)$. This time $w_0$ is an even element, i.e. belongs to $ W(E_8)^+$. The quotient group $W(E_8)^+/\la w_0\ra$ is isomorphic to the simple group $\Or(8,\bbF_2)^+$.

Since $Q$ is a unique quadric cone containing $B$, the group $\Aut(B)$  is a subgroup of $\Aut(Q)$.  Consider the natural homomorphism
$$r:\Aut(B)\to  \Aut(\bbP^1).$$

Let $G$ be a subgroup of $\Aut(B)$ and $P$ be its image in $\Aut(\bbP^1)$. We assume that elements from $G$ act on the variables $T_0,T_1$ by linear transformations with determinant 1. The polynomials $F_4$ and $F_6$ are the relative  invariants of  the binary group $\bar{P}$.  They are polynomials in Gr\"undformen which were listed in section \ref{5.5}. Let $\chi_4,\chi_6$ be the corresponding characters of $\bar{P}$ defined by the binary forms $F_4, F_6$. Let $\chi_2, \chi_3$ be the characters of $G$ defined by the action on the variables $T_2,T_3$.  Assume that $F_4\ne 0$. Then
$$\chi_4\chi_2 = \chi_6 = \chi_3^3 =\chi_3^2.$$
If $g\in G\cap \Ker(r)\setminus \{1\}$, then $g$ acts on the variables $T_0,T_1$ by either the identity or  the minus identity. Thus $\chi_4(g) = \chi_6(g) = 1$ and we must have $\chi_2(g)  = \chi_3(g)^2 = 1$. This shows that $g = [t_0,t_1,t_3,-t_3] = [-t_0,-t_1,t_2,-t_3] = \beta$.

If $F_4 = 0$, then we must have only $\chi_2(g)^3 = \chi_3(g)^2 = 1$. Since $[-t_0,-t_1,t_2,-t_3]$ is the identity transformation, we may assume that $\chi_3(g) = 1$ and represent $g$ by
$g =  [t_0,t_1,\epsilon_3 t_3,\pm t_3]$. Thus $G\cap \Ker(r)  =\la \beta,\alpha\ra \cong 6$.

Conversely, start with a polyhedral group $P$ such that its lift to a binary polyhedral group $\bar{P}$ acts on the variables $T_0,T_1$ leaving $V(F_4)$ and $V(F_6)$ invariant. Let $\chi_4, \chi_6$ be the corresponding characters. Assume that  there exist  character $\chi_2,\chi_3:\bar{P} \to \bbC^*$ such that
\begin{equation}\label{char3}
\chi_0^2 = \chi_4\chi_1 = \chi_6 = \chi_1^3.
\end{equation}
Then $g = [at_0+bt_1,ct_0+t_1]\in \bar{P}$ acts on $S$ by
$[at_0+bt_1,ct_0+t_1,\chi_2(g)t_2,\chi_3(t_3)].$ This transformation is the identity in $\Aut(S)$ if and only if $g = [-t_0,-t_1]$ and $\chi_2(-1) = 1, \chi_3(-1) = -1$. This shows that $\bar{P}$ can be identified  with a subgroup of $\Aut(S)$ with $-I_2 = \beta$   if and only if $\chi_3(-1) = -1$. If $\chi_3(-1) = 1$, then $P$ can be identified with a subgroup  of $\Aut(S)$ not containing $\beta$. In the latter case,
$$r^{-1}(P) = \begin{cases} P\times \la \beta\ra&\text{if}\   F_4 \ne 0\\
P\times \la \beta,\alpha\ra& \text{otherwise}.\end{cases}
 $$
 In particular, if $F_4 = 0$, there are three subgroups of $\Aut(S)$ which are mapped surjectively to $P$.

 In the former case
$$r^{-1}(P) = \begin{cases} \bar{P}&\text{if}\   F_4 \ne 0,\\
 \bar{P}\times \la \alpha\ra& \text{otherwise}.\end{cases}
 $$

Of course it could happen that  neither $P$ nor $ \bar{P}$ lifts to a subgroup of $\Aut(S)$. In this case
$r^{-1}(P) \cong 2P \not\cong \bar{P}$ or $r^{-1}(P) \cong 3\times 2P$ (if $F_4 = 0$).

In the following list we give a nontrivial subgroup $P$ of $\Aut(\bbP^1)$ as a group of automorphisms of $B$ and a smallest lift $\widetilde{P}$ of $P$ to a subgroup of $r^{-1}(P)$.  If $F_4\ne 0$, then we will see that $\tilde{P} = r^{-1}(P)$ or $\tilde{P} \cong P$. In the latter case  $r^{-1}(P) \cong 2\times P$. If $F_4 = 0$, and $\tilde{P} \cong P$, then   $r^{-1}(P) \cong 6\times \tilde{P}$. Otherwise $r^{-1}(P) \cong 3\times \widetilde{P}$.

Also we give generators of $\tilde{P}$ to $\Aut(S)$ as a group acting on $t_0,t_1$ with determinant 1 and the Bertini involution as an element of the lift.

\begin{enumerate}\label{enu}
\item {\it Cyclic groups $P$}

\begin{itemize}
\item[(i)] $P = \{1\}, F_4 = 0, r^{-1}(P) = \la \beta,\alpha\ra \cong 6$.
\item[(ii)] $P \cong  2,  g =  [it_0,-it_1,-t_2,it_3],  $
$$F_4 = F_2(T_0^2,T_1^2)\ne 0, \quad F_6 = F_3(T_0^2,T_1^2);$$
\item[(iii)] $P \cong  2, \widetilde{P} \cong 4,   g =  [it_0,-it_1,t_2,t_3], \beta = g^2, $
$$F_4 = a(T_0^4+T_1^4)+bT_0^2T_1^2, \quad F_6 = T_0T_1F_2(T_0^2,T_1^2);$$
\item[(iv)] $P \cong  2,  \widetilde{P} = 4, g =  [-t_1,t_0,t_2,t_3], \beta = g^2, $
{\small $$ F_4\  \text{as in (iii)}, \ F_6 = (T_0^2+T_1^2)(a(T_0^4+T_1^4)+T_0T_1(bT_0T_1+c(T_0^2-T_1^2));$$}
\item[(v)] $P \cong  2,   g =  [-t_1,t_0,-t_2,it_3], $
$$ F_4\ \text{as in (iii)} , \  F_6 = a(T_0^6-T_1^6)+bT_0T_1(T_0^4+T_1^4);$$
\item[(vi)] $P \cong  3,    g =  [\epsilon_3t_0,\epsilon_3^2t_1,\epsilon_3^2t_2,t_3], $
$$F_4 = T_0(aT_0^3 +bT_1^3), \quad F_6= F_2(T_0^3,T_1^3);$$
\item[(vii)] $P \cong 3, g = [\epsilon_3t_0,\epsilon_3^2t_1,t_2,t_3], $
$$F_4 = aT_0^2T_1^2, \quad F_6 =  F_2(T_0^3,T_1^3);$$
\item[(viii)] $P \cong  4, g =  [\epsilon_8t_0,\epsilon_8^{-1}t_1,it_2,\epsilon_8^3t_3], $
$$F_4 = aT_0^4+bT_1^4, \quad F_6 = T_0^2(cT_0^4+dT_1^4);$$
\item[(ix)] $P \cong 4, \widetilde{P} \cong 8,  g =    [\epsilon_8t_0,\epsilon_8^{-1}t_1,-t_2,t_3], \beta = g^4$,
$$F_4 = aT_0^2T_1^2, \quad F_6 = T_0T_1(T_0^4+T_1^4);$$
\item[(x)] $P \cong  5,  g =  [\epsilon_{10}t_0,\epsilon_{10}^{-1}t_1,\epsilon_5t_2,\epsilon_{10}^3t_3]$,
$$F_4 = aT_0^4, \quad F_6 = T_0(T_0^5+T_1^5);$$
\item[(xi)] $P \cong 6,   g =  [\epsilon_{12}t_0,\epsilon_{12}^{-1}t_1,\epsilon_6t_2,it_3], $
$$F_4 = aT_0^4, \quad F_6 = bT_0^6+T_1^6, b\ne 0;$$
\item[(xii)] $P \cong 6,  g =  [\epsilon_{12}t_0,\epsilon_{12}^{-1}t_1,-t_2,it_3], $
$$F_4 = aT_0^2T_1^2, \quad F_6 = T_0^6+T_1^6;$$
\item[(xiii)] $P \cong 10, \widetilde{P} \cong 20,  g =  [\epsilon_{20}t_0,\epsilon_{20}^{-1}t_1,\epsilon_{10}^8t_2,\epsilon_{10}^{-1}t_3], g^{10} = \beta$,
$$F_4 = T_0^4, \quad F_6 = T_0T_1^5;$$
\item[(xiv)] $P \cong 12,  g =  [\epsilon_{24}t_0,\epsilon_{24}^{-1}t_1,\epsilon_{12}t_2,\epsilon_{24}t_3], $
$$F_4 = T_0^4, \quad F_6 = T_1^6.$$

\end{itemize}
\item {\it Dihedral groups}
\begin{itemize}
\item[(i)] $P \cong 2^2, \widetilde{P} \cong D_8, g_1= [it_1,it_0,-t_2,it_3], \ g_2= [-t_1,t_0,-t_2,it_3], \beta = (g_1g_2)^2, g_1^2=g_2^2 = 1,$
$$ F_4 = a(T_0^4+T_1^4)+bT_0^2T_1^2,  \quad F_6 = T_0T_1[c(T_0^4+T_1^4)+dT_0^2T_1^2];$$
\item[(ii)] $P \cong  2^2, \widetilde{P} \cong Q_8, g_1= [it_1,it_0,t_2,t_3], \ g_2= [-t_1,t_0,t_2,t_3], \beta = g_1^2=g_2^2$,
$$ F_4 = a(T_0^4+T_1^4)+bT_0^2T_1^2,  \quad F_6 = T_0T_1(T_0^4-T_1^4);$$
\item[(iii)] $P \cong  D_6,  g_1= [\epsilon_6t_0,\epsilon_6^{-1}t_1,t_2,-t_3], g_2=
[it_1,it_0,-t_2,it_3], $
$$F_4 = aT_0^2T_1^2, \quad F_6 = T_0^6+T_1^6+bT_0^3T_1^3; $$
\item[(iv)] $P \cong D_8, \widetilde{P} \cong  D_{16},\  g_1 = [\epsilon_8t_0,\epsilon_8^{-1}t_1,-t_2,it_3], 
g_2= [-t_1,t_0,-t_2,it_3],\\
g_1^4 = \beta, g_2^2 = 1,$
$$F_4 = aT_0^2T_1^2, \quad F_6 = T_0T_1(T_0^4+T_1^4);$$
\item[(v)] $P \cong D_{12}, \widetilde{P}  \cong 2D_{12} \cong (2\times 6)^\bullet 2, g_1= [\epsilon_{12}t_0,\epsilon_{12}^{-1}t_1,-t_2,it_3], g_2= [-t_1,t_0,t_2,t_3],  g_1^6 = 1,  g_2^2 = \beta,$
$$F_4 = aT_0^2T_1^2, \quad F_6 = T_0^6+T_1^6.$$

\end{itemize}
\item {\it Other groups}
\begin{itemize}
\item[(i)] $P \cong A_4, \widetilde{P}  \cong \overline{T}, g_1 =  [\epsilon_8^7t_0+\epsilon_8^7t_1,\epsilon_8^5t_0+\epsilon_8t_1,\sqrt{2}\epsilon_3t_2,2t_3],\
g_2= [it_0,-it_1,t_2,t_3], \  g_1^3 = g_2^2 = \beta,$
$$F_4 = T_0^4+2\sqrt{-3}T_0^2T_1^2+T_2^4, \quad F_6 = T_0T_1(T_0^4-T_1^4),$$
\item[(ii)] $P \cong O, \widetilde{P} \cong \overline{T}:2,  g_1 =  [\epsilon_8^7t_0+\epsilon_8^7t_1,\epsilon_8^5t_0+\epsilon_8t_1,\sqrt{2}\epsilon_3t_2,2t_3],\  g_2= [\epsilon_8t_0,\epsilon_8^{-1}t_1,-t_2,it_3], \  g_3 = [-\epsilon_8t_1,\epsilon_8^7t_0,-t_2,it_3], g_1^3 = g_2^4 = \beta, g_3^2 = 1, r^{-1}(P) = 3\times \bar{O},$
$$F_4 = 0, \quad F_6 = T_0T_1(T_0^4-T_1^4),$$

\end{itemize}
\end{enumerate}
Table 6 gives the list of the full automorphism groups of Del Pezzo surfaces of degree 1.

\begin{table}[h]
\begin{center}
\begin{tabular}{| l|r | r | r | r|r|r|}
\hline
\footnotesize{Type}&\footnotesize{Order}&\footnotesize{Structure} &\footnotesize{$F_4$}&\footnotesize{$F_6$}&\footnotesize{Parameters}\\ \hline
I&144&\small{$3\times (\overline{T}:2)$}&\footnotesize{$0$}&\footnotesize{$T_0T_1(T_0^4-T_1^4)$}&\\ \hline
II&72&\small{$3\times 2D_{12}$}&\footnotesize{$0$}&\footnotesize{$T_0^6+T_1^6$}&\\ \hline
III&36&\small{$6\times D_{6}$}&\footnotesize{$0$}&\footnotesize{$T_0^6+aT_0^3T_1^3+T_1^6$}&\footnotesize{$a\ne 0$}\\ \hline
IV&30&\small{$30$}&\footnotesize{$0$}&\footnotesize{$T_0(T_0^5+T_1^5)$}&\\ \hline
V&24&\small{$\overline{T}$}&\footnotesize{$a(T_0^4+\alpha T_0^2T_1^2+T_1^4)$}&\footnotesize{$T_0T_1(T_0^4-T_1^4)$}&\footnotesize{$\alpha = 2\sqrt{-3}$}\\ \hline
VI&24&\small{$2D_{12}$}&\footnotesize{$aT_0^2T_1^2$}&\footnotesize{$T_0^6+T_1^6$}&\footnotesize{$a\ne 0$}\\ \hline
VII&24&\small{$2\times 12$}&\footnotesize{$T_0^4$}&\footnotesize{$T_1^6$}&\\ \hline
VIII&20&\small{$20$}&\footnotesize{$T_0^4$}&\footnotesize{$T_0T_1^5$}&\\ \hline
IX&16&\small{$D_{16}$}&\footnotesize{$aT_0^2T_1^2$}&\footnotesize{$T_0T_1(T_0^4+T_1^4)$}&\footnotesize{$a\ne 0$}\\ \hline
X&12&\small{$D_{12}$}&\footnotesize{$T_0^2T_1^2$}&\footnotesize{$T_0^6+aT_0^3T_1^3+T_1^6$}&\footnotesize{$a\ne 0 $}\\ \hline
XI&12&\small{$2\times 6$}&\footnotesize{$0$}&\footnotesize{$G_3(T_0^2,T_1^2)$}&\\ \hline
XII&12&\small{$2\times 6$}&\footnotesize{$T_0^4$}&\footnotesize{$aT_0^6+T_1^6$}&\footnotesize{$a\ne 0$}\\ \hline
XIII&10&\small{$10$}&\footnotesize{$T_0^4$}&\footnotesize{$T_0(aT_0^5+T_1^5)$}&$a\ne 0$\\ \hline
XIV&8&\small{$Q_8$}&\footnotesize{$T_0^4+T_1^4+aT_0^2T_1^2$}&\footnotesize{$bT_0T_1(T_0^4-T_1^4)$}&\footnotesize{$a\ne 2\sqrt{-3}$}\\ \hline
XV&8&\small{$2\times 4$}&\footnotesize{$aT_0^4+T_1^4$}&\footnotesize{$T_0^2(bT_0^4+cT_1^4)$}& \\ \hline
XVI&8&\small{$D_8$}&\footnotesize{$T_0^4+T_1^4+aT_0^2T_1^2$}&\footnotesize{$T_0T_1(b(T_0^4+T_1^4)+cT_0^2T_1^2)$}&\footnotesize{$b\ne 0$}\\ \hline
XVII&6&\small{$6$}&\footnotesize{$0$}&\footnotesize{$F_6(T_0,T_1)$}&\\ \hline
XVIII&6&$6$&\footnotesize{$T_0(aT_0^3+bT_1^3)$}&\footnotesize{$cT_0^6+dT_0^3T_1^3+T_1^6$}&\\ \hline
XIX&4&\small{$4$}&\footnotesize{$G_2(T_0^2,T_1^2$}&
\footnotesize{$T_0T_1F_2(T_0^2,T_1^2)$}&\\ \hline
XX&4&\small{$2^2$}&\footnotesize{$G_2(T_0^2,T_1^2$}&
\footnotesize{$G_3(T_0^2,T_1^2)$}&\\ \hline
XXI&2&\small{$2$}&\footnotesize{$F_4(T_0,T_1)$}&\footnotesize{$F_6(T_0,T_1)$}&\\ \hline

\end{tabular}
\end{center}
\caption{Groups of automorphisms of Del Pezzo surfaces of degree 1}\label{tab7}
\end{table}

\bigskip
The following is the list of cyclic minimal groups $\la g\ra $ of automorphisms of Del Pezzo surfaces $V(F)$ of degree 1.
\begin{enumerate}
\item Order 2
\begin{itemize}
\item $A_1^8$ (the Bertini involution)
  $g = [t_0,t_1,t_2,-t_3]$
$$F = T_3^2+T_2^3+F_4(T_0,T_1)T_2+F_6(T_0,T_1),$$
\end{itemize}
\item Order 3
\begin{itemize}
\item $4A_2$ \ $g = [t_0,t_1,\epsilon_3t_2,t_3]$
$$F =T_3^2+T_2^3+F_6(T_0,T_1),$$
\end{itemize}
\item Order 4
\begin{itemize}
\item $2D_4(a_1)$\  $g = [t_0,-t_1,-t_2,\pm it_3]$
$$ F =T_3^2+T_2^3+(aT_0^4+bT_0^2T_1^2+cT_1^4)T_2+T_0T_1(dT_0^4+eT_1^4),$$
\end{itemize}
\item Order 5
\begin{itemize}
\item $2A_4$\  $g = [t_0,\epsilon_5t_1,t_2,t_3]$
$$F =T_3^2+T_2^3+aT_0^4T_2+T_0(bT_0^5+T_1^5),$$
\end{itemize}
\item Order 6
\begin{itemize}
\item $E_6(a_2)+A_2$\ $g = [t_0,-t_1,\epsilon_3t_2,t_3]$
$$F = T_3^2+T_2^3+G_3(T_0^2,T_1^2),$$
\item $E_7(a_4)+A_1$\  $g = [t_0,\epsilon_3t_1,t_2,-t_3]$
$$F =T_3^2+T_2^3+(T_0^4+aT_0T_1^3)T_2+bT_0^6+cT_0^3T_1^3+dT_1^6, $$
\item $2D_4$\  $g = [\epsilon_6t_0,\epsilon_6^{-1}t_1,t_2,t_3]$
$$F =T_3^2+T_2^3+aT_0^2T_1^2T_2+bT_0^6+cT_0^3T_1^3+eT_1^6, $$
\item $E_8(a_8)$\  $g = [t_0,t_1,\epsilon_3t_2,-t_3]$
$$F =T_3^2+T_2^3+F_6(T_0,T_1), $$
\item $A_5+A_2+A_1$\ $g =  [t_0,\epsilon_6t_1,t_2,t_3]$
$$F =T_3^2+T_2^3+aT_0^4T_2+T_0^6+bT_1^6, $$
\end{itemize}
\item Order 8
\begin{itemize}
\item $D_8(a_3)$\ $g= [it_0,t_1,-it_2,\pm \epsilon_8t_3]$
$$F = T_3^2+T_2^3+aT_0^2T_1^2T_2+T_0T_1(T_0^4+T_1^4), $$
\end{itemize}
\item Order 10
\begin{itemize}
\item $E_8(a_6)$\ $g = [t_0,\epsilon_5t_1,t_2,-t_3]$
$$F = T_3^2+T_2^3+aT_0^4T_2+T_0(bT_0^5+T_1^5), $$
\end{itemize}
\item Order 12
\begin{itemize}
\item $E_8(a_3)$\ $g = [-t_0,t_1,\epsilon_6t_2,it_3]$
$$F = T_3^2+T_2^3+T_0T_1(T_0^4+aT_0^2T_1^2+T_1^4), $$
\end{itemize}
\item Order 15
\begin{itemize}
\item $E_8(a_5)$\ $g = [t_0,\epsilon_5t_1,\epsilon_3t_2,t_3]$
$$F = T_3^2+T_2^3+T_0(T_0^5+T_1^5), $$
\end{itemize}
\item Order 20
\begin{itemize}
\item $E_8(a_2)$\ $g = [t_0,\epsilon_{10}t_1,-t_2,it_3]$
$$F = T_3^2+T_2^3+aT_0^4T_2+T_0T_1^5, $$
\end{itemize}
\item Order 24
\begin{itemize}
\item $E_8(a_1)$\ $g = [it_0,t_1,\epsilon_{12}t_2,\epsilon_8t_3]$
$$F = T_3^2+T_2^3+T_0T_1(T_0^4+T_1^4), $$
\end{itemize}
\item Order 30
\begin{itemize}
\item $E_8$\ $g = [t_0,\epsilon_5t_1,\epsilon_{3}t_2,-t_3]$
$$F = T_3^2+T_2^3+T_0(T_0^5+T_1^5).$$
\end{itemize}
\end{enumerate}

To list all minimal subgroups of $\Aut(S)$ is very easy. We know that any subgroup in $\Ker(r)$ contains contains  one of the elements $\alpha,\beta, \alpha\beta$ which are all minimal of types $8A_1,4A_2,E_8(a_8)$.  So, a subgroup is not minimal only if its image $P$ in $\Aut(B)$ can be lifted isomorphically to $\Aut(S)$.

We will use the following lemma.

\begin{lemma}\label{use} Let $P\subset \Aut(\bbP^1)$ and $G \subset \Aut(S)$  be contained in $r^{-1}(P)$. Then $G$ is a minimal group unless $G = \tilde{P} \cong P$ and $G$ is a nonminimal cyclic group  or nonminimal dihedral  group $D_6$.
\end{lemma}

\begin{proof} It follows from above classification of possible subgroups of $\Aut(B)$ and its lifts to $\Aut(S)$ that any non-isomorphic lift contains $\beta$ or $\alpha$, or $\beta\alpha$ which generate minimal cyclic groups. If the lift is isomorphic to $P$ then $P$ is either a cyclic group or $P\cong D_6$. The group $D_6$ contains 3 elements of type $4A_1$ and 2 elements of type $2A_2$. Adding up the traces we see that the group is not minimal.
\end{proof}

Let us classify minimal groups of automorphisms of a Del Pezzo surface of degree 1. As in the previous cases, to find a structure of such groups is enough to consider the types of surfaces which are not specialized to surfaces of other types. The following notation Type A$ \to $Type B indicates that a surface of type A specializes to a surface of type B.
$$V, IX, XIV, XVI, XVII, XIX, XXI     \to I, $$
$$III, VI, X, XI, XII, XVI, XVII, XVIII, XX, XXI \to II$$
$$XIII ,XXI\to IV, \quad XIII ,XXI\to VIII,$$
$$XII, XX, XXI \to VII, \ XX,XI \to XV.$$
It remains to consider surfaces of types
$$I,\  II,  \ IV,  \ VII, \ VIII, \ XV.$$

\medskip\noindent
Type I. $P \cong  S_4$.

Possible conjugacy classes of subgroups $H$ are $\{1\}, 2, 2, 3, 2^2,4, D_8, D_6, A_4, S_4$. Groups of order 2 have two conjugacy classes in $P$ represented by
$[it_0,-it_1]$ and $[-t_1,t_0]$. The groups are realized in cases (iii) and (iv). None of them lifts isomorphically. An cyclic group of order 3 is generated by a nonminimal element realized in case (vii). Its isomorphic minimal lift is not minimal. A cyclic group of order 4 does not admit an isomorphic lift. The dihedral subgroup $2^2$ is of type (ii). This information, together with Lemma \ref{use} allows us to classify all minimal subgroups.

\smallskip
\begin{itemize}
\item $P = \{1\}$:\  $\la \beta\alpha \ra \cong 6, \la \alpha \ra \cong 3, \la \alpha \ra \cong 2;$
\item $P = 2$: \  $4$,   $12$;
\item $P = 2$: \  $4$,  $12$;
\item $P = 3$: \   $3^2$, $3\times 6;$
\item $P = 2^2$: \   $Q_8$, $Q_8\times 3$;
\item $P = 2^2$: \   $D_8$, $D_8\times 3$;
\item $P = 4$: \  $8$, $8\times 3$;
\item $P = D_8$: \  $D_{16}$, $D_8\times 3$;
\item $P = D_6$: \  $D_{6}\times 2$, $D_{6}\times 3, D_6\times 6$;
\item $P = A_4$: \  $\overline{T}$, $\overline{T}\times 3$;
\item $P = S_4$: \  $\overline{T}:2$, $3\times (\overline{T}:2)$.
\end{itemize}

Surfaces specializing to a surface of type I have the following minimal subgroups.

V: $4$, $6$, $Q_8$, $\overline{T}.$

IX:  $4(2)$,  $8$,  $D_{16}$.

XIV: $4$,  $Q_8$.

XVI: $D_8$.

XVII:$2$, $3$, $6$.

XIX: $2$, $4$.

XXI:\ $2$.

\smallskip
Type II: $P = D_{12}$.

Possible subgroups are $\{1\}, 2, 2, 3, 2^2, 6, D_6, D_{12}$. Cyclic subgroups of order 2, 3 and 6 admit isomorphic nonminimal lifts. All these groups are not minimal. There are  two conjugacy classes of subgroups of order 2 in $P$ represented by $[it_0,-it_1]$ and $[-t_1,t_0]$.  One subgroups lifts isomorphically, other one does not.  The cyclic group of order 6 admits an isomorphic lift and not minimal. The dihedral group $D_6$ admits a nonminimal isomorphic lift.
\begin{itemize}
\item $P = \{1\}$: $\la \beta\alpha \ra \cong 6, \la \alpha \ra \cong 3, \la \alpha \ra \cong 2;$
\item $P = 2$: \  $4$,  $12$;
\item $P = 2$:\   $2^2$,\  $2^2\times 3$, $6$;
\item $P = 3$: \   $3^2$,\  $3^2\times 2;$
\item $P = 2^2$: \  $Q_8$,\  $Q_8\times 3$;
\item $P = 6$: \  $2\times 6$,
\item $P = D_6$: \  $2\times D_{6}$,\  $D_{6}\times 3, \ D_6\times 6$;
\item $P = D_{12}$:\   $2D_{12} $,\  $3\times 2D_{12}$.
\end{itemize}

Surfaces specializing to a surface of type II have the following subgroups:

III: $4$,  $12$,   $2^2$, $2^2\times 3$, $6$, $3^2$, $3^2\times 2$, $Q_8$, $Q_8\times 3$,  $2\times D_{6}$, $D_{6}\times 3, D_6\times 6$.

VI: \  $4$,   $2^2$,  $3^2$,  $Q_8$,  $2\times 6$,  $2\times D_{6}$,  $2D_{12} \cong (2\times 6)^\bullet 2$.

X: \ $2$, $2\times D_6$.

XI: \ $2, \ 3, \ 6, \ 2^2,\  2\times 6$.

XII: \ $6\times 2, 6, 2^2, 2\times 6.$

XVI:\ $2,\  4, \ D_8$.

XVII:\ $2, \ 3,\  6$.

XVIII:\ $2, \ 6$.

XX: \ $2, \ 2^2$.

XXI:\ $2$.

\smallskip
Type IV: $P = 5$
This is easy.  We have $P \cong 5$. It admits an isomorphic lift to a nonminimal subgroup. \begin{itemize}
\item $P = \{1\}$:\ $\la \beta\alpha \ra \cong 6, \la \alpha \ra \cong 3, \la \alpha \ra \cong 2;$
\item $P = 5$: \  $5$,  $10, 15, 30$;
\end{itemize}

Surfaces specializing to a surface of type IV have the following subgroups:

XIII:\ $5, \ 10$.

XXI: \ $2$.

\smallskip
Type VII:
 $P \cong 12$.

 \begin{itemize}
\item $P = 2$:\ $2^2$.
\item $P = 3$:\ $6$;
\item $P = 4$:\ $2\times 4$;
\item $P = 6$:\ $2\times 6$;
\item $P = 12$:\ $2\times 12$.
\end{itemize}

Surfaces specializing to a surface of type VII have the following subgroups:

XII:\ $2, \  6, \ 2\times 6$.
XX: \ $2, \ 2^2$.

XXI:\ $2$.

\smallskip
Type VIII:  $P \cong 10$.
  \begin{itemize}
\item $P = \{1\}$: \ $\la \beta\alpha \ra \cong 6, \la \alpha \ra \cong 3, \la \alpha \ra \cong 2;$
\item $P = 2$:\ $2^2$.
\item $P = 5$:\ $10$;
\item $P = 10$:\ $20$.
\end{itemize}

Surfaces specializing to a surface of type VIII  have the following subgroups:

XIII:\ $5, \ 10$.

XXI: \ $2$.

\smallskip
Type XV:  $P \cong 4$.

 \begin{itemize}
\item $P = \{1\}$:\  $\la \beta\alpha \ra \cong 6, \la \alpha \ra \cong 3, \la \alpha \ra \cong 2;$
\item $P = 2$:\ $2^2$.
\item $P = 4$:\ $2\times 4$.
\end{itemize}

Surfaces specializing to a surface of type VIII  have the following subgroups:

XX: \ $2, \ 2^2$.

XXI:\ $2$.

\section{Elementary links and factorization theorem}
\subsection{Noether-Fano inequality}  Let
$|d\ell-m_1x_1-\cdots -m_Nx_N|$ be a homaloidal net in $\bbP^2$. The following is a well-known classical result.

\begin{lemma} (Noether's inequality)\label{noetherineq}  Assume
$d> 1, m_1\ge \cdots\ge m_N\ge 0$. Then
$$m_1+m_2+m_3\ge d+1,$$
and the equality holds if and only if either $m_1=\cdots=m_N$ or $m_1=n-1,m_2=\cdots=m_N$.
\end{lemma}

\begin{proof}
We have
$$m_1^2+\cdots+m_N^2 = d^2-1,\quad m_1+\cdots+m_N = 3d-3.$$
Multiplying the second equality by $m_3$ and subtracting from the first one, we get
$$m_1(m_1-m_3)+m_2(m_2-m_3)-\sum_{i\ge 4}m_i(m_3-m_i) = d^2-1-3m_3(d-1).$$
From this we obtain
$$(d-1)(m_1+m_2+m_3-d-1) = (m_1-m_3)(d-1-m_1)$$
$$+(m_2-m_3)(d-1-m_2)+\sum_{i\ge 4}m_i(m_3-m_i).$$
Since $d-1-m_i\ge 0$, this obviously proves the assertion.
\end{proof}

\begin{corollary}\label{fano2}
$$m_1 > d/3.$$
\end{corollary}

Let us generalize Corollary \ref{fano2} to birational maps of any rational surfaces. The same idea works even for higher-dimensional varieties. Let $\chi:S\to S'$ be a birational map  of surfaces. Let $\sigma:X\to S, \phi:X\to S'$ be its resolution. Let $|H'|$ be a linear system on $S'$ without base points.
Let
$$\phi^*(H') \sim \sigma^*(H)-\sum_i m_i \calE_i$$
for some divisor $H$ on $S$ and exceptional curves $\calE_i$ of the map $\sigma$.

\begin{theorem}(Noether-Fano inequality)\label{NF} Assume that there exists some integer $m_0\ge 0$ such that $|H'+mK_{S'}|=\emptyset$ for $m\ge m_0$. For any $m\ge m_0$ such that $|H+mK_S|\ne \emptyset$ there exists $i$ such that
$$m_i > m.$$
\end{theorem}

\begin{proof} We know that
$K_X = \sigma^*(K_S)+\sum_i \calE_i.$ Thus we have the equality in $\Pic(X)$
$$\phi^*(H')+mK_X = (\sigma^*(H+mK_S))+\sum (m-m_i)\calE_i.$$
Applying $\phi_*$ to the left-hand side we get the divisor class
$H'+mK_{S'} $ which, by assumption, cannot be effective. Since $|\sigma^*(H+mK_S)|\ne\emptyset$, applying $\phi_*$ to the right-hand side, we get the sum of an effective divisor and the image of the divisor $\sum_i(m-m_i)\calE_i$. If all $m-m_i$ are nonnegative, it is also an effective divisor, and we get a contradiction. Thus there exists $i$ such that $m-m_i < 0$.
\end{proof}

\begin{example} Assume $S=S' = \bbP^2$, $H = d\ell$ and $H' = \ell$. We have
$|H'+K_{S'}| = |-2\ell| = \emptyset$. Thus we can take $m_0 = 1$. If $d\ge 3$, we have for any
$1\le a\le d/3$,
$|H'+aK_S| = |(d-3a)\ell|\ne \emptyset$. This gives $m_i >  d/3$ for some $i$. This is Corollary \ref{fano2}.
\end{example}

\begin{example}\label{ex7.5}
Let $\chi:S- \to S'$ be a birational map of Del Pezzo surfaces. Assume that $S'$ is not a quadric or  the plane.  Consider the complete  linear system $H' = |-K_{S'}|$. Then $|H'+mK_{S'}| =\emptyset$ for $m \ge 2$. Let $\chi^{-1}(H') = |D-\eta|$ be its proper transform on $S$.
Choose a standard basis $(e_0,\ldots,e_k)$ in $\Pic(S)$ corresponding to the blowup $S\to \bbP^2$. Since $K_S = -3e_0+e_1+\cdots+e_k$, we can write
$\chi^{-1}(H')= |-aK_S-\sum m_ix_i|$, where $a\in \frac{1}{3}\bbZ.$ Assume that
$\chi^{-1}(H') = -aK_S$. Then there exists a point with multiplicity $> a$ if $a > 1$ that we assume. \end{example}

\begin{remark} The Noether inequality is of course well-known (see, for example, \cite{Alberich}, \cite{Hudson}). We give it here to express  our respect of  this important and beautiful result of classical algebraic geometry. Its generalization from Theorem \ref{NF} is also well-known (see \cite{Isk3}, 1.3). Note that the result can be also applied to $G$-equivariant maps $\chi$ provided that the linear system $|H'|$ is $G$-invariant. In this case the linear system $|H-\eta|$ is also $G$-invariant and the bubble cycle $\eta = \sum m_ix_i$ consists of the sum of $G$-orbits.

The existence of base points of high multiplicity in the linear system $|H-\eta| =\chi^{-1}(H')$ follows from the classical theory of termination of the adjoint system for rational surfaces which goes back to G. Castelnuovo. Nowadays this theory has  an elegant interpretation in the Mori theory which we give in the next section.
\end{remark}

\subsection{Elementary links}
We will be dealing with minimal Del Pezzo $G$-surfaces or minimal conic bundles $G$-surfaces. In the $G$-equivariant version of the Mori theory they are interpreted as extremal contractions $\phi:S\to C$, where $C = \text{pt}$ is a point in the first case and $C\cong \bbP^1$ in the second case. They are also two-dimensional analogs of rational Mori $G$-fibrations.

A birational $G$-map between  Mori fibrations are diagrams
\begin{equation}\label{link1}
\xymatrix{S{}\ar[d]_\phi\ar@{-->}[r]^\chi&S'\ar[d]_{\phi'}\\
C&C'}
\end{equation}
which in general do not commute with the fibrations. These maps are decomposed into \emph{elementary links}. These links are divided into the four following types.

\begin{itemize}
\item Links of type I:
\end{itemize}
They are commutative diagrams of the form
\begin{equation}\label{link2}
\xymatrix{S{}\ar[d]_\phi&Z=S'\ar[d]_{\phi'}\ar[l]_\sigma\\
C=\text{pt}&C'= \bbP^1\ar[l]_\alpha}
\end{equation}
Here $\sigma:Z\to S$ is the blowup of a $G$-orbit, $S$ is a minimal Del Pezzo surface, $\phi':S'\to \bbP^1$ is a minimal conic bundle $G$-fibration, $\alpha$ is the constant map. For example, the blowup of a $G$-fixed point on $\bbP^2$ defines a minimal conic $G$-bundle $\phi':\bfF_1\to \bbP^1$ with a $G$-invariant exceptional section.

\begin{itemize}
\item Links of type II:
\end{itemize}
They are commutative diagrams of the form
\begin{equation}\label{link2}
\xymatrix{S{}\ar[d]_\phi&Z\ar[l]_\sigma\ar[r]^\tau&S'\ar[d]_{\phi'}\\
C&=&C'}
\end{equation}
Here $\sigma:Z\to S, \tau:Z\to S'$ are the blowups of $G$-orbits such that
$\rank~\Pic(Z)^G = \rank~\Pic(S)^G+1 = \rank~\Pic(S')^G+1$, $C=C' $ is either a point or $\bbP^1$. An example of a link of type II is the Geiser (or Bertini) involution of $\bbP^2$, where one blows up 7 (or 8) points in general position which form one $G$-orbit. Another frequently used link of type II is an elementary transformation of minimal ruled surfaces and conic bundles.

\begin{itemize}
\item Links of type III:
\end{itemize}
These are the birational maps which are the inverses of links of type I.

\begin{itemize}
\item Links of type IV:
\end{itemize}
They exist when $S$ has two different structures of $G$-equivariant conic bundles. The link is the exchange of the two conic bundle structures
\begin{equation}\label{link4}
\xymatrix{S{}\ar[d]_\phi&=&S'\ar[d]_{\phi'}\\
C&&C'}
\end{equation}
One uses these links to relate elementary links with respect to one conic fibration to elementary links with respect to another conic fibration. Often the change of the conic bundle structures is realized via an involution in $\Aut(S)$, for example, the switch of the factors of $S= \bbP^1\times \bbP^1$ (see the following classification of elementary links).

\subsection{The factorization theorem}
Let $\chi:S-\to S'$ be a $G$-equivariant birational map of minimal $G$-surfaces. We would like to decompose it into a composition of elementary links. This is achieved with help of $G$-equivariant theory of log-pairs $(S,D)$, where $D$ is a $G$-invariant $\bbQ$-divisor on $S$. It is chosen as follows. Let us fix a $G$-invariant very ample linear system $H'$ on $S'$. If $S'$ is a minimal Del Pezzo surface we take $\calH'  = |-a'K_{S'}|, a'\in \bbZ_+$. If $S'$ is a conic bundle we take $\calH'= |-a'K_{S'}+b'f'|$, where $f'$ is the class of a fibre of the conic bundle, $a',b'$ are some appropriate positive integers.

Let $\calH= \calH_S= \chi^{-1}(\calH')$ be the proper transform of $\calH'$ on $S$. Then
$$\calH = |-aK_S-\sum m_xx|,$$
if $S$ is a Del Pezzo surface, $a\in \half\bbZ_+\cup \frac{1}{3}\bbZ_+$, and
$$\calH = |-aK_S+bf-\sum m_xx|,$$
if $S$ is a conic bundle, $a\in \half \bbZ_+, b \in \half \bbZ.$ The linear system $\calH$ is $G$-invariant, and the $0$-cycle $\sum m_xx$ is a sum of  $G$-orbits with integer multiplicities. One uses the  theory of log-pairs  $(S,D)$, where $D$ is a general divisor from the linear system $\calH$, by  applying some ``untwisting links'' to $\chi$ in order to decrease the number $a$, the algebraic degree of $\calH$. Since $a$ is a rational positive number with bounded denominator, this process terminates after finitely many steps (see \cite{Corti},\cite{Isk3}).

\begin{theorem}\label{fact} Let $f:S-\to S'$ be a birational map of minimal $G$-surfaces. Then $\chi$  is equal to a composition of elementary links.
\end{theorem}

The proof of this theorem is the same as in the arithmetic case (\cite{Isk3}, Theorem 2.5). Each time one chooses a link to apply and the criterion used for termination of the process is based on the following version of Noether's inequality in  Mori theory.

\begin{lemma}\label{newN} In the notation from above, if $m_x\le a$ for all base points $x$ of $\calH$ and $b\ge 0 $ in the case of conic bundles, then $\chi$ is an isomorphism.
\end{lemma}

The proof of this lemma is the same as in the arithmetic case (\cite{Isk3}, Lemma 2.4).

We will call a base points $x$ of $\calH$ with $m_x > a$ a \emph{maximal singularity} of $\calH$. It  follows from 3.2  that if $\calH$ has a maximal singularity of height $> 0$, then it also has a maximal singularity of height 0. We will be applying the ``untwisting links'' of types I-III to these points. If $\phi:S \to \bbP^1$ is a conic bundle with all its maximal singularities untwisted with helps of links of type II, then either the algorithm terminates, or $b < 0$. In the latter case the linear system
$|K_S+\frac{1}{a}\calH| = |\frac{b}{a}f|$ is not nef and has \emph{canonical singularities} (i.e. no maximal singularities). Applying the theory of log-pairs to the pair $(S,|\frac{b}{a}f|)$ we find an extremal contraction $\phi':S\to \bbP^1$, i.e. another conic bundle structure on $S$.  Rewriting $\calH$ in a new basis $-K_S,f'$ we find the new coefficient $a' < a$. Applying the link of type $IV$ relating $\phi$ and $\phi'$, we start over the algorithm with decreased $a$.

It follows from the proofs of Theorem \ref{fact} and Lemma \ref{newN} that all maximal singularities of $H$ are in general position in the following sense.

\begin{itemize}
\item[(i)] If $S$ is a minimal Del Pezzo $G$-surface, then the blowup of all maximal singularities of $\calH$ is a Del Pezzo surface (of course this agrees with the description of points in general position at the end of section 3.8).
\item[(ii)] If $\phi:S\to \bbP^1$ is a conic bundle, then none of the maximal singularities lie on a singular fibre of $\phi$ and no two lie on one fibre.
\end{itemize}

The meaning of these assertions is that the linear system $|H|$ has no fixed components. In the case of Del Pezzo surfaces with an orbit of maximal singular points we can find a link by blowing up this orbit to obtain a surface $Z$ with $\Pic(Z)^G \cong \bbZ\oplus \bbZ$ and two extremal rays. By applying Kleiman's criterion this implies that $-K_Z$ is ample. The similar situation occurs in the case of conic bundles (see \cite{Isk3}, Comment 2).

Let $S$ be a minimal Del Pezzo $G$-surface of degree $d$. Let us write $\calH_S = \\|-aK_S-\sum m_\kappa\kappa|$ as in \eqref{kappa}.

\begin{lemma}\label{7.11}  Let $\kappa_1,\ldots,\kappa_n$ be the $G$-orbits of maximal multiplicity. Then
$$\sum_{i=1}^n d(\kappa_i) < d.$$
\end{lemma}

\begin{proof}  Let $D_1,D_2\in \calH_S$ be two general divisors from $\calH_S$. Since $\calH_S$ has no fixed components, we have
$$0\le D_1\cdot D_2 = a^2d-\sum m_\kappa^2d(\kappa) \le a^2d-\sum_{i=1}^nm_{\kappa_i}^2d(\kappa_i) = $$
$$a^2(d-\sum_{i=1}^nd(\kappa_i))-\sum_{i=1}^n(m_{\kappa_i}^2-a^2)d(\kappa_i).$$
It follows from Example \ref{ex7.5} that $m_{\kappa_i}> a$ for all $i = 1,\ldots,n$. This implies
that  $d > \sum_{i=1}^nd(\kappa_i)$.

\end{proof}

\begin{definition}\label{rigid} A minimal Del Pezzo $G$-surface is called \emph{superrigid} (resp. \emph{rigid})  if any birational $G$-map $\chi:S-\to S'$ is a $G$-isomorphism (resp. there exists a birational $G$-automorphism $\alpha:S-\to S$ such that $\chi\circ \alpha$ is a $G$-isomorphism).

A minimal conic bundle $\phi:S\to \bbP^1$ is called \emph{superrigid} (resp. \emph{rigid}) if for any birational $G$-map $\chi:S-\to S'$, where $\phi':S'\to \bbP^1$ is a minimal conic bundle, there exists an isomorphism $\delta:\bbP^1\to \bbP^1$ such that the following diagram is commutative
\begin{equation}\label{link}
\xymatrix{S\ar[d]_\phi\ar@{-->}[r]^\chi&S'\ar[d]_{\phi'}\\
\bbP^1\ar[r]^\delta&\bbP^1}
\end{equation}
(resp. there exists a birational $G$-automorphism $\alpha:S-\to S'$ such that the diagram is commutative after we replace $\chi$ with $\chi\circ \alpha$).
\end{definition}

Applying Lemma \ref{newN} and Lemma \ref{7.11}, we get the following.

\begin{corollary}\label{c7.11}
Let $S$ be a minimal Del Pezzo $G$-surface of degree $d = K_S^2$. If $S$ has no $G$-orbits $\kappa$ with $d(\kappa) < d$, then $S$ is superrigid. In particular, a Del Pezzo surface of degree 1  is always superrigid and a Del Pezzo surface of degree 2 is superrigid unless $G$ has a fixed point.

A minimal conic $G$-bundle with $K_S^2\le 0$ is superrigid.
\end{corollary}

The first assertion is clear. To prove the second one, we untwist all maximal base points of $\calH_S$ with help of links of type II to get a conic bundle $\phi':S'\to \bbP^1$ with $b' < 0$. Since
$H_{S'}^2 = a^2K_{S'}^2+4ab'-\sum m_x'{}^2 \ge 0$ and $K_{S'}^2 = K_S^2 \le 0, 4ab' < 0$, we get a contradiction with Lemma \ref{newN}. Thus $\chi$ after untwisting maximal base points terminates at an isomorphism (see \cite{Isk1}, \cite{Isknew}, \cite{Isk3}, Theorem 1.6).

\subsection{Classification of elementary links}

Here we consider an elementary link $f:S- \to S'$ defined by a resolution $(S\overset{\sigma}{\leftarrow} Z\overset{\tau}{\rightarrow} S')$. We take $H_{S'}$ to be the linear system $|-aK_{S'}|$ if $S'$ is a Del Pezzo surface and $|f|$ if $S'$ is a conic bundle, where $f$ is the divisor class of a fibre. It is assumed that the point which we blow up are in general position in sense of the previous subsection.

We denote by $\calD P_k$  (resp. $\calC_k$) the set of isomorphism classes of minimal Del Pezzo surfaces (resp. conic bundles) with $k = K_S^2$ (resp. $k = 8 -K_S^2$).

\begin{proposition}\label{7.12} Let $S,S'$ be as in Link I of type I. The map $\sigma:Z=S'\to S$ is the blowing up of a $G$-invariant bubble cycle $\eta$ with $\htt(\eta) = 0$ of some degree $d$. The proper transform of the linear system $|f|$ on $S'$ is equal to the linear system $\calH_S = |-aK_S-m\eta|$.  Here  $f$ is the class of a fibre of the conic bundle structure on $S'$. The following cases are possible:
\begin{enumerate}
\item $K_S^2 = 9$
\begin{itemize}
\item $S = \bbP^2, S' = \bfF_1, d = 1, m = 1, a = \frac{1}{3}.$
\item  $S = \bbP^2, S' \in \calC_3, d = 4, m = 1,  a = \frac{2}{3}$.
\end{itemize}
\item $K_S^2 = 8$
\begin{itemize}
\item $S = \bfF_0, \pi:S' \to \bbP^1$ a conic bundle with $k = 2$, $d= 2$, $m = 1$, $a = \frac{1}{2}$.
\end{itemize}
\item $K_S^2 = 4$
\begin{itemize}
\item $S \in \calD P_4, p:S' \to \bbP^1$ a conic bundle with $f = -K_{S'}-l$, where $l$ is a (-1)-curve,
$d = a = 1, m = 2$.
\end{itemize}
\end{enumerate}
\end{proposition}

\begin{proof} Let $\calH_S = |-aK_S-b\eta|$, where $\eta$ is a $G$-invariant bubble cycle of degree $d$. We have
$$(-aK_S-b\eta)^2 = a^2K_S^2-b^2d = 0,\quad (-aK_S-b\eta,-K_S) = aK_S^2-bd = 2.$$
Let $t = b/a$. We have
$$(td)^2= dK_S^2,  \quad K_S^2-td = 2/a > 0.$$
The second inequality, gives $td < K_S^2$, hence $d < K_S^2$. Giving the possible values for $K_S^2$  and using that $a\in \frac{1}{3}\bbZ$,  we check that the only possibilities are:
$$(K_S^2,d,t)  = (9,1,3),\  (8, 2, 2), \ (4,1,2), \ (4,2,1).$$
This gives our cases and one extra case $(4,2,1)$. In this case $a = 2$ and $\calH_S = |-2K_S-2x_1|$ contradicting the primitiveness of $f$. Note that this case is realized in the case when the ground field is not algebraically closed (see \cite{Isk3}).
\end{proof}

\begin{proposition}\label{7.16} Let $S,S'$ be as in Link of type II. Assume that $S,S'$ are both minimal Del Pezzo surfaces. Then $(S\overset{\sigma}{\leftarrow} Z\overset{\tau}{\rightarrow} S')$, where $\sigma$  is the blowup of a $G$-invariant bubble cycle $\eta$ with $\htt(\eta) = 0$ and some degree $d$. The proper transform of the linear system $|-K_{S'}|$ on $S$ is equal to $|-aK_S-m\eta|$.  And similarly defined $d',m',a'$ for $\tau$. The following cases are possible:
\begin{enumerate}
\item $K_S^2 = 9$
\begin{itemize}
\item $S' \cong S = \bbP^2,  d = d' =8, m = m' = 18, a = a' = 17$ ($S\leftarrow Z\rightarrow S')$ is a minimal resolution of a Bertini involution).
\item  $S' \cong  S = \bbP^2,  d = d' = 7,  m = m' = 9, a = a' = 8$ ($S\leftarrow Z\rightarrow S')$ is a minimal resolution of a Geiser transformation).
\item $S' \cong S = \bbP^2,  d = d' = 6, m = m' = 6, a = a'= 5$ ($S\leftarrow Z\rightarrow S')$ is a minimal resolution of a  Cremona transformation given by the linear system  $|5\ell-2p_1-2p_2-2p_3-2p_4-2p_5|$,
\item  $S \cong \bbP^2,  S' \in \calD P_5, d = 5, m' = 6, a = \frac{5}{3}, d' = 1, m = 2, a' = 3$.
\item $S \cong S' = \bbP^2, d = d' = 3, m = m' = 1, a = a' = \frac{2}{3},$
$(S\leftarrow Z\rightarrow S')$ is a minimal resolution of a standard quadratic transformation.
\item $S =  \bbP^2, S' =\bfF_0,  d = 2, m = 3, a' = \frac{3}{2},  d' = 1, a = \frac{4}{3}$.
 \end{itemize}
 \item $K_S^2 = 8$
 \begin{itemize}
 \item  $S \cong    S' \cong \bfF_0,  d =   d' = 7, a = a' = 15, m = m' = 16$.
 \item  $S \cong   S' \cong \bfF_0,  d =   d' = 6, a = a' = 7, m = m' = 8$.
  \item  $S \cong   \bfF_0,  S'\in \calD P_5, d = 5,  d' = 2, a = \frac{5}{2}, m = 4, a' = 4, m' =  6$.
  \item $S \cong  \bfF_0, S' \cong  \bfF_0, d =  d' = 4, a = a' = 3, m = m' = 4$.
  \item  $S \cong  \bfF_0, S' \in \calD P_6, d = 3,  d' = 1, a = \frac{3}{2}, m = 2, m' = 4, a' = 2$.
 \item $S \cong  \bfF_0, S' \cong \bbP^2, d = 1,  d' = 2, a = \frac{3}{2}, m = 3, a' = \frac{4}{3}, m' = 2$. This link is the inverse of the last case from the preceding list.
 \end{itemize}
 \item $K_S^2 = 6$
\begin{itemize}
\item  $S \cong S' \in \calD P_6,  d =  d' = 5, a = 11, m = 12$.
\item  $S \cong  S' \in \calD P_6, d =  d' = 4, a = 5, m = 6$.
\item  $S \cong  S' \in \calD P_6,  d =  d' = 3, a = 3, m = 4$.
\item  $S \cong S' \in \calD P_6,  d =  d' = 2, a = 2, m = 3$.
\item $S  \in \calD P_6, S' = \bfF_0, d =  1, d' = 3, a = \frac{3}{2}, m = 2$. This link is the inverse of the link from the preceding list with $S'\in \calD P_6, d = 3$.
\end{itemize}
\item $K_S^2 = 5$
\begin{itemize}
\item  $S \cong  S' \in \calD P_5,  d =  d' = 4, m= m' = 10, a=a'= 5.$.
\item  $S = S' \in \calD P_5,  d =  d' = 3, m= m' = 5, a=a'= 4$.
\item $S \in \calD P_5, S' = \bfF_0, d =  2, d' = 5$. This link is inverse of the link with $S = \bfF_0, S'\in \calD P_5, d = 5$.
\item $S  \in \calD P_5,  S' = \bbP^2, d =   1, d' = 5$. This link is inverse of the link with $S = \bbP^2, S'\in \calD P_5, d = 5$.

\end{itemize}
\item $K_S^2 = 4$
\begin{itemize}
\item $S \cong S' \in \calD P_4, d =   d' = 3$. This is an analog of the Bertini involution.
\item $S \cong  S'\in \calD P_4,  d =   d' = 2$. This is an analog of the Geiser involution.
\end{itemize}
\item $K_S^2 = 3$
\begin{itemize}
\item $S \cong S'\in \calD P_3, d = d' = 2.$ This is an analog of the Bertini involution.
\item $S \cong S'\in \calD P_3, d = d' = 1.$ This is an analog of the Geiser involution.
\end{itemize}
\item $K_S^2 = 2$
\begin{itemize}
\item $S = S'\in \calD P_2, d = d' = 1$. This is an analog of the Bertini involution.
\end{itemize}
\end{enumerate}
\end{proposition}

\begin{proof}
Similar to the proof  of the previous proposition, we use that
$$\calH_S^2 = a^2K_S^2-b^2d = K_{S'}^2, \quad aK_S^2-bd = K_{S'}^2,$$
$$H_{S'}^2 = a'{}^2K_{S'}^2-b'{}^2d = K_{S}^2, \quad a'K_S^2-b'd = K_{S}^2.$$
Since the link is not a biregular map, by Noether's inequality we have $a > 1, a' > 1, b > a, b' > a'$. This implies
$$d < K_S^2-\frac{1}{a}K_{S'}^2, \  d' < K_{S'}^2-\frac{1}{a'}K_{S}^2.$$
It is not difficult to list all solutions. For example, assume $K_S^2 = 1$. Since $d$ is a positive integer, we see that there are no solutions. If $K_S^2 = 2$, we must have $d= d' =1$.
\end{proof}

\begin{proposition}\label{7.17} Let $S,S'$ be as in Link of type II. Assume that $S,S'$ are both minimal conic bundles. Then $(S\leftarrow Z\rightarrow S')$ is a composition of  elementary transformations
$\elm_{x_1}\circ \ldots \circ \elm_{x_s}$, where $(x_1,\ldots,x_s)$ is a  $G$-orbit of points not lying on a singular fibre with no two points lying on the same fibre.
\end{proposition}

We skip the classification of links of type III. They are the inverses of links of type I.

\begin{proposition}\label{7.18} Let $S,S'$ be as in Link of type IV. Recall that they consist of changing the conic bundle structure. The following cases are possible:
\begin{itemize}
\item $K_S^2 = 8, S' = S, f' = -K_{S'}-f$, it is represented by a switch automorphism;
\item $K_S^2 = 4, S' = S, f' = -K_{S'}-f$;
\item  $K_S^2 = 2, S' = S, f '= -2K_{S'}-f$; it is represented by a Geiser involution;
\item $K_S^2 = 1, S' = S , f '= -4K_{S'}-f $; it is represented by a Bertini involution;
\end{itemize}
\end{proposition}

\begin{proof} In this case $S$ admits two extremal rays and $\rank~\Pic(S)^G  = 2$ so that $-K_S$ is ample. Let $|f'|$ be the second conic bundle. Write $f'\sim -aK_S+bf$. Using that  $f'{}^2 = 0, f\cdot K_S = f'\cdot K_S = -2$, we easily get $b = -1$ and $aK_S^2= 4$. This gives all possible cases from the assertion.
\end{proof}

\section{Birational classes of minimal $G$-surfaces}
\subsection{} Let $S$ be a minimal $G$-surface $S$ and $d = K_S^2$. We will classify all isomorphism classes of $(S,G)$ according to the increasing parameter $d$. Since the number of singular fibres of a minimal conic bundle is at least 4, we have $d\le 4$ for conic bundles.

\begin{itemize}
\item $d\le 0$.
\end{itemize}

By Corollary \ref{c7.11}, $S$ is a superrigid conic bundle with $k = 8-d$ singular fibres. The number $k$ is a birational invariant. The group $G$ is of de Jonqui\`eres type and its conjugacy class in $\Cr(2)$ is determined uniquely by Theorem \ref{itog} or Theorem \ref{autoex}.

Also observe that if $\phi:S\to \bbP^1$ is an exceptional conic bundle and $G_0 =\Ker(G\to \Or(\Pic(S))$ is non-trivial, then no links of type II is possible. Thus the conjugacy class of $G$ is uniquely determined by the isomorphism class of $S$.

\begin{itemize}
\item $d = 1$, $S$ is a Del Pezzo surface.
\end{itemize}

By Corollary \ref{c7.11}, the surface $S$ is superrigid. Hence the conjugacy class of $G$ is determined uniquely by its conjugacy class in $\Aut(S)$. All such conjugacy classes are listed in subsection \ref{clasdp1}.

\begin{itemize}
\item $d = 1$, $S$ is a conic bundle.
\end{itemize}

Let $\phi:S\to \bbP^1$ be a minimal conic bundle on $S$. It has $k = 7$ singular fibres.
If $-K_S$ is ample, i.e. $S$ is a (nonminimal) Del Pezzo  surface, then the center of $\Aut(S)$ contains the Bertini involution $\beta$. We know that  $\beta$ acts as $-1$ on $\calR_S$, thus $\beta$ does not act identically on $\Pic(S)^G$, hence $\beta\not\in G$. Since $k$ is odd, the conic bundle is not exceptional, so we cam apply Theorem \ref{itog}, Case (2). It follows that $G$ must contain a subgroup isomorphic to $2^2$, adding $\beta$ we get that $S$ is a minimal Del Pezzo $2^3$-surface of degree 1. However, the classification shows that there are no such surfaces.

Thus $-K_S$ is not ample. It follows from Proposition \ref{7.16} that the structure of a conic bundle on $S$ is unique. Any other conic bundle birationally $G$-isomorphic to $S$ is obtained from $S$ by elementary transformations with $G$-invariant set of centers.

\begin{itemize}
\item $d = 2$, $S$ is a Del Pezzo surface.
\end{itemize}

By Corollary \ref{c7.11}, $S$ is superrigid unless $G$ has a fixed point on $S$. If $\chi:S-\to S'$ is a birational $G$-map, then $H_S$ has only one maximal base point and it does not lie on a $(-1)$-curve. We can apply an elementary link  $Z\to S, Z\to S$ of type II which together with the projections $S\to \bbP^2$ resolves the Bertini involution. These links together with the $G$-automorphisms (including the Geiser involution) generate the group of birational $G$-automorphisms of $S$ (see \cite{Isk3}, Theorem 4.6). Thus the surface is rigid. The conjugacy class of $G$ in $\Cr(2)$ is determined uniquely by the conjugacy class of $G$ in $\Aut(S)$. All such conjugacy classes are listed in Table \ref{tab6} and Theorem \ref{clasdp2}.

\begin{itemize}
\item $d = 2$, $\phi:S \to \bbP^1$ is a conic bundle.
\end{itemize}
If $-K_S$ is ample, then  $\phi$ is not exceptional. The center of $\Aut(S)$ contains the Geiser involution $\gamma$. Since $\gamma$ acts non-trivially on $\Pic(S)^G = \bbZ^2$, we see that $\gamma\not\in G$.   Applying $\gamma$ we obtain another conic bundle structure. In other words, $\gamma$ defines an elementary link of type IV. Using the factorization theorem we show that the group of birational $G$-automorphisms of $S$ is generated by links of type II, the Geiser involution, and $G$-automorphisms (see \cite{Isknew}, \cite{Isk4}, Theorem 4.9). Thus $\phi:S\to \bbP^1$ is a rigid conic bundle.

If $S$ is not a Del Pezzo surface, $\phi$ could be an exceptional conic bundle with $g= 2$. In any case the group $G$ is determined in Theorem \ref{autoex}. We do not know whether $S$ can be mapped to a conic bundle with $-K_S$ ample (see \cite{Isknew}).

Applying Proposition \ref{cdp}, we obtain that any conic bundle with $d\ge 3$ is a nonminimal Del Pezzo surface, unless $d = 4$ and $S$ is an exceptional conic bundle. In the latter case,  the group $G$ can be found in Theorem \ref{autoex}. It is not known whether it is birationally $G$-isomorphic to a Del Pezzo surface. It is true in the arithmetic case.

\begin{itemize}
\item $d = 3$, $S$ is a minimal Del Pezzo surface.
\end{itemize}

The classification of elementary links shows that $S$ is rigid. Birational $G$-automorphisms are generated by links of type (6) from Proposition \ref{7.12}. The conjugacy class of $G$ in $\Cr(2)$ is determined by the conjugacy class of $G$ in $\Aut(S)$.

\begin{itemize}
\item $d = 3$, $S$ is a minimal conic bundle.
\end{itemize}
 Since $k = 5$ is odd, $G$ has 3 commuting involutions, the fixed-point locus of one of them must be a rational 2-section of the conic bundle. It is easy to see that it is a $(-1)$-curve $C$ from the divisor class $-K_S-f$. The other two fixed-point curves are of genus 2. The group $G$ leaves the curve $C$ invariant. Thus blowing it down, we obtain a minimal Del Pezzo $G$-surface $S'$ of degree 4. The group $G$ contains a subgroup isomorphic to $2^2$. Thus $G$ can be found in the list of minimal groups of degree 4 Del Pezzo surfaces with a fixed point. For example, the group $2^2$ has 4 fixed points.

\begin{itemize}
\item $d = 4$, $S$ is a minimal Del Pezzo surface.
\end{itemize}

If $S^G =\emptyset$, then $S$ admits only self-links of type II, so it is rigid or superrigid. The conjugacy class of $G$ in $\Cr(2)$ is determined by the conjugacy class of $G$ in $\Aut(S)$ and we can apply Theorem \ref{clasdp4}. If $x$ is a fixed point of $G$, then we can apply a link of type I, to get a minimal conic bundle with $d = 3$. So, all groups with $S^G\ne \emptyset$ are conjugate to groups of de Jonqui\`eres type realized on a conic bundle $S\in \calC_5$.

\begin{itemize}
\item $d = 4$, $S$ is a  minimal conic bundle.
\end{itemize}

Since $k = 4$, it follows from Lemma \ref{ineqq} that either $S$ is an exceptional conic bundle with $g = 1$, or $S$ is a Del Pezzo surface with two  sections with self-intersection $-1$ intersecting at one point. In the latter case, $S$ is obtained by regularizing a de Jonqu\'eres involution $IH_3$ (see section \ref{dejon}).  In the case when $S$ is an exceptional conic  bundle  the groups of automorphisms  are described in Example \ref{exept}. They are minimal if and only if  the kernel of the map $G \to \PGL(2)$ contains an involution not contained in $G_0= \Ker(G\to \Or(\Pic(S))$.  If  $G_0$  is not trivial, then no elementary transformation is possible. So, $S$ is not birationally isomorphic to a Del Pezzo surface.

\begin{itemize}
\item $d = 5$, $S$ is a Del Pezzo surface, $G \cong 5$.
\end{itemize}

Let us show that $(S,G)$ is birationally isomorphic to $(\bbP^2,G)$. Since rational surfaces are simply-connected, $G$ has a fixed point $x$ on $S$. The anti-canonical model of $S$ is a surface of degree $5$ in $\bbP^5$. Let $P$ be the tangent plane of $S$ at $x$. The projection from $P$ defines a birational $G$-equivariant map  $S-\to \bbP^2$ given by the linear system of anti-canonical curves with double point at $x$. It is an elementary  link of type II.

\begin{itemize}
\item $d = 5$, $S$ is a Del Pezzo surface, $G \cong 5:2, 5:4$.
\end{itemize}

The cyclic subgroup of order 5 of $G$ has two fixed points on $S$. This immediately follows from the Lefschetz fixed-point formula. Since it is normal in $G$, the groups $G$ has an orbit $\kappa$  with $d(\kappa) = 2$. Using an elementary link of type II with $S' =\bfF_0$, we obtain that $G$ is conjugate to a group acting on $\bfF_0$.

\begin{itemize}
\item $d = 5$, $S$ is a Del Pezzo surface, $G \cong A_5,S_5$.
\end{itemize}

It is clear that $S^G = \emptyset$ since otherwise $G$ admits a faithful 2-dimensional linear representation. It is known that it does not exist. Since $A_5$ has no index 2 subgroups $G$ does not admit orbits $\kappa$ with $d(\kappa) = 2$. The same is obviously true for $G = S_5$. It follows from the classification of links that $(S,G)$ is superrigid.

\begin{itemize}
\item $d = 6 $.
\end{itemize}
One of the groups from this case, namely $G \cong  2\times S_3$  was considered in \cite{Isk4}, \cite{Isk5} (the papers also discuss the relation of this problem to some questions  in the theory of algebraic groups raised in \cite{Lemire}).  It is proved there that $(S,G)$ is not birationally isomorphic to $(\bbP^2,G)$ but birationally isomorphic to minimal $(\bfF_0,G)$. The birational isomorphism is easy to describe. We know that $G$ contains the lift of the standard Cremona involution. It has 4 fixed points in $S$, the lifts of the  points  given in affine coordinates by $(\pm 1,\pm 1)$. The group $S_3$ fixes $(1,1)$ and permutes the remaining points $p_1,p_2,p_3$. The proper transforms of the lines $\la p_i,p_j\ra$ in $S$ are disjoint $(-1)$-curves $E_i$. The anti-canonical model of $S$ is a surface of degree 6 in $\bbP^6$. The projection from the tangent plane to $S$ at the fixed point, is a link of type II with $S' = \bfF_0$. It blows-up the fixed point, and then blows down the preimages of the curves $E_i$. Now the group $G$ acts on $\bfF_0$ with $\bfF_0^G =\emptyset$.

If minimal $G$ contains some non-trivial imprimitive projective transformations, then $G$ has no fixed points. It follows from the classification of links that $S$ is rigid. If $G \cong 6$ or $S_3$, then it acts on $\bfF_0$ with a fixed point. The projection from this point defines a birational isomorphism $(S,G)$ and $(\bbP^2,G)$. Thus the only groups which are not conjugate to a group of projective transformations are the groups which are mapped surjectively  onto  $W_S = S_3\times 2$. Those of them which are mapped isomorphically are conjugate to subgroups of $\bfF_0$.

\begin{itemize}
\item $d = 8$.
\end{itemize}
In this case $S = \bfF_0$ or $\bfF_n, n> 1$. In the first case $(S,G)$ is birationally isomorphic to $(\bbP^2,G)$ if $S^G\ne\emptyset$ (via the projection from the fixed point). This implies that the subgroup $G'$ of $G$ belonging to the connected component of the identity of $\Aut(\bfF_0)$ is an extension of cyclic groups. As we saw in Theorem \ref{clas} this implies that $G'$ is an abelian group of transformations
$(x,y) \mapsto (\epsilon_{nk}^ax,\epsilon_{mk}^by),$ where $a = sb \mod k$ for some $s$ coprime to $k$.
If  $G\ne G'$, then we must have $m = n = 1$ and $s = \pm 1\mod k$.

If $\bfF_0^G = \emptyset$ and $\Pic(\bfF_0)^G \cong \bbZ$, then the classification of links shows that links of type II with $d= d' = 7, 6,5, d = 3, d' = 1$ map $\bfF_0$ to $\bfF_0$ or to minimal Del Pezzo surfaces of degrees 5 or 6.  These cases  have been already considered. If $\rank~\Pic(S)^G = 2$, then any birational $G$-map $S-\to S'$ is composed of elementary transformations with respect to one of the conic bundle fibrations.  They do not change $K_S^2$ and do not give rise a fixed points. So, $G$ is not conjugate to any subgroup of $\Aut(\bbP^2)$.

Assume $n > 1$. Then $G = A.B$, where $A\cong n$ acts identically on the base of the fibration and $B\subset \PGL(2)$. The subgroup $B$ fixes pointwise two disjoint sections, one of them is the exceptional one. Let us consider different cases corresponding to possible groups $B$.

\smallskip
\emph{$B \cong C_n$}. In this case  $B$ has two fixed points on the base, hence $G$ has 2 fixed points on the non-exceptional section. Performing an elementary transformation with center at one of these points we descend $G$ to a subgroup of $\bfF_{n-1}$. Proceeding in this way, we arrive to the case $n = 1$, and then obtain that $G$ is a group of automorphisms of $\bbP^2$.

\smallskip
\emph{$B \cong D_n$}. In this case $B$  has an orbit of cardinality 2 in $\bbP^1$. A similar argument shows that $G$ has an orbit of cardinality 2 on the non-exceptional section. Applying the product of the elementary transformations at these points we descend $G$ to a subgroup of automorphisms of $\bfF_{n-2}$. Proceeding in this way we obtain that $G$ is a conjugate to a subgroup of $\Aut(\bbP^2)$ or of $\Aut(\bfF_0)$. In the latter case it does not have fixed points, and hence is not conjugate to a linear subgroup of $\Cr(2)$.

\smallskip
$B \cong  T$. The group $B$  has an orbit of cardinality 4 on the non-exceptional section. A similar argument shows that $G$ is conjugate to a group of automorphisms of $\bfF_2, \bfF_0,$ or $\bbP^2$.

\smallskip
 $B \cong O$. The group $B$  has an orbit of cardinality 6. As in the previous case we show that $G$ is conjugate to a group of automorphisms of $\bbP^2$, or $\bfF_0$, or $\bfF_2$, or $\bfF_3$.

\smallskip
$B \cong I$. The group $B$  has an orbit of cardinality 12. We obtain that $G$ is conjugate to a group of automorphisms of $\bbP^2$ or of $\bfF_n$, where $n = 0,2,3,4,5,6$\footnote{A better argument due to I. Cheltsov allows one to show that in this and the previous  cases $B$ is conjugate to a group of automorphisms of $\bbP^2$, or $\bfF_0$, or $\bfF_2$.} 

\begin{itemize}
\item $d = 9$.
\end{itemize}
In this case $S = \bbP^2$ and $G$ is a finite subgroup of $\PGL(3)$.
The methods of the representation  theory allows us to classify them up to conjugacy in the group $\PGL(3)$. However, some of non-conjugate groups can be still conjugate inside the Cremona group.

For example   all cyclic subgroups of $\PGL(3)$ of the same order $n$ are conjugate in $\Cr(2)$.
Any element $g$ of order $n$ in $\PGL(3)$ is conjugate to a transformation $g$ given in affine coordinates by the formula $(x,y) \mapsto (\epsilon_nx,\epsilon_n^ay)$. Let $T\in \dJJ(2)$ be  given by the formula
$(x,y) \mapsto (x,x^a/y)$. Let $g':(x,y)\mapsto (\epsilon_n^{-1}x,y)$. We have
$$g'\circ T\circ g:(x,y) \mapsto (\epsilon_n x, \epsilon_n^a y) \mapsto (\epsilon_n x, x^a/y) \mapsto (x,x^a/y) = T.$$
This shows that $g'$ and $g$ are conjugate.

We do not know whether any two isomorphic non-conjugate subgroups  of $\PGL(3)$ are conjugate in $\Cr(2)$.

\section{What is left?}
Here we list some problems which  have not been yet resolved.

\begin{itemize}
\item Find the conjugacy classes in $\Cr(2)$ of subgroups of $\PGL(3)$. For example, there are two non-conjugate subgroups of $\PGL(3)$ isomorphic to $A_5$ and three to $A_6$ which differ by an outer automorphism of the groups. Are they conjugate in $\Cr(2)$?
\end{itemize}

\begin{itemize}
\item Find a finer classification of the conjugacy classes of de Jonqui\`eres groups.
\end{itemize}
We already know that the number of singular fibres in a minimal conic bundle $G$-surface is an invariant. Even more, the projective equivalence class of the corresponding $k$ points on the base of the conic fibration is an invariant. Are there other invariants? In the case when $G_K\cong 2$, we know that the quotient of the conic bundle by the involution generating $G_K$ is a minimal ruled surface $\bfF_e$. Is the number $e$ a new  invariant?

\begin{itemize}
\item  Give a finer geometric description of the algebraic variety parametrizing conjugacy classes.
\end{itemize}
Even in the case of Del Pezzo surfaces we give only normal forms. What is precisely the moduli space of   Del Pezzo surfaces with a fixed isomorphism class of a minimal automorphism group?

\smallskip
We know that conic bundles $(S,G)$ with $k\ge 8$ singular fibres are superrigid, so any finite subgroup $G'$ of $\Cr(2)$ conjugate to $G$  is realized as an automorphism group of a conic bundle obtained from $S$ by a composition of elementary transformations with $G$-invariant  centers. If $S$ is not exceptional and $G\cong 2.P$, then the invariant of the conjugacy class is the hyperelliptic curve of fixed points of the central involution. If $G\cong 2^2.P$, then we have three commuting involutions and their curves of fixed points are the invariants of the conjugacy class. Do they determine the conjugacy class?

When $k= 6, 7$ we  do not know whether $(S,G)$  is birationally isomorphic to $(S',G)$, where $S'$ is a Del Pezzo surface. This is true when $k\in \{0,1,2,3,5\}$ or $k=4$ and $S$ is not exceptional.

\begin{itemize}
\item  Find more explicit description of groups $G$ as subgroups of $\Cr(2)$.
\end{itemize}

This has been done in the case of abelian groups in
\cite{Blanc}. For example one may ask to  reprove and revise Autonne's  classification of groups whose elements are quadratic transformations \cite{Autonne}. An example of such non-cyclic group is the group of automorphisms $S_5$ of a Del Pezzo surface of degree 5.

\begin{itemize}
\item   Finish the classical work on the birational classification of rational cyclic planes $z^n=f(x,y)$.
\end{itemize}

More precisely, the quotient $S/G$ of a rational surface $S$ by a cyclic group of automorphisms defines a cyclic extension of the fields of rational functions. Thus there exists a rational function $R(x,y)$ such that there exists an isomorphism of fields $\bbC(x,y)(\sqrt[n]{R(x,y)}) \cong \bbC(x,y)$, where $n$ is the order of $G$. Obviously we may assume that $R(x,y)$ is a polynomial $f(x,y)$, hence we obtain an affine model of $S$ in the form $z^n = f(x,y)$. A birational isomorphism of $G$-surfaces replaces the branch curve $f(x,y) = 0$ by a Cremona equivalent curve $g(x,y)$. The problem is to describe the Cremona equivalence classes of the branch curves which define rational cyclic planes.

For example, when $(S,G)$ is birationally equivalent to $(\bbP^2,G)$, we may take $f(x,y) = x$ since all cyclic groups of given order are conjugate in $\Cr(2)$. When $n = 2$, the problem was solved by M. Noether \cite{Noether} and later G. Castelnuovo and F. Enriques \cite{CE} had realized that the classification follows from  Bertini's classification of involutions in $\Cr(2)$.  When $n$ is prime the problem was studied by A. Bottari \cite{Bottari}. We refer for modern work on this problem to \cite{Calabri1}, \cite{Calabri2}.

\begin{itemize}
\item Extend the classification to the case of non-algebraically closed fields, e.g. $\bbQ$, and algebraically closed fields of positive characteristic.
\end{itemize}

Note that there could be more automorphism groups in the latter case. For example, the Fermat cubic surface
$T_0^3+T_1^3+T_2^3+T_3^3 = 0$ over a field of characteristic 2 has the automorphism  group isomorphic to $U(4,\bbF_4)$, which is a subgroup of index 2 of the Weyl group $W(E_6)$.

\section{Tables}
In the following tables we give the tables of conjugacy classes of subgroups in $\Cr(2)$ which are realized as minimal automorphism groups of Del Pezzo surfaces of degree $d \le 5$ and not conjugate to subgroups of automorphisms of minimal rational surfaces or conic bundles.   The information about other  groups can be found in the corresponding sections of the paper. The tables contain the order of a group $G$, its structure, the type of a surface on which the group is realized as a minimal group of automorphisms, the equation of the surface, and the number of conjugacy classes, if finite, or the dimension of the variety parametrizing the conjugacy classes.

\clearpage
\begin{table}[t]
\begin{center}
\begin{tabular}{| l |r| r | r |r|}\hline
Order&Type&Surface&Equation&Conjugacy\\ \hline
2&$A_1^7$&$\calD P_2$&XIII&$\infty^6$\\ \hline
2&$A_1^8$&$\calD P_1$&XXII&\\ \hline
3&$3A_2$&$\calD P_3$&I,III,IV&$\infty^1$\\ \hline
3& $4A_2$&$\calD P_1$&XVIII&$\infty^3$\\ \hline
4&$2A_3+A_1$&$\calD P_2$&II,III,V&$\infty^1$\\ \hline
4&$2D_4(a_1)$&$\calD P_1$&I,VI, X,XVII,XX&$\infty^5$\\ \hline
5&$2A_4$&$\calD P_1$&XIV&$\infty^2$\\ \hline
6&$E_6(a_2)$&$\calD P_3$&I,VI&$\infty^1$\\ \hline
6&$A_5+A_1$&$\calD P_3$&I, III, IV& $\infty^1$\\ \hline
6&$E_7(a_4)$&$\calD P_2$&XI&$\infty^1$ \\ \hline
6&$A_5+A_2$&$\calD P_2$&VIII&$\infty^1$\\ \hline
6&$D_6(a_2)+A_1$&$\calD P_2$&II,III,IV,IX&$\infty^1$\\ \hline
6&$A_5+A_2+A_1$&$\calD P_1$&II,VIII,XIII& $\infty^2$\\ \hline
6&$E_6(a_2)+A_2$&$\calD P_1$&II,XII&$\infty^2$\\ \hline
6&$E_8(a_8)$&$\calD P_1$&XVIII&$\infty^3$\\ \hline
6&$2D_4$&$\calD P_1$&VII,XI&$\infty^1$\\ \hline
6&$E_7(a_4)+A_1$&$\calD P_1$&II,VIII,XIX&$\infty^4$\\ \hline
8&$D_5$&$\calD P_4$&\eqref{6.6}&1\\ \hline
8&$D_8(a_3)$&$\calD P_1$&X&$1$\\ \hline
9&$E_6(a_1)$&$\calD P_3$&I&$1$\\ \hline
10&$E_8(a_6)$&$\calD P_1$&IV,IX,XIV&$\infty^2$\\ \hline
12&$E_6$&$\calD P_3$&III&1\\ \hline
12& $E_7(a_2)$&$\calD P_2$&III&1\\ \hline
12& $E_8(a_3)$&$\calD P_1$& I,V&$\infty^2$\\ \hline
14&$E_7(a_1)$&$\calD P_2$&I&$1$\\ \hline
15&$E_8(a_5)$&$\calD P_1$&IV&1\\ \hline
18& $E_7$ &$\calD P_2$&VI&1\\ \hline
20& $E_8(a_2)$&$\calD P_1$&IX&1\\ \hline
24&$E_8(a_1)$ &$\calD P_1$&I&1\\ \hline
30& $E_8$&$\calD P_1$&IV&1\\ \hline

\end{tabular}
\end{center}
\label{}
\caption{Cyclic subgroups}
\end{table}

\clearpage
\begin{table}[t]
\begin{center}
{\small \begin{tabular}{| l |r| r | r |r|}
\hline
Order&Structure&Surface&Equation&Conjugacy classes\\ \hline
4&$2^2$&$\calD P_4$&&$\infty^2$\\ \hline
4&$2^2$&$\calD P_2$&XII&$\infty^5$\\ \hline
4&$2^2$&$\calD P_1$&XXI&$\infty^5$\\ \hline
4&$2^2$&$\calD P_1$&V, VI,X,XV,XVII&$\infty^3$\\ \hline
8&$2\times 4$&$\calD P_4$&\eqref{6.3}&$\infty^1$\\ \hline
8&$2\times 4$&$\calD P_2$&V&$2\times \infty^1$\\ \hline
8&$2\times 4$&$\calD P_2$&I-V&$\infty^1$\\ \hline
8&$2\times 4$&$\calD P_2$&VII&$\infty^2$\\ \hline
8&$2\times 4$&$\calD P_1$&VIII,XVI&$\infty^2$\\ \hline
8&$2^3$&$\calD P_4$&&$\infty^2$\\ \hline
8&$2^3$&$\calD P_2$&I-V,X&$\infty^3$\\ \hline
9&$3^2$&$\calD P_3$&I&1\\ \hline
9&$3^2$&$\calD P_3$&I,IV&$2\times \infty^1$\\ \hline
9&$3^2$&$\calD P_3$&III&$1$\\ \hline
9&$3^2$&$\calD P_1$&I,II,III&$\infty^1$\\ \hline
12&$2\times 6$&$\calD P_4$&\eqref{6.5}&1\\ \hline
12&$2\times 6$&$\calD P_2$&III, VIII&$\infty^1$\\ \hline
12&$2\times 6$&$\calD P_1$&II,VIII,XIII&$\infty^2$\\ \hline
12&$2\times 6$&$\calD P_1$&III,XII&$\infty^2$\\ \hline
12&$2\times 6$&$\calD P_1$&II,VII&$\infty^1$\\ \hline
16&$2^4$&$\calD P_4$&&$\infty^2$\\ \hline
16&$2^2\times 4$&$\calD P_2$&II,III,V&$\infty^1$\\ \hline
16&$4^2$&$\calD P_2$&II&1\\ \hline
16&$2\times 8$&$\calD P_2$&II&1\\ \hline
18&$3\times 6$&$\calD P_3$&I&1\\ \hline
18&$3\times 6$&$\calD P_1$& III&$\infty^1$\\ \hline
18&$3\times 6$&$\calD P_1$&II&$1$\\ \hline
24&$2\times 12$&$\calD P_1$&VIII&$\infty^1$\\ \hline
24&$2\times 12$&$\calD P_2$&III&$1$\\ \hline
27&$3^3$&$\calD P_3$&I&1\\ \hline
32&$2\times 4^2$&$\calD P_2$&II&1\\ \hline
36&$6^2$&$\calD P_1$&II&1\\ \hline

\end{tabular}}
\end{center}
\label{}
\caption{Abelian non-cyclic subgroups}
\end{table}

\clearpage
\begin{table}[h]
\begin{center}
\begin{tabular}{| l |r| r | r |r|}
\hline
Order&Structure&Surface&Equation&Conjugacy classes\\ \hline
6&$D_{6}$&$\calD P_3$&III,IV,VIII,XI&$\infty^2$\\ \hline
6&$D_{6}$&$\calD P_3$&I&$1$\\ \hline
8&$D_{8}$&$\calD P_4$&\eqref{6.4}&$\infty^1$\\ \hline
8&$D_{8}$&$\calD P_2$&II,III,V,VII&$\infty^2$\\ \hline
 8&$D_{8}$&$\calD P_1$&I, XVI&$\infty^3$\\ \hline
8&$Q_8$&$\calD P_1$&II,V,VI,XIV&$\infty^2$\\ \hline
 8&$3\times D_{8}$&$\calD P_1$&I, XVI&$\infty^3$\\ \hline
8&$3\times Q_8$&$\calD P_1$&II,III&$\infty^1$\\ \hline
12&$2\times D_6$&$\calD P_3$&I,VI&$\infty^1$\\ \hline
12&$T$&$\calD P_2$&II&1\\ \hline
12&$2\times D_6$&$\calD P_2$&I,II,IV,IX&$\infty^2$\\ \hline
12&$2\times D_6$&$\calD P_2$&II&$1$\\ \hline
12&$2\times D_6$&$\calD P_1$&I,II,III,VI&$\infty^1$\\ \hline
16&$2\times D_{8}$&$\calD P_2$&II,III,V,VII&$\infty^2$\\ \hline
16&$D_{16}$&$\calD P_1$&I,IX&$\infty^1$\\ \hline
18&$3\times D_6$&$\calD P_3$&III&1\\ \hline
18&$3\times D_6$&$\calD P_3$&I,IV&$2\times \infty^1$\\ \hline
18&$3\times D_6$&$\calD P_1$&I,II,III&$\infty^1$\\ \hline
24&$\overline{T}$&$\calD P_1$&I,V&$\infty^1$\\ \hline
24&$2\times T$&$\calD P_4$&\eqref{6.5}&1\\ \hline
24&$S_4$&$\calD P_3$&I&3\\ \hline
24&$S_4$&$\calD P_3$&II&1\\ \hline
24&$S_4$&$\calD P_2$&II&1\\ \hline
24&$S_4$&$\calD P_2$&I, II,IV&$\infty^1$\\ \hline
24&$2\times T$&$\calD P_2$&II&1\\ \hline
24&$3\times Q_8$&$\calD P_1$&I&$1$\\ \hline
24&$3\times Q_8$&$\calD P_1$&II&$1$\\ \hline
24&$3\times D_8$&$\calD P_1$&I&$\infty^2$\\ \hline
36&$6\times D_6$&$\calD P_1$&I,II,III&$\infty^1$\\ \hline
48&$2\times O$&$\calD P_2$&I,II,III,V&$\infty^1$\\ \hline
48&$2\times \overline{T}$&$\calD P_1$&I&$1$\\ \hline
72&$3\times \overline{T}$&$\calD P_1$&I&$1$\\ \hline

\end{tabular}
\end{center}
\label{}
\caption{Products of cyclic groups and polyhedral or binary polyhedral non-cyclic group}
\end{table}

\begin{table}[h]
\begin{center}
 \begin{tabular}{| l |r| r | r |r|}
\hline
Order&Structure&Surface&Equation&Conjugacy classes\\ \hline
$8$&$D_{8}$&$\calD P_4$&\eqref{6.3}&$\infty^1$\\ \hline
$8$&$D_{8}$&$\calD P_2$&II,III,IV,V,VII&$\infty^2$\\ \hline
$8$&$Q_8$&$\calD P_1$&I,V,XIV&$\infty^2$\\ \hline
$8$&$D_8$&$\calD P_1$&I,V,XVI&$\infty^2$\\ \hline
$16$&$L_{16}$&$\calD P_4$&\eqref{6.3}&$2\times \infty^1$\\ \hline
$16$&$2\times D_8$&$\calD P_2$&II,III,IV,V,VII&1\\ \hline
$16$&$2\times D_8$&$\calD P_2$&I&1\\ \hline
$16$&$AS_{16}$&$\calD P_2$&II,III,V&$2\times \infty^1$\\ \hline
$16$&$M_{16}$&$\calD P_2$&II&1\\ \hline
$16$&$D_{16}$&$\calD P_1$&I,IX&$\infty^1$\\ \hline
$32$&$2^2:8$&$\calD P_4$&\eqref{6.4}&1\\ \hline
$32$&$2^4:2$&$\calD P_4$&\eqref{6.3}&$\infty^1$\\ \hline
$32$&$D_8:4$&$\calD P_2$&II&1\\ \hline
$32$&$2\times AS_{16}$&$\calD P_2$&III,V&$3\times \infty^1$\\ \hline
$32$&$2\times M_{16}$&$\calD P_2$&II&1\\ \hline
$64$&$2^4:4$&$\calD P_4$&\eqref{6.4}&1\\ \hline
$64$&$2\times (D_8:4)$&$\calD P_2$&II&1\\ \hline
$27$&$H_3(3)$&$\calD P_3$&I,III,IV&$\infty^1$\\ \hline
$81$&$3^3:3$&$\calD P_3$&I&1\\ \hline

\end{tabular}
\end{center}
\label{}
\caption{Non-abelian $p$-groups}
\end{table}
\clearpage

\begin{table}[h]
\begin{center}
\begin{tabular}{| l |r| r | r |r|}
\hline
Order&Structure&Surface&Equation&Conjugacy classes\\ \hline
$18$&$3^2:2$&$\calD P_3$&I&2\\ \hline
$24$&$2D_{12}$&$\calD P_1$&II,VI&$\infty^1$\\ \hline
$24$&$D_8:3$&$\calD P_2$&III&1\\ \hline
$36$&$3^2:2^2$&$\calD P_3$&I&1\\ \hline
$42$&$2\times (7:3)$&$\calD P_1$&II&1\\ \hline
$48$&$2^4:3$&$\calD P_4$&\eqref{6.5}&1\\ \hline
$48$&$2\times D_8:3$&$\calD P_2$&III&1\\ \hline
$48$&$\overline{T}:2$&$\calD P_1$&I&1\\ \hline
$48$&$4^2:3$&$\calD P_2$&II&1\\ \hline
$54$&$H_3(3):2$&$\calD P_3$&IV&$\infty^1$\\ \hline
$54$&$3^3:2$&$\calD P_3$&I&2\\ \hline
$60$&$A_5$&$\calD P_5$&&1\\ \hline
$72$&$3\times 2D_{12}$&$\calD P_1$&II&1\\ \hline
$80$&$2^4:5$&$\calD P_4$&\eqref{6.6}&1\\ \hline
$96$&$2^4:S_3$&$\calD P_4$&\eqref{6.5}&1\\ \hline
$96$&$4^2:S_3$&$\calD P_2$&I&2\\ \hline
$96$&$2\times (4^2:3)$&$\calD P_2$&II&1\\ \hline
$108$&$3^3:4$&$\calD P_3$&I&2\\ \hline
$108$&$3^3:2^2$&$\calD P_3$&I&2\\ \hline
$120$&$S_5$&$\calD P_5$&&1\\ \hline
$120$&$S_5$&$\calD P_3$&II&1\\ \hline
$144$&$3\times (\overline{T}:2)$&$\calD P_1$&I&1\\ \hline
$160$&$2^4:D_{10}$&$\calD P_4$&\eqref{6.6}&1\\ \hline
$162$&$3^3:S_3$&$\calD P_3$&I&2\\ \hline
$168$&$L_2(7)$&$\calD P_2$&I&1\\ \hline
$192$&$2\times (4^2:S_3)$&$\calD P_3$&I&1\\ \hline
$216$&$3^3:D_8$&$\calD P_3$&I&2\\ \hline
$336$&$2\times L_2(7)$&$\calD P_2$&I&1\\ \hline
$648$&$3^3:S_4$&$\calD P_3$&I&2\\ \hline

\end{tabular}
\end{center}
\label{}
\caption{Other groups}
\end{table}

\bibliographystyle{amsplain}

\begin{thebibliography}{100}

\bibitem{Adem}  A.  Adem, J.  Milgram, \textit{Cohomology of finite groups}, Grundlehren der Mathematischen Wissenschaften,  {\bf 309},  Springer-Verlag, Berlin, 1994.

\bibitem{Alberich} M. Alberich-Carrami–ana, \emph{Geometry of the plane Cremona maps}, Lecture Notes in Mathematics, 1769. Springer-Verlag, Berlin, 2002. xvi+257 pp.

\bibitem{Autonne} L. Autonne, \textit{
Recherches sur les groupes d'ordre fini contenus dans le groupe Cremona, Prem. M\'em. G\'en\'eralit\'es et groupes quadratiques},
J. Math. Pures et Appl., (4) {\bf 1} (1885), 431--454.

\bibitem{BT} S. Bannai, H. Tokunaga, \textit{A note on embeddings of $S_4$ and $A_5$ into the Cremona group and versal Galois covers}, Publ. Res. Inst. Math. Sci. {\bf 43} (2007),  1111--1123.

\bibitem {Bayle-Beauville} L. Bayle, A.  Beauville,
\textit{Birational involutions of $P^2$},
Asian J. Math. {\bf 4} (2000),  11--17.


\bibitem{Blanc} J. Blanc, \textit{Finite abelian subgroups of the Cremona group of the plane}, Thesis, Univ. of Geneva, 2006.

\bibitem{Bl} J. Blanc, \textit{Elements and cyclic subgroups of finite order of the Cremona group}, arXiv:0809.4673.\footnote{This reference is absent in the published version of the paper}

\bibitem{BB} A. Beauville, J. Blanc, \textit{On Cremona transformations of prime order}, C. R. Math. Acad. Sci. Paris {\bf 339} (2004),  257--259.

\bibitem{Beauville} A. Beauville, \textit{p-elementary subgroups of the Cremona group}, J. of Algebra {\bf 314} (2007), 553--564.


\bibitem{Bertini} E. Bertini, \textit{
Ricerche sulle trasformazioni univoche involutorie nel piano},
Annali di Mat. Pura Appl. (2) {\bf 8} (1877), 254--287.

\bibitem {Blichfeldt} H. Blichfeldt, \textit{Finite collineation groups, with an introduction to the theory of operators and substitution groups}, Univ. of Chicago Press, Chicago, 1917.

\bibitem{Bottari} A. Bottari, \textit{Sulla razionalit\`a dei piani multipli $\{x,y,\sqrt[n]{F(x,y)}\}$}, Annali di Mat. Pura Appl. (3) {\bf 2} (1899), 277--296.


\bibitem{Calabri1} A. Calabri,
\textit{Sulle razionalit\`a dei piani doppi e tripli cyclici}, Ph. D. thesis, Univ. di Roma ``La Sapienza'', 1999.

\bibitem{Calabri2} A. Calabri,
\textit{On rational and ruled double planes},
Annali. di Mat. Pura Appl. (4) {\bf 181} (2002),  365--387.

\bibitem{Carter} R. Carter, \textit{Conjugacy classes in the Weyl group}. in `` Seminar on Algebraic Groups and Related Finite Groups'', The Institute for Advanced Study, Princeton, N.J., 1968/69, pp. 297--318, Springer, Berlin.


\bibitem{Castelnuovo} G. Castelnuovo,  \textit{Sulle razionalit\`a delle involutioni piani}, Math. Ann, {\bf 44} (1894), 125--155.

\bibitem{CE} G. Castelnuovo, F. Enriques, \textit{Sulle condizioni di razionalit\`a dei piani doppia}, Rend. Circ. Mat. di Palermo, {\bf 14} (1900), 290--302.

\bibitem  {Coble}  A. Coble, \textit{Algebraic geometry and theta functions} (reprint of the 1929 edition), A. M. S. Coll. Publ., v. 10. A. M. S., Providence, R.I., 1982. 

\bibitem{Atlas} J. Conway, R.  Curtis, S. Norton,  R.  Parker,  R. Wilson, \textit{Atlas of finite groups}, Oxford Univ.  Press, Eynsham, 1985.

\bibitem{Conway2} J. Conway, D. Smith,
\textit{On quaternions and octonions: their geometry, arithmetic, and symmetry},  A K Peters, Ltd., Natick, MA, 2003.

\bibitem {Coolidge} J. Coolidge, \textit{A treatise on algebraic plane curves}, Dover Publ. New York. 1959.

\bibitem{Corti} A. Corti,
\textit{Factoring birational maps of threefolds after Sarkisov},
J. Algebraic Geom. {\bf 4} (1995), 223--254.






\bibitem{deFernex} T. de Fernex, \textit{On planar Cremona maps of prime order}, Nagoya Math. J. {\bf 174} (2004), 1--28.

\bibitem{dFe-Ein} T. de Fernex, L. Ein, \textit{Resolution of indeterminacy of pairs}, in ``Algebraic geometry'', pp. 165--177, de Gruyter, Berlin, 2002.

\bibitem{Demazure} M. Demazure, \textit{Surfaces de Del Pezzo, I-V}, in  ``S\'eminaire sur les Singularit\'es des Surfaces'',  Ed.  by M. Demazure, H. Pinkham and B. Teissier. Lecture Notes in Mathematics, 777.  Springer, Berlin, 1980, pp. 21--69.

\bibitem{Dolgachev} I. Dolgachev, \textit{Weyl groups and Cremona transformations}. in ``Singularities, Part 1 (Arcata, Calif., 1981)'', 283--294, Proc. Sympos. Pure Math., 40, Amer. Math. Soc., Providence, RI, 1983.

\bibitem{Topics} I. Dolgachev, \emph{Topics in classical algebraic geometry, Part I}, manuscript in preparations, see  www.math.lsa.umich.edu/~idolga/lecturenotes.html.


\bibitem{DuVal}  P. Du Val, \emph{On the Kantor group of a set of points in a plane}, Proc. London Math. Soc. {\bf 42} (1936), 18--51.

\bibitem{Enriques} F. Enriques, \textit{Sulle irrazionalita da cui puo farsi dipendere la resoluzione d'un' equazione algebrica $f(x,y,z) = 0$ con funzioni razionali di due parametri}, Math. Ann. {\bf 49} (1897), 1--23.

\bibitem {Geiser} C. Geiser,
\textit{\"Uber zwei geometrische Probleme}, J. Reine Angew. Math. {\bf 67} (1867), 78--89.




\bibitem {Gorenstein} D. Gorenstein, \textit{Finite groups},  Chelsea Publ. Co., New York, 1980.

\bibitem{Goursat} \'E. Goursat, \textit{Sur les substitutions orthogonales et les divisions rŽgulires de l'espace},
Ann. de l'\'Ecole Norm. Sup. (3) {\bf 6} (1889),  9--102.

\bibitem {Hartshorne}  R. Hartshorne, \textit{Algebraic geometry}, Graduate Texts in Mathematics, No. 52. Springer-Verlag, New York-Heidelberg, 1977.

\bibitem{Hosoh1} T. Hosoh, \textit{Automorphism groups of cubic surfaces}, J. Algebra {\bf 192} (1997), 651--677.


\bibitem{Hudson} H. Hudson, \textit{Cremona transformations in plane and space}, Cambridge Univ. Press. 1927.
\bibitem{Isk1} V. A. Iskovskih, \textit{Rational surfaces with a pencil of rational curves} (Russian) Mat. Sb. (N.S.) {\bf 74} (1967),  608--638.

\bibitem{Isknew} V. A. Iskovskikh, \textit{Rational surfaces with a pencil of rational curves with positive square of the canonical class}, Math. USSR Sbornik, {\bf 12} (1970), 93--117.

\bibitem{Isk2}  V. A. Iskovskih, \textit{Minimal models of rational surfaces over arbitrary fields} (Russian) Izv. Akad. Nauk SSSR Ser. Mat. {\bf 43} (1979), 19--43,

\bibitem{Isk3} V. A. Iskovskikh, \textit{Factorization of birational mappings of rational surfaces from the point of view of Mori theory} (Russian) Uspekhi Mat. Nauk {\bf 51} (1996), 3--72; translation in Russian Math. Surveys {\bf 51} (1996),  585--652.

\bibitem{Isk4} V. A. Iskovskikh, \textit{Two nonconjugate embeddings of the group $S_3\times Z_2$ into the Cremona group} (Russian) Tr. Mat. Inst. Steklova {\bf 241} (2003), Teor. Chisel, Algebra i Algebr. Geom., 105--109; translation in Proc. Steklov Inst. Math. {\bf 241} (2003), 93--97.

\bibitem{Isk5} V. A. Iskovskikh, \textit{Two non-conjugate embeddings of $S_3\times \bbZ_2$ into the Cremona group II}, Adv. Study in Pure Math. {\bf 50} (2008), 251--267.

\bibitem{Kantor} S.  Kantor, \textit{Theorie der endlichen Gruppen von eindeutigen Transformationen in der Ebene},
Berlin. Mayer \& MŸller. 111 S. gr. $8^\circ$. 1895.


\bibitem{Kollar-Mori} J. Koll\'ar, S.  Mori,
\textit{Birational geometry of algebraic varieties},
Cambridge Tracts in Mathematics, 134.
Cambridge University Press, Cambridge, 1998.


\bibitem{Lemire}  N. Lemire, V. Popov, Z. Reichstein, \textit{Cayley groups}, 
J.  Amer. Math. Soc. {\bf 19} (2006), 921--967.

\bibitem{Manin} Yu. I. Manin, \textit{Rational surfaces over perfect fields. II} (Russian) Mat. Sb. (N.S.) {\bf 72} (1967), 161--192.

\bibitem{Manin2} Yu. I. Manin, \textit{
Cubic forms: algebra, geometry, arithmetic}
Translated from Russian by M. Hazewinkel. North-Holland Mathematical Library, Vol. 4.
North-Holland Publishing Co., Amsterdam-London; American Elsevier Publishing Co., New York, 1974.

\bibitem{Noether} M. Noether, \textit{\"Uber die ein-zweideutigen Ebenentransformationen}, Sitzungberichte der physic-medizin, Soc. zu Erlangen, 1878.

\bibitem{Segre} B. Segre, \textit{The non-singular cubic surface}, Oxford Univ. Press. Oxford. 1942.

\bibitem{Springer} T. Springer,  \textit{Invariant theory}.  Lecture Notes in Mathematics, Vol. 585. Springer-Verlag, Berlin-New York, 1977.

\bibitem{Wiman} A. Wiman, \textit{Zur Theorie endlichen Gruppen von birationalen Transformationen in der Ebene}, Math. Ann. {\bf 48} (1896), 195--240.

 \bibitem{Zas} H. Zassenhaus, \textit{The theory of groups}, Chelsea Publ. New York, 1949.

\bibitem{Zhang} D.-Q. Zhang, \textit{Automorphisms of finite order on rational surfaces. With an appendix by I. Dolgachev},  J. Algebra {\bf 238} (2001),  560--589.
\end{thebibliography}

\end{document}